\documentclass[a4paper,11pt]{article}
\usepackage{a4wide}
\usepackage{theorem}
\usepackage{amsmath}
\usepackage{array}
\usepackage{amssymb}
\usepackage{amsfonts}
\usepackage[english]{babel}
\usepackage{epsf}
\usepackage{epsfig}
\usepackage{graphicx}

\newtheorem{theo}{\indent Theorem \newline}[section]
\newtheorem{defi}[theo]{\indent Definition\newline}
\newtheorem{rem}[theo]{\noindent Remark}
\newtheorem{prop}[theo]{\indent Proposition\newline}
\newtheorem{lemma}[theo]{\indent Lemma\newline}
\newtheorem{cor}[theo]{\indent Corollary \newline}

 \def\N{{\mathbb{N}}}
\def\Z{{\mathbb{Z}}}
\def\Q{{\mathbb{Q}}}
\def\R{{\mathbb{R}}}
\def\C{{\mathbb{C}}}
\def\H{{\mathbb{H}}}
\def\A{{\mathbb{A}}}

\usepackage{eufrak}
\newcommand{\goth}[1]{\EuFrak{#1}}

\newcommand{\p}{\mathop{\goth p}\nolimits}

\newcommand{\ind}{\mathop{\rm ind}\nolimits}

\newcommand{\coker}{\mathop{\rm coker}\nolimits}

\newlength{\indentation}%
\makeatletter
\newcommand\@makefntextsans[1]{%
    \parindent 0em%
    \noindent%
    \hb@xt@0em{\hss}%
    #1}
\def\footnotetextsans{%
     \@ifnextchar [\@xfootnotenextsans%
       {\@footnotetextsans}}
\def\@xfootnotenextsans[#1]{%
  \begingroup%
     \csname c@\@mpfn\endcsname #1\relax%
  \endgroup%
  \@footnotetextsans}
\long\def\@footnotetextsans#1{\insert\footins{%
    \reset@font\footnotesize%
    \interlinepenalty\interfootnotelinepenalty%
    \splittopskip\footnotesep%
    \splitmaxdepth \dp\strutbox \floatingpenalty \@MM%
    \hsize\columnwidth \@parboxrestore%
    \color@begingroup%
      \@makefntextsans{%
        \rule\z@\footnotesep\ignorespaces#1\@finalstrut\strutbox}
    \color@endgroup}}
\makeatother

\begin{document}

\cleardoublepage

\title{Open strings, Lagrangian conductors and Floer functor}
\author{Jean-Yves Welschinger}
\maketitle

\makeatletter\renewcommand{\@makefnmark}{}\makeatother
\footnotetextsans{Keywords: Open strings, Lagrangian submanifolds, Floer cohomology, homological mirror symmetry.}
\footnotetextsans{AMS Classification : 53D40.
}

{\bf Abstract:}

 We introduce a contravariant functor, called Floer functor, from the category 
 of Lagrangian conductors of a symplectic manifold to the homotopy category
  of bounded chain complexes of open strings in this manifold. The latter two categories are 
 defined for all symplectic manifolds, whereas Floer functor is defined for semipositive manifolds which are either closed
 or convex at infinity. We then prove that when the first Chern class of the symplectic manifold vanishes, Lagrangian spheres
 define Lagrangian conductors so that in particular their integral Floer cohomology is well defined. This requires the
 introduction of singular almost-complex structures given by symplectic field theory. 

\section*{Introduction}

The present paper deals with Lagrangian Floer theory, with a view towards problems of mirror symmetry
arising from string theory physics. We introduce a notion of strings living in a symplectic manifold $(X, \omega)$
and focus our attention to open strings in dimension greater than two. These open strings are the objects of
a small preadditive category, denoted by ${\cal OS} (X, \omega)$, which has a duality, a tensor product defined on 
objects and morphisms as well as an action of the integers. This structure is only partially studied here, we postpone
a detailed study of duality and of closed strings. Nevertheless, we introduce a functor ${\cal C}$,
called functor coefficients, from the category of open strings to the category of free modules over
the Novikov ring $\Z ((t^\R))$.
This functor actually satisfies the axioms of a $1+1$ topological field theory. We also study the homotopy category
of bounded chain complexes of ${\cal OS} (X, \omega)$, or rather of its quotient ${\cal OS}_1 (X, \omega)$ by the $\Z$-action,
which we denote by $K^b ({\cal OS}_1 (X, \omega))$. Then, we turn attention to closed Lagrangian submanifolds. 
We introduce a notion
of Lagrangian conductors which form the objects of a small category ${\cal CL}^\pm (X, \omega)$ associated to any symplectic
manifold $(X, \omega)$. This category also carry more structure, exact sequences make sense in this category, objects
have subobjects as well as extensions. Roughly speaking, a Lagrangian conductor is a finite collection of closed Lagrangian
submanifolds of $(X, \omega)$ which are transversal to each other together with a generic function defined on the universal
curve over  Stasheff's associahedron of the appropriate dimension with values in the space ${\cal J}_\omega$ of compatible almost complex
structures of $(X, \omega)$. The latter function has to be compatible with the structure of the associahedron.  
Lagrangian Floer theory \cite{Floer}, \cite{Fuk2}, then provides a contravariant functor ${\cal F} : {\cal CL}^\pm (X, \omega) \to
K^b ({\cal OS}_1 (X, \omega))$, called Floer functor. At this point, the symplectic manifold  $(X, \omega)$ is supposed to be
semipositive and either closed or convex at infinity. Twisted complexes introduced by Kontsevich in \cite{Kont} appear as
augmentations of the image chain complexes of ${\cal F}$, so that the derived category $D^b ({\cal F} (X, \omega))$ introduced
in \cite{Kont} fits well in the present formalism. We finally restrict our attention to symplectic manifolds for which the first Chern
class vanish and show that the category of Lagrangian conductors then contains Lagrangian spheres. For this purpose, 
we extend the space of almost complex structures by introducing singular ones and get a space 
$\overline{\cal J}_\omega = {\cal J}_\omega \sqcup \partial {\cal J}_\omega$. These singular almost complex structures are
actually the split almost complex structures given by symplectic field theory. The ones associated to Lagrangian spheres are called
$A_1$-singular. The $A_1$-singular almost complex structures $J$ for which there exist $J$-holomorphic disks with boundary
on the Lagrangian sphere turn out to be localized on strata of greater codimensions than expected. This localization phenomenon
(Theorem \ref{theolocalization}) makes it possible to define integral Floer cohomology for such spheres, to include them
in the present  formalism and hence to overcome in this case the problem pointed out in the note added in proof of \cite{Kont}. 
A systematic study of the latter obstruction is found in \cite{FOOO}. 

The first part of the paper is devoted to open strings and the functor coefficients ${\cal C}$, the second part 
to Lagrangian conductors and singular almost complex structures, the third part to Floer functor and the last one
to manifolds with vanishing first Chern class and the localization Theorem. \\

{\bf Acknowledgements: }

I am grateful to the French Agence nationale de la recherche for its support. I wish also to acknowledge
Paul Biran and Dietmar Salamon for fruitful discussions on the last part of the paper, ETH Z\"urich for its
hospitality during a one week stay at the end of $2007$, the London mathematical society for giving
me the opportunity to discuss the last part of the paper at several places in England and finally
Denis Auroux, Paul Seidel and Ivan Smith for discussions on the homological mirror symmetry conjecture. 

\tableofcontents

\section{Open strings and functors coefficients}

In this first paragraph, $(X, \omega)$ denotes a symplectic manifold of dimension $2n \geq 4$ ;
we won't consider the special case of surfaces. Denote by $\pi_{\cal L} : {\cal L} \to X$ the bundle
of Lagrangian subspaces of $TX$, so that its fiber over every point $x \in X$ is the grassmannian
${\cal L}_x = \pi_{\cal L}^{-1} (x) = \{ l_x \in Vect(T_x X) \, \vert \, \text{dim}_\R (l_x) = n \text{ and } \omega \vert_{l_x} = 0 \}$.
Denote by $\widetilde{GL}_n^+ (\R)$ and $\widetilde{GL}_n^- (\R)$ the two Lie group structures on
the nowhere trivial double cover of $GL_n (\R)$ which turn
the covering map into a morphism. Any lift of a reflection is of order two in $\widetilde{GL}_n^+ (\R)$
and four in $\widetilde{GL}_n^- (\R)$, see \cite{ABS}. 

\begin{defi}
A $\widetilde{GL}_n^\pm (\R)$-structure on $l \in {\cal L}$ is a $\widetilde{GL}_n^\pm (\R)$-principal bundle
$\p^\pm$ which lifts the $GL_n (\R)$-principal bundle of frames of $l$.
\end{defi}

Denote by ${\cal L}^\pm = \{ (l , \p^\pm) \, \vert \, l \in {\cal L} \text{ and } \p^\pm \text{ is a } \widetilde{GL}_n^\pm (\R)\text{-structure on } l \}$.
We now have to choose once for all between ${\cal L}^+$ and ${\cal L}^-$. Since we do not mind which one to choose,
we leave this choice to the reader and denote by ${\cal L}^\pm$ the chosen bundle. It is however understood that
${\cal L}^\pm$ should denote throughout the paper either ${\cal L}^+$ or  ${\cal L}^-$ but not both.

The first two subparagraphs are devoted to the construction of the category of open strings ${\cal OS} (X, \omega)$. 
The third one is devoted to the construction of the functor coefficients ${\cal C}$ and the last one to the
homotopy category of bounded chain complexes of open strings and the notion of augmentations of these complexes. 

\subsection{Open strings}
\label{subsectionopenstring}

\begin{defi}
An elementary open string of $(X, \omega)$ is a homotopy class with fixed extremities transversal to each other 
of paths $\lambda : [-1 , 1] \to {\cal L}^\pm$ such that $\pi_{\cal L} \circ \lambda$ is a constant path in $X$. 
\end{defi}
An elementary closed string of $(X, \omega)$ is likewise a homotopy class of loops 
$\lambda : S^1 \to {\cal L}^\pm$ such that $\pi_{\cal L} \circ \lambda$ is a constant loop
in $X$. However, we
will focus attention on open strings throughout this paper and postpone discussion on
closed string.

\begin{defi}
The index of an elementary open string $\lambda : t \in  [-1 , 1] \mapsto (l (t) , \p^\pm (t)) \in {\cal L}^\pm$
is the quantity $\mu (\lambda) = \frac{n}{2} - \mu(l_{-1} , l) \in \Z$, where $l_{-1}  : [-1 , 1] \to {\cal L}$
denotes the constant path $l(-1)$ and $ \mu(l_{-1} , l)$ the relative Maslov index defined
in \cite{RobSal} (see also \cite{Vit}, \cite{Kont}, \cite{Seid}, \cite{Fuk2}).
\end{defi}

{\bf Examples:}

1) If $n=1$ and $l : t \in  [-1 , 1] \mapsto \exp{(i \frac{\pi}{4}t)} \R \subset \C$, then 
$\mu(l_{-1} , l) = -\frac{1}{2}$ and $\mu (\lambda) = 1$.\\

2) As a consequence, if $(X, \omega)$  is the cotangent bundle of some manifold, $L \subset X$
the zero section, $L_f \subset X$ the graph of the differential of a Morse function $f : L \to \R$,
$x$ a critical point of $f$ with Morse index $i_M (x)$ and $l : t \in  [-1 , 1] \mapsto L_{(t+1)f} \subset X$,
then $\mu (\lambda) = n - i_M (x)$.\\

Let $\lambda'$ and $\lambda''$ be two elementary open strings such that $\lambda' (1) = \lambda'' (-1)$, we
denote by $\lambda' * \lambda''$ the concatenation 
$t \in  [-1 , 1] \mapsto \left\{ \begin{array}{l}
\lambda' (2t + 1) \text{ if } -1 \leq t \leq 0 \\
\lambda'' (2t - 1) \text{ if } 0 \leq t \leq 1 \\
\end{array}
\right.
$.
Likewise, if  $\lambda$ is an elementary string and $\tau \in  [-1 , 1]$, we set $\lambda_{\tau}' :  t \in  [-1 , 1] \mapsto 
\lambda ( (\frac{1 + \tau}{2})t + (\frac{\tau - 1}{2}))$ and $\lambda_{\tau}'' :  t \in  [-1 , 1] \mapsto 
\lambda ( (\frac{1 - \tau}{2})t + (\frac{1 + \tau}{2}))$. When the index of $\lambda_{\tau}''$ is positive (resp. negative),
we say that the string $\lambda_{\tau}'$ stretches to $\lambda$ (resp. shrinks to $\lambda$). Note however that the
index is not additive under concatenation and does not increase in general when the string stretches. It jumps
only when the string no more has transversal extremities. 

For every elementary open string $\lambda$, let the configuration space $\text{Conf} (\lambda) $ of $\lambda$ be the space of
elementary open strings $ \lambda' $ such that $ \lambda' (-1) = \lambda (-1)$, $\lambda' (1) = \lambda (1)$
and $\mu(\lambda') = \mu (\lambda')$.

\begin{lemma}
\label{lemmaconfig}
Configuration spaces of elementary open strings are connected. They are
moreover simply connected as soon as $n > 2$.
\end{lemma}

{\bf Proof:}

Let $\lambda_x = (l_x , \p^\pm)$ be an elementary open string at the point $x \in X$. The configuration space
$\text{Conf} (\lambda_x)$ is a double cover of the configuration space $\text{Conf} (l_x) = 
\{ l'  :  [-1 , 1] \to {\cal L} \text{ such that } l' (-1) = l_x (-1), l' (1) = l_x (1)
\text{ and } \mu(l') = \mu (l_x) \}  $. The connectedness of $\text{Conf} (l_x)$ is well known, see \cite{RobSal}
for example, we have to prove that the double cover is not trivial. The fundamental group of
$\text{Conf} (l_x)$ is isomorphic to $\pi_2( GL_n (\C)/ GL_n (\R)) \cong \pi_1 (GL_n (\R))$, hence to
$\Z/2\Z$ when $n>2$ and to $\Z$ when $n = 2$. Let $u : S^2 \to {\cal L}_x \cong GL_n (\C)/ GL_n (\R)$
be a map generating $\pi_2( {\cal L}_x)$ and $u^*{\cal L}_x \to S^2$ be the associated tautological bundle. 
Then, the second Stiefel-Whitney class of $u^*{\cal L}_x$ does not vanish in $H^2 (S^2 ; \Z/2\Z)$ since the
trivializations of $u^*{\cal L}_x$ over the two hemispheres of $S^2$ differ from a generator of $\pi_1 (GL_n (\R))$.
It follows that $u^*{\cal L}_x$ does not carry any $\widetilde{GL}_n^\pm (\R)$-structure, see \cite{KT},
so that the double cover $\text{Conf} (\lambda_x) \to \text{Conf} (l_x) $ is non-trivial (compare \S $3.6$ of \cite{OS}). $\square$\\

It follows from Lemma \ref{lemmaconfig} that $\text{Conf} (\lambda)$ is indeed the space of positions that
can take the string $\lambda$ when it vibrates. If $\lambda$ is an elementary open string and $e \in \Z$, we
denote by $\lambda^{+e}$ the elementary open string having same extremities but such that
$\mu (\lambda^{+e}) = \mu (\lambda) + e$. This defines an action of the group of integers on the set of elementary
strings by concatenation of closed strings. Let $\lambda_1 , \dots , \lambda_q$ be elementary strings of
$(X, \omega)$, $q \geq 1$, we denote by $\lambda_1 \otimes \dots \otimes \lambda_q$ their ordered union
modulo the relation $\lambda_1^{+e_1} \otimes \dots \otimes \lambda_q^{+e_q} = \lambda_1 \otimes \dots \otimes \lambda_q$
if and only if $e_1 + \dots +e_q = 0 \in \Z$. We set $\mu (\lambda_1 \otimes \dots \otimes \lambda_q) = \mu (\lambda_1) + \dots + \mu (\lambda_q)$.

\begin{defi}
\label{defopenstring}
An open string is either the empty set $\emptyset$ or an element $\lambda_1 \otimes \dots \otimes \lambda_q$
where $q \geq 1$ and $\lambda_1 , \dots , \lambda_q$ are elementary open strings. The integer $q \geq 1$ 
is called the cardinality of the string. The cardinality of $\emptyset$ vanishes.
\end{defi}

The set of open strings of $(X, \omega)$ will be denoted by $\text{Ob} ({\cal OS} (X, \omega))$, it inherits an action of
the integers from the one defined on elementary open strings. 

\begin{defi}
\label{defdualstring}
The dual string of an elementary string $\lambda$ is the string $\lambda^* :  t \in  [-1 , 1] \mapsto \lambda (-t)  \in {\cal L}^\pm$.
The dual string of a string $\lambda_1 \otimes \dots \otimes \lambda_q$ is the string $\lambda_q^* \otimes \dots \otimes \lambda_1^*$.
\end{defi}

The indices of an elementary open string and its dual are related by the formula $\mu ( \lambda) + \mu ( \lambda^*) = n$,
see  \cite{RobSal}, \cite{Kont}, \cite{Seid}, \cite{Fuk2}. In particular, $\mu ( \emptyset^*) = n$ ; we could decide that
$\emptyset^*$ is also an open string but this does not seem necessary.  We set
$\text{Ob} ({\cal OS}^* (X, \omega)) = \{ \emptyset^* \} \cup (\text{Ob} ({\cal OS} (X, \omega)) \setminus \{ \emptyset \})$.

\subsection{Propagation of open strings}

\begin{figure}[ht]
\begin{center}
\includegraphics[scale=0.5]{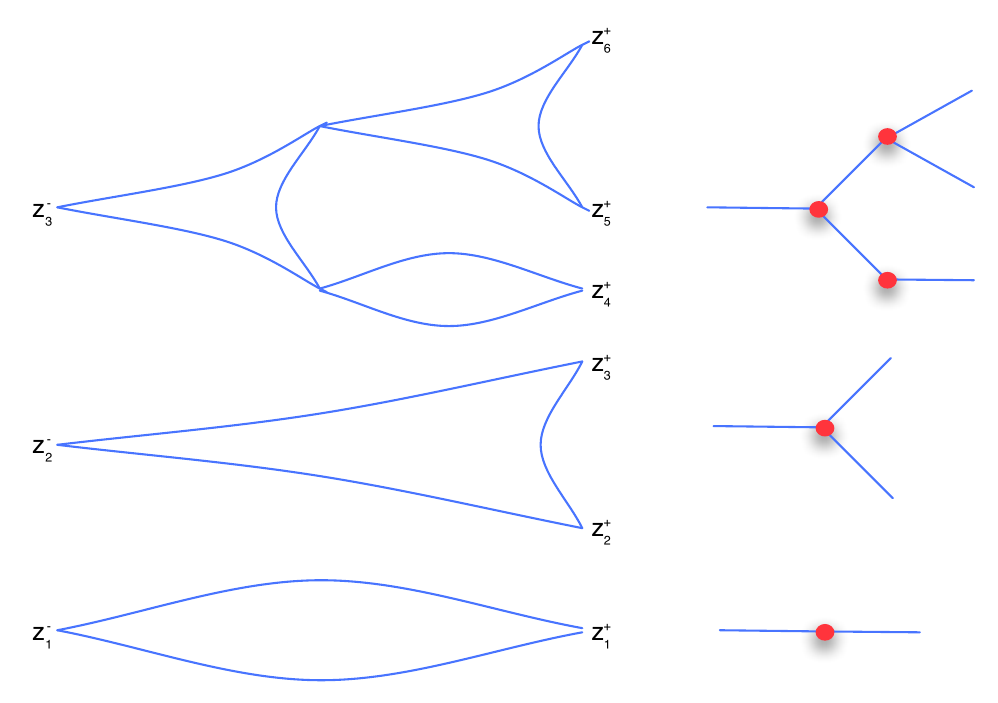}
\end{center}
\caption{A nodal disk $D_{\underline{z}}$ with $q^- = 3$, $q^+ = 6$ and its forest.}
\label{figuretrajectory}
\end{figure}

Let $\lambda_1^- \otimes \dots \otimes \lambda_{q^-}^-$ and $\lambda_1^+ \otimes \dots \otimes \lambda_{q^+}^+$
be two open strings of cardinalities $q^- \geq 1$ and $q^+ \geq 0$ such that either $q^- \leq q^+$ or $q^+ = 0$.
Let $(u, D_{\underline{z}} , \lambda_{\partial u} , \overline{\partial}_u )$ be a quadruple such that:

1) $D_{\underline{z}}$ is a punctured nodal disk without spherical component having $q^- $ possibly
reducible connected components. Its set of punctures $\underline{z} =\{ z_1^- , \dots , z_{q^-}^- , z_1^+ , \dots , z_{q^+}^+ \}$
is made of one negative puncture on every connected component of $D_{\underline{z}}$. Hence, $D_{\underline{z}}$
can be encoded by a forest $F(D_{\underline{z}})$ made of $q^-$ finite connected trees with one free negative edge on every connected component
and $q^+$ free positive edges, such that each vertex represents a disk and each edge adjacent to a vertex a puncture
on this disk, see Figure \ref{figuretrajectory}.

2) $u : D_{\underline{z}} \to X$ is a map which pulls back $\omega$ to a non-negative two-form on 
$D_{\underline{z}}$. It converges
to $\pi_{\cal L} \circ \lambda_i^\pm$ at every puncture $z_i^\pm$ and to a point $u(e) \in X$ at every 
pairs of punctures encoded by a non-free edge $e$ of $F(D_{\underline{z}})$.

3) $\lambda_{\partial u} : \partial D_{\underline{z}} \to u^* {\cal L}^\pm$ is a section which extends to a section
$\lambda_{u} : D_{\underline{z}} \to u^* {\cal L}^\pm$ over the whole $D_{\underline{z}}$ such that its limit at 
the puncture $z_i^\pm$ is the string $\lambda_i^\pm$ and its limit at every pairs of punctures encoded 
by a non-free edge $e$ of $F(D_{\underline{z}})$ is an elementary open string $\lambda (e)$.

4) $\overline{\partial}_u : L^{k,p} ( D_{\underline{z}} ; u^* TX ,  \lambda_{\partial u}) \to L^{k-1,p} ( D_{\underline{z}} ; 
\Lambda^{0,1} D_{\underline{z}} \otimes u^* TX )$ is an oriented Cauchy-Riemann operator. The space of Cauchy-Riemann
operators is contractible, because the space of complex structures of $D_{\underline{z}}$ compatible with its orientation
is contractible, the space of complex structures of $TX $ compatible with $\omega$ is contractible and for every such
structures, the space of associated Cauchy-Riemann operators is an affine space. However, there are two different 
orientations of its determinant line $\text{Det} (\overline{\partial}_u) =  \Lambda^{\text{max}} \ker (\overline{\partial}_u)
\otimes \Lambda^{\text{max}} \coker (\overline{\partial}_u)^*$, see the appendix of \cite{FloerHofer}. Note that here
$L^{k,p} ( D_{\underline{z}} ; u^* TX ,  \lambda_{\partial u})$ denotes the Banach space of sections of $u^* TX$
having $k$ derivatives in $L^p$ and with boundary values in the image of $ \lambda_{\partial u}$ whereas
$L^{k-1,p} ( D_{\underline{z}} ; 
\Lambda^{0,1} D_{\underline{z}} \otimes u^* TX )$ denotes the Banach space of complex antilinear one-forms of  $D_{\underline{z}}$
with value in $u^* TX$ having $k-1$ derivatives in $L^p$, $2 < p < +\infty$, $k \gg 0$. The Lebesgue measure fixed on $D_{\underline{z}}$ has infinite volume, it is near each puncture
the Lebesgue measure of $\R^+ \times [-1 , 1]$ read in a local chart, compare Definition \ref{defstriplikeend}. 
Finally, a Cauchy-Riemann operator here is an operator
$\overline{\partial}_u$ such that $\overline{\partial}_u (fv) = f \overline{\partial}_u (v) + \overline{\partial} f \otimes v$
for every real valued function $f$ on $D_{\underline{z}}$ and every section $v$. Such an operator is often called
generalized Cauchy-Riemann operator since we do not require it to be complex linear. 

\begin{defi}
\label{defelemtrajectory}
An elementary trajectory from the open string $\lambda_1^+ \otimes \dots \otimes \lambda_{q^+}^+$ to the open 
string $\lambda_1^- \otimes \dots \otimes \lambda_{q^-}^-$ is a homotopy class among stable maps having
fixed energy $a(\gamma) = \int_{D_{\underline{z}}} u^*\omega$ of such quadruples 
$\gamma = (u, D_{\underline{z}} , \lambda_{\partial u} , \overline{\partial}_u )$.
\end{defi}

Every open string has an empty trajectory going to itself, it has vanishing energy. An open string of positive cardinality cannot
propagate to an open string of greater cardinality. In particular, no trajectory has the empty string as a target except the empty
trajectory. One may think of the homotopy class appearing in Definition \ref{defelemtrajectory} as the space of configurations
of a given trajectory. The parameter $h (\gamma) = \exp(- a(\gamma))$ is less than one except for the empty trajectory or
classical trajectories, that is paths in $X$ between based points of elementary strings
or more generally maps in $X$ of forests $F(D_{\underline{z}})$, compare \cite{FukOh}, \cite{CorLal} and Remark \ref{remquantum}.

\begin{defi}
\label{defindextrajectory}
The index of an elementary trajectory $\gamma$ from the open string $\lambda_1^+ \otimes \dots \otimes \lambda_{q^+}^+$ to the open 
string $\lambda_1^- \otimes \dots \otimes \lambda_{q^-}^-$ is the difference
$\mu (\gamma) =  \mu (\lambda_1^- \otimes \dots \otimes \lambda_{q^-}^-) - \mu (\lambda_1^+ \otimes \dots \otimes \lambda_{q^+}^+) \in \Z$.
Its cardinality is the difference $q (\gamma) =  q^+ - q^-$.
\end{defi}

Let $(\gamma_i)_{i \in I}$ be a family of elementary trajectories between two open strings, a same trajectory being allowed
to appear several times in the family. We denote by $\oplus_{i \in I} \gamma_i$ the union of these trajectories counted with multiplicities,
so that for every trajectory $\gamma$, $\gamma \oplus \gamma = 2\gamma \neq \gamma$.

\begin{defi}
\label{deftrajectory}
A  trajectory from the open string $\lambda_1^+ \otimes \dots \otimes \lambda_{q^+}^+$ to the open 
string $\lambda_1^- \otimes \dots \otimes \lambda_{q^-}^-$ is a sum $\oplus_{i \in I} \gamma_i$, where
$(\gamma_i)_{i \in I}$ is a family of elementary trajectories from $\lambda_1^+ \otimes \dots \otimes \lambda_{q^+}^+$ to
 $\lambda_1^- \otimes \dots \otimes \lambda_{q^-}^-$ which satisfies the Novikov condition
 $\forall C \in \R, \, \# \{ i \in I \, \vert \, a(\gamma_i)  < C \} < +\infty.$
\end{defi}

It would be nice to replace the Novikov condition in Definition \ref{deftrajectory} by the stronger condition
$\sum_{i \in I} h (\gamma_i) < +\infty$. However, this is related to a still unsolved conjecture of convergence
in Lagrangian Floer theory, so cannot be done yet. The set of trajectories between two open strings has
the structure of an Abelian group. The group law is given by the operation $\oplus$ whereas the opposite
of a trajectory is obtained by flipping all the orientations of the associated Cauchy-Riemann operators. 

The set of trajectories given by Definition \ref{deftrajectory} will be denoted by $\text{Hom} ({\cal OS} (X, \omega))$.
If $\gamma^1 = \oplus_{i \in I} \gamma_i^1 : \lambda_1^+ \to  \lambda_1^-$ and $\gamma^2 = \oplus_{j \in J} \gamma_j^2 : \lambda_2^+ \to  \lambda_2^-$ are two trajectories, we denote by $\gamma^1 \otimes \gamma^2 :  \lambda_1^+  \otimes \lambda_2^+ \to \lambda_1^-  \otimes \lambda_2^-$ the trajectory $\oplus_{(i,j) \in I \times J}  (\gamma_i^1 \cup \gamma_j^2)$. Likewise, if
$\gamma^1 = \oplus_{i \in I} \gamma_i^1 : \lambda^+ \to  \lambda^0$ and $\gamma^2 = \oplus_{j \in J} \gamma_j^2 : \lambda^0 \to  \lambda^-$ are two trajectories, we denote by $\gamma^2 \circ \gamma^1 :  \lambda^+ \to \lambda^-$ the trajectory
$\oplus_{(i,j) \in I \times J}  (\gamma_i^1 \cup \gamma_j^2)$.

\begin{defi}
\label{defcotrajectory}
The dual cotrajectory of a trajectory $\gamma : \lambda_+ \to \lambda_-$ is the cotrajectory $\gamma^* = \gamma :
\lambda_-^* \to \lambda_+^*$
\end{defi}
The dual of an empty trajectory $\emptyset_\lambda : \lambda \to \lambda$ is the empty cotrajectory
$\emptyset_\lambda^* : \lambda^* \to \lambda^*$. We denote by $\text{Hom} ({\cal OS}^* (X, \omega))$ the
set of cotrajectories given by Definition \ref{defcotrajectory}. We then denote by  ${\cal OS} (X, \omega) $
(resp. ${\cal OS}^* (X, \omega) $) the pair $\big( \text{Ob} ({\cal OS} (X, \omega)) , \text{Hom} ({\cal OS} (X, \omega)) \big)$
(resp. $\big( \text{Ob} ({\cal OS}^* (X, \omega)) , \text{Hom} ({\cal OS}^* (X, \omega)) \big)$).

\begin{prop}
\label{propcategoryos}
Let $(X, \omega) $ be a symplectic manifold of dimension at least four. Then, ${\cal OS} (X, \omega)$
and ${\cal OS}^* (X, \omega)$
have the structure of small preadditive categories dual to each other. Moreover, they are equipped with a tensor product
defined on objects and morphisms together with an action of the integers.
\end{prop}

{\bf Proof:}

The composition of morphisms is given by the operation $\circ$, it is associative. The identity morphism of every object
is given by the empty (co)trajectory. Objects and morphisms of these categories are sets and morphisms between
two objets are Abelian groups. The tensor product is given by the operation $\otimes$ whereas
the group of integers acts on $\text{Ob} ({\cal OS} (X, \omega))$ and $\text{Ob} ({\cal OS}^* (X, \omega)) $ by
$(e ,  \lambda_1 \otimes \dots \otimes \lambda_q) \mapsto \lambda_1^{+e} \otimes \lambda_2 \otimes \dots \otimes \lambda_q$ for every
$e \in \Z$. The elements of $\text{Hom} ({\cal OS} (X, \omega))$ and $\text{Hom} ({\cal OS}^* (X, \omega))$ are
equivariant for these actions. Hence the result. $\square$\\

For every $N \in \N$, we denote by ${\cal OS}_N (X, \omega)$ and ${\cal OS}_N^* (X, \omega)$ the quotient
of ${\cal OS} (X, \omega)$, ${\cal OS}^* (X, \omega)$ by the subgroup $N\Z$ of $\Z$, so that
$\text{Ob} ({\cal OS}_N (X, \omega)) = \text{Ob} ({\cal OS} (X, \omega))/N\Z$,
$\text{Hom} ({\cal OS}_N (X, \omega)) = $ $ \text{Hom} ({\cal OS} (X, \omega)) $, 
${\cal OS}_0 (X, \omega) = {\cal OS} (X, \omega)$ and ${\cal OS}_0^* (X, \omega) = {\cal OS}^* (X, \omega)$.
The objects of ${\cal OS}_1 (X, \omega)$ are principal spaces over the integers, canonically isomorphic to $\Z$.
Morphisms of ${\cal OS}_1 (X, \omega)$ are morphisms of principal spaces. 

\begin{rem}
\label{remquantum}
All the categories of the family ${\cal OS} (X, t \omega)$, $t \in \R^*_+$, are isomorphic to each other, only the parameters
$h (\gamma) = \exp(- a(\gamma))$ of trajectories $\gamma$ vary. When $t$ goes to infinity, only the trajectories $\gamma$
with $h (\gamma)$ infinitely close to one survive, so that ${\cal OS} (X, t \omega)$ converges to the classical category
whose objects are finite collections of points in the space $X$ and morphisms are homotopy classes of embedded forests having one negative 
root on every tree. We postpone a detailed description of the latter.
\end{rem}

\subsection{Functors coefficients}

Set $\Z ((t^\R))  = \{ \sum_{a \in \R} n_a t^a \, \vert \, \forall C \in \R, \, \# \{ a < C \, \vert \, n_a \neq 0 \in \Z \} < \infty \}$.
Let ${\cal M}od^1_{\Z ((t^\R))}$ be the category of free modules of rank one over this Novikov ring. Denote by
$\overline{{\cal M}od}^1_{\Z ((t^\R))}$ the extended category defined by
$\text{Ob} (\overline{{\cal M}od}^1_{\Z ((t^\R))}) = \text{Ob} ({\cal M}od^1_{\Z ((t^\R))}) \cup \{ \Z \}$ and
$\text{Hom} (\overline{{\cal M}od}^1_{\Z ((t^\R))}) = \text{Hom} ({\cal M}od^1_{\Z ((t^\R))}) \cup 
\cup_{M \in \text{Ob} (\overline{{\cal M}od}^1_{\Z ((t^\R))})} \text{Hom}_\Z (\Z , M)$.
The aim of this paragraph is to define dual functors ${\cal C} :  {\cal OS} (X, \omega) \to \overline{{\cal M}od}^1_{\Z ((t^\R))}$ and
${\cal C}^* :  {\cal OS}^* (X, \omega) \to \overline{{\cal M}od}^1_{\Z ((t^\R))}$. 

Let $ \lambda$ be an elementary open string and $ \text{Hom}_{\cal OS} (\emptyset , \lambda)$ be the space of
elementary trajectories $\gamma : \emptyset \to \lambda$. This space is equipped with two involutions
$$\begin{array}{rcrcl}
c_{\p} & : & \text{Hom}_{\cal OS} (\emptyset , \lambda) & \to & \text{Hom}_{\cal OS} (\emptyset , \lambda)\\
&& (u, D_{\underline{z}} , \lambda_{\partial u} , \overline{\partial}_u ) & \mapsto & (u, D_{\underline{z}} , 
\overline{\lambda}_{\partial u} , \overline{\partial}_u ),
\end{array}$$
and 
$$\begin{array}{rcrcl}
c_{\overline{\partial}} & : & \text{Hom}_{\cal OS} (\emptyset , \lambda) & \to & \text{Hom}_{\cal OS} (\emptyset , \lambda)\\
&& (u, D_{\underline{z}} , \lambda_{\partial u} , \overline{\partial}_u ) & \mapsto & (u, D_{\underline{z}} , 
\lambda_{\partial u} , - \overline{\partial}_u ),
\end{array}$$
where $- \overline{\partial}_u$ stands for the same Cauchy-Riemann operator $\overline{\partial}_u$ but with
opposite orientation and $\overline{\lambda}_{\partial u}$ stands for the same path of Lagrangian subspaces
but switching the $\widetilde{GL}_n^\pm (\R)$-structure $\p^\pm$. There are indeed exactly two different
$\widetilde{GL}_n^\pm (\R)$-structures on a vector bundle over the circle or likewise two different extensions
over an interval of $\widetilde{GL}_n^\pm (\R)$-structures given on its boundary. Denote by 
$\vert  \text{Hom}_{\cal OS} (\emptyset , \lambda) \vert$ the quotient of $ \text{Hom}_{\cal OS} (\emptyset , \lambda)$ 
by the composition $c_{\p} \circ c_{\overline{\partial}} $ of these involutions.

\begin{lemma}
\label{lemmac}
For every elementary open string $ \lambda$ of a symplectic manifold $(X, \omega)$,
the space $\vert  \text{Hom}_{\cal OS} (\emptyset , \lambda) \vert$ has two connected components exchanged
by both involutions $c_{\p}$ and $c_{\overline{\partial}} $.
\end{lemma}

{\bf Proof:}

The space of supports $(u, D_{\underline{z}})$ of trajectories $\gamma : \emptyset \to \lambda$ retracts on the
constant map at the base point $\pi_{\cal L} \circ \lambda$ of the string. From Lemma \ref{lemmaconfig},
the quotient of $ \text{Hom}_{\cal OS} (\emptyset , \lambda)$ by the involution $c_{\overline{\partial}} $ is connected,
it is a nontrivial double cover of the quotient $ \text{Hom}_{\cal OS} (\emptyset , \lambda)/ \left< c_{\p} , c_{\overline{\partial}} \right>$.
We have to prove that the same holds for the double cover 
$ \text{Hom}_{\cal OS} (\emptyset , \lambda)/ \left< c_{\p} \right> \to  \text{Hom}_{\cal OS} (\emptyset , \lambda)/ \left< c_{\p} , c_{\overline{\partial}} \right>$.
For this purpose, we have to prove the existence of a loop $(l_{\partial u}^t)_{t \in S^1}$ of boundary Lagrangian conditions
such that the associated rank one real vector bundle $\det (\overline{\partial}_u^t) \to S^1$ is not orientable. The existence
of such a loop was observed in \cite{FOOO}, we propose here a different proof. First of all, it suffices to prove
this result for $n=2$ and from the linear gluing Theorem, see for example Theorem $10$ of \cite{FloerHofer}, the thesis
\cite{Schw} or \cite{BirCor}, it suffices to prove this result for a loop of Lagrangian boundary conditions on the
closed disk, gluing a fixed trajectory $\emptyset \to \lambda^*$. Let $F_t \to \C P^1$, $t \in \C$, be the non-trivial Kodaira deformation
of rank two vector bundles with special fiber $F_0 \cong {\cal O}_{\C P^1} (-1) \oplus {\cal O}_{\C P^1} (1) \to \C P^1$.
This deformation is trivializable over $\C^*$, with fiber $F_t \cong {\cal O}_{\C P^1}  \oplus {\cal O}_{\C P^1}  \to \C P^1$,
$t \in \C^*$, so that it extend to a deformation $F_t \to  \C P^1$, $t \in  \C P^1$. Fix a hemisphere $H$ of $ \C P^1$ and denote
by $\R F_t \to \partial H$ the real locus of $F_t$, $t \in  \R P^1$. The loop of associated Cauchy-Riemann operators
$\overline{\partial}_t : L^{k,p} (H ; F_t , \R F_t) \to L^{k-1,p} (H ; \Lambda^{0,1} H \otimes F_t )$, $t \in  \R P^1$, is not
orientable. Indeed, these operators are all surjectives and their kernels are two-dimensional. For every $t \in  \R P^1 \setminus \{ 0 \}$
and every $z \in \partial H$, the evaluation map at the point $z$ provides an isomorphism between $\ker (\overline{\partial}_t)$
and $\R F_t \vert_z$. Let us look at the behavior of this evaluation map in the neighborhood of $t=0$.
Let $U_0$, $U_1$ be the standard affine charts of $\C P^1$. The vector bundle $F_t$ is trivial over $U_0$, $U_1$,
the change of trivializations is given by the matrix 
$\left[
\begin{array}{cc}
z & t \\
0 & \frac{1}{z}
\end{array}
\right]$, $t \in \C$, $z \in U_0 \setminus \{ 0 \}$.
The kernel of $\overline{\partial}_t$ is generated by the sections $z  \in U_0  \mapsto (t , -z) \in \C^2$ and
$z  \in U_0  \mapsto (0 , 1) \in \C^2$. Hence, for every $z \in \R U_0 \setminus \{ 0 \}$, $t= 0$ is a non-degenerated critical point
of the evaluation map $\ker (\overline{\partial}_t) \to \R F_t \vert_z$. Since this is the only critical point, the family
$(\overline{\partial}_t )_{t \in  \R P^1}$ is not
orientable. $\square$\\

Let $\Z [\vert  \text{Hom}_{\cal OS} (\emptyset , \lambda) \vert ]$ be the free Abelian group of rank two generated
by the connected components of $\vert  \text{Hom}_{\cal OS} (\emptyset , \lambda) \vert$ and $\Z_\lambda$
be the kernel of the morphism 
$$a[\gamma] + b [c_{\p} (\gamma)] \in \Z [\vert  \text{Hom}_{\cal OS} (\emptyset , \lambda) \vert ]
\mapsto a + b \in \Z.$$ 
Set 
$\Z_\lambda ((t^\R))  = \{ \sum_{a \in \R} n_a t^a \, \vert \, \forall C \in \R, \, \# \{ a < C \, \vert \, n_a \neq 0 \in \Z_\lambda \} < \infty \}$,
it is a free module of rank one over $\Z ((t^\R))$. Now for every elementary trajectory $\gamma : \emptyset \to \lambda$, denote by
$\vert \gamma \vert$ its projection in $\vert  \text{Hom}_{\cal OS} (\emptyset , \lambda) \vert$ and by ${\cal C} (\gamma)(1)$
the element $\vert \gamma \vert t^{a(\gamma)} \in \Z_\lambda ((t^\R))$. For every trajectory $\gamma = \oplus_{i \in I} \gamma_i : \emptyset \to \lambda$,
denote by ${\cal C} (\gamma)(1)$ the element $\sum_{i \in I} {\cal C} (\gamma_i)(1)$ and then by ${\cal C} (\gamma)$ the morphism
$n \in \Z \mapsto n {\cal C} (\gamma)(1) = {\cal C} (\gamma)(n) \in \Z_\lambda ((t^\R))$. Likewise, for every elementary trajectory
$\gamma : \lambda_1^+ \otimes \dots \otimes \lambda_{q^+}^+ \to \lambda_1^- \otimes \dots \otimes \lambda_{q^-}^-$,
the composition with $\gamma$ in  ${\cal OS} (X, \omega)$ provides from the linear gluing Theorem (see Theorem $10$ \cite{FloerHofer}, the thesis
 \cite{Schw} or \cite{BirCor}) a map $\vert  \text{Hom}_{\cal OS} (\emptyset , \lambda_1^+) \vert \times \dots \times
\vert  \text{Hom}_{\cal OS} (\emptyset , \lambda_{q^+}^+) \vert \to \vert  \text{Hom}_{\cal OS} (\emptyset , \lambda_1^-) \vert \times \dots \times
\vert  \text{Hom}_{\cal OS} (\emptyset , \lambda_{q^-}^-) \vert$. We deduce from this map a morphism
 $\Z_{\lambda_1^+} \otimes \dots \otimes \Z_{ \lambda_{q^+}^+} \to \Z_{\lambda_1^-} \otimes \dots \otimes \Z_{\lambda_{q^-}^-}$
 and hence after multiplication by $ t^{a(\gamma)} $ a morphism 
 $\Z_{\lambda_1^+} ((t^\R)) \otimes \dots \otimes \Z_{ \lambda_{q^+}^+} ((t^\R)) \to \Z_{\lambda_1^-} ((t^\R)) \otimes \dots \otimes 
 \Z_{\lambda_{q^-}^-} ((t^\R))$ denoted by $ {\cal C} (\gamma)$. Once more, if $\oplus_{i \in I} \gamma_i $ is a trajectory
 $ \lambda_1^+ \otimes \dots \otimes \lambda_{q^+}^+ \to \lambda_1^- \otimes \dots \otimes \lambda_{q^-}^-$, we set
${\cal C} (  \oplus_{i \in I} \gamma_i ) =   \sum_{i \in I} {\cal C} (\gamma_i)$ and $ {\cal C} (\emptyset) = id$.

\begin{defi}
\label{deffoncteurc}
The functors coefficients are the functors ${\cal C} :  {\cal OS} (X, \omega) \to \overline{{\cal M}od}^1_{\Z ((t^\R))}$ and
${\cal C}^* :  {\cal OS}^* (X, \omega) \to \overline{{\cal M}od}^1_{\Z ((t^\R))}$ defined by

$\begin{array}{ccccc}
{\cal C} & : &  \text{Ob} ( {\cal OS} (X, \omega) ) & \to & \text{Ob} ( \overline{{\cal M}od}^1_{\Z ((t^\R))} ) \\
&& \emptyset & \mapsto & \Z \\
&&  \lambda_1 \otimes \dots \otimes \lambda_q & \mapsto & \Z_{\lambda_1} ((t^\R)) \otimes \dots \otimes \Z_{ \lambda_{q}} ((t^\R)), \\
{\cal C} & : &  \gamma \in \text{Hom} ( {\cal OS} (X, \omega) ) & \mapsto & {\cal C} (\gamma) \in \text{Hom} ( \overline{{\cal M}od}^1_{\Z ((t^\R))} ) \\
\end{array}$

 and

$\begin{array}{ccccc}
{\cal C}^* & : &  \text{Ob} ( {\cal OS}^* (X, \omega) ) & \to & \text{Ob} ( \overline{{\cal M}od}^1_{\Z ((t^\R))} ) \\
&& \emptyset^* & \mapsto & \Z \\
&&  \lambda_1 \otimes \dots \otimes \lambda_q & \mapsto & \Z_{\lambda_1} ((t^\R)) \otimes \dots \otimes \Z_{ \lambda_{q}} ((t^\R)), \\
{\cal C}^* & : &  \gamma^* \in \text{Hom} ( {\cal OS}^* (X, \omega) ) & \mapsto & {\cal C} (\gamma)^* \in \text{Hom} ( \overline{{\cal M}od}^1_{\Z ((t^\R))} ) \\
\end{array}$ 
\end{defi}

This Definition \ref{deffoncteurc} gets justified by the following Proposition \ref{propfunctorc}.

\begin{prop}
\label{propfunctorc}
The coefficients functors ${\cal C} $, ${\cal C}^*$ given by Definition \ref{deffoncteurc} are indeed functors
of small preadditive categories.
\end{prop}

{\bf Proof:}

This Proposition \ref{propfunctorc} follows tautologically from the definitions of ${\cal C} $ and ${\cal C}^*$. $\square$\\

The duality between $ {\cal OS} (X, \omega) $ and $ {\cal OS}^* (X, \omega) $ writes as a map from $ {\cal OS} (X, \omega) \times
 {\cal OS}^* (X, \omega) $ to the category of trivial closed strings. These trivial closed strings admit an analogous functor coefficient
 which attach to them a free module of rank one $\tilde{\Z} ((t^\R))$ over $\Z ((t^\R)) $, where $\tilde{\Z} $ is a torsion free cyclic group
 independant of the string analogous to $\Z_\lambda$'s. Since we won't use this duality here, we postpone a detailed description of it.
 
 \subsection{Chain complexes and augmentations}
 \label{subsectionchain}
 
 For every $N \in \N$, denote by $K^b ({\cal OS}_N (X, \omega))$ the homotopy category of bounded
 chain complexes of ${\cal OS}_N (X, \omega)$. We want the objects and morphisms of this category to be compatible
 with the tensor product and duality defined on ${\cal OS}_N (X, \omega)$. This leads to a list of axioms which are
 described in this paragraph.
 
 Let $(\Lambda , \delta^\Lambda)$ be a bounded chain complex of ${\cal OS}_N (X, \omega)$. Such a complex writes
 $$\Lambda^0 \stackrel{\delta^0}{\to} \Lambda^1 \stackrel{\delta^1}{\to} \dots  \stackrel{\delta^{d(\Lambda) - 1}}{\to} \Lambda^{d(\Lambda)},$$
 where for every $0 \leq j \leq d(\Lambda)$, $\Lambda^j$ is a finite sum of objects of ${\cal OS}_N (X, \omega)$ and
 $\delta^j$ a finite sum of morphisms. It is indeed convenient to restrict ourselves to  chain complexes starting at grading $0$.
 It may also sometimes be convenient to allow infinite sum of objects and morphisms, see \S \ref{subsubsectbased}. In particular,
 a finite sum of objects in ${\cal OS}_1 (X, \omega)$ reads as an infinite sum of objects in ${\cal OS} (X, \omega)$.
 Each factor $\Lambda^j$ gets filtrated by the cardinalities of its open strings, so that $\Lambda^j = \sum_{q=0}^{l_j} \Lambda^j_q$,
 where $ \Lambda^j_q$ is a finite sum of open strings of cardinality $q$, $l_j \in \N$ and  $0 \leq j \leq d(\Lambda)$.
 Likewise, every morphism $\delta^j$ decomposes $\delta^j = \sum_{q \in \Z} \delta^j_q$, where
 $\delta^j_q$ is a finite sum of trajectories of cardinality $q$, so that $\delta^j_q$ increases the graduation of the complex by one
 and decreases its filtration by $q$. Denote by $\text{Ob} ( K^b ({\cal OS}_N (X, \omega)))$ the set of bounded chain complexes
 $(\Lambda , \delta^\Lambda)$ of ${\cal OS}_N (X, \omega)$ which satisfy the following three axioms $A_1$, $A_2$, $A_3$.\\
 
 $A_1$ : For every $0 \leq j \leq d(\Lambda)$ and $1 \leq i \leq q$, $\lambda_1 \otimes \dots \otimes \lambda_q \in  \Lambda^j
 \implies \lambda_i \in  \Lambda^j$.\\
 
 $A_2$ : For every $0 \leq j \leq d(\Lambda)$ and $q \in \N$, $\mu (\delta^j_q) = 1-q$.\\
 
 $A_3$ : For every $0 \leq j \leq d(\Lambda)-1$, if 
 $\lambda_1 \otimes \lambda_2 \in (\Lambda^{j+1})^*$, 
 $\delta^j$ satisfies the following Leibniz formula
$$(\delta^j)^* (\lambda_1 \otimes \lambda_2) = \Big( \sum_{q=0}^{+ \infty} (-1)^{q(\lambda_2)} (\delta^j_q)^*  \otimes id_{q(\lambda_2)}
+ (-1)^{q(\lambda_1) (q+1)} id_{q(\lambda_1)}  \otimes (\delta^j_q)^* \Big) (\lambda_1 \otimes \lambda_2).$$

\begin{rem}
\label{remgraduation}
1) One cannot hope that $\delta^\Lambda$ satisfies a Leibniz formula of the form
$(\delta^\Lambda)^* (\lambda_1 \otimes \lambda_2) = \big( \sum_{q=0}^{+ \infty} (\delta^\Lambda_q)^*  \otimes id_{q(\lambda_2)}
+ (-1)^{s (q)} id_{q(\lambda_1)}  \otimes (\delta^\Lambda_q)^* \big) (\lambda_1 \otimes \lambda_2)$ for some $s(q)$,
since then $\delta^\Lambda \circ \delta^\Lambda$ would contain the term $2 (-1)^{s (q)} (\delta^\Lambda_0)^* \otimes  (\delta^\Lambda_1)^*$.

2) The differentials $(\delta^j)^*$ are completely determined by their values on elementary open strings, so that from
axioms $A_1$, $A_3$,  for every $q \in \N^*$, $l \in \N$, 
$(\delta^l)^*  (\lambda_1 \otimes \dots \otimes \lambda_q) = \sum_{i=1}^q (-1)^{(i-1)(l+1) + q-i} id_{i-1}  \otimes (\delta^l)^*
\otimes id_{q-i} (\lambda_1 \otimes \dots \otimes \lambda_q)$.

3) From axiom $A_2$, the levels of the function $\mu + q$ provide a graduation of $\Lambda$ defined modulo
$N$. The differential $ \delta^\Lambda$ is also of degree one for this graduation.
\end{rem}

Now, let $(\Lambda , \delta^\Lambda)$, $(\Lambda' , \delta^{\Lambda'})$ be two elements of $\text{Ob} ( K^b ({\cal OS}_N (X, \omega)))$
and $H : (\Lambda , \delta^\Lambda) \to (\Lambda' , \delta^{\Lambda'})$ be a degree zero chain map, so that for every
$0 \leq j \leq d(\Lambda)-1$, $\delta'^j \circ H^j = H^{j+1} \circ \delta^j$. This chain map is a finite sum of trajectories and
decomposes as $H^j =  \sum_{q \in \N} H^j_q$ where $H^j_q$ is a finite sum of trajectories of cardinality $q$. Such a chain map
is said to be a morphism from $(\Lambda , \delta^\Lambda)$ to $(\Lambda' , \delta^{\Lambda'})$ if and only if it satisfies the following
two axioms $B_1$, $B_2$.\\

$B_1$ : For every $0 \leq j \leq d(\Lambda)$ and $q \in \N$, $\mu (H^j_q) = -q$.\\

$B_2$ : For every $0 \leq j \leq d(\Lambda')$, if 
$\lambda_1 \otimes \lambda_2 \in (\Lambda'^{j})^*$, then  
 
 $$(H^j)^* (\lambda_1 \otimes \lambda_2) = \Big( \sum_{l_1 + l_2 = q} (-1)^{q(\lambda_1) l_2} (H^j_{l_1})^*  \otimes  (H^j_{l_2})^* \Big) 
 (\lambda_1 \otimes \lambda_2).$$
 
 \begin{rem}
 \label{remh}
 From axiom $B_2$, the morphism $H^*$ is completely determined by its values on elementary strings, so that
 $H^*  (\lambda_1 \otimes \dots \otimes \lambda_q) = \big( \sum_{l_1 , \dots , l_q \in \N} (-1)^{\sum_{i=1}^{q-1} i l_{i+1}}  H_{l_1}^*  \otimes \dots 
\otimes H_{l_q}^* \big) (\lambda_1 \otimes \dots \otimes \lambda_q)$.
 \end{rem}
 
 \begin{lemma}
 The axioms $A$ and $B$ are consistent to each other.
 \end{lemma}
 
 {\bf Proof:}
 $$\delta^* \circ H^* =  \sum_{q, l_1 , l_2 \in \N} (-1)^{q(\lambda_1) l_2} \big(
 (-1)^{q(\lambda_2) + l_2} (\delta_q^* \circ H_{l_1}^*)  \otimes  H_{l_2}^*  + (-1)^{q(\lambda_1)(q+1)}
 H_{l_1}^* \otimes (\delta_q^* \circ H_{l_2}^*) \big),$$
and
 $$H^* \circ  \delta^*  =  \sum_{q, l_1 , l_2 \in \N} \big( (-1)^{q(\lambda_2) + (q(\lambda_1) + 1) l_2} 
 (H_{l_1}^* \circ \delta_q^* )  \otimes  H_{l_2}^*  + (-1)^{q(\lambda_1)(q+1 + l_2)}
 H_{l_1}^* \otimes (H_{l_2}^* \circ \delta_q^* ) \big),$$
 from the commutation relations $(id_{q(\lambda_1)} \otimes \delta_q^*) \circ (H_{l_1}^*  \otimes H_{l_2}^*) =
 (-1)^{l_1 (q+1)} H_{l_1}^*  \otimes (\delta_q^* \circ H_{l_2}^*)$ and $(H_{l_1}^*  \otimes  H_{l_2}^*)   \circ (\delta_q^* \otimes
 id_{q(\lambda_2)}) = (-1)^{l_2 (q+1)} (H_{l_1}^* \circ \delta_q^*) \otimes  H_{l_2}^*$. $\square$\\

Let finally $K : (\Lambda , \delta^\Lambda) \to (\Lambda' , \delta^{\Lambda'})$ be a degree $-1$ map,  which once more has
a decomposition of the form  $K^j =  \sum_{q \in \N} K^j_q$ where $K^j_q$ is a finite sum of trajectories of cardinality $q$,
$0 \leq j \leq d(\Lambda)$. Such a map is said to be a homotopy from $(\Lambda , \delta^\Lambda)$ to $(\Lambda' , \delta^{\Lambda'})$ if and only if 
it satisfies the following
two axioms $C_1$, $C_2$.\\

$C_1$ : For every $0 \leq j \leq d(\Lambda)$ and $q \in \N$, $\mu (K^j_q) = -1-q$.\\

$C_2$ : For every $0 \leq j \leq d(\Lambda) + 1$, if 
$\lambda_1 \otimes \lambda_2 \in (\Lambda'^{j})^*$, there exists a morphism $H :  (\Lambda , \delta^\Lambda) \to (\Lambda' , \delta^{\Lambda'})$
 such that 
 $$K_q^* (\lambda_1 \otimes \lambda_2) = \sum_{l_1 + l_2 = q} \big( (-1)^{q(\lambda_2) + (q(\lambda_1) + 1) l_2} K^*_{l_1}  \otimes  H_{l_2}^*
 + (-1)^{q(\lambda_1) (l_2 + 1)} H^*_{l_1}  \otimes  K_{l_2}^*\big) (\lambda_1 \otimes \lambda_2).$$
 
 \begin{rem}
 \label{remhC}
 From axiom $C_2$, the morphism $K^*$ is completely determined by its values on elementary strings, so that
 $K^*  (\lambda_1 \otimes \dots \otimes \lambda_q) = \sum_{l_1 , \dots , l_q \in \N} \sum_{i=1}^q
 (-1)^{q-1 + \sum_{j=1}^{q-1} j l_{j+1} + \sum_{j=i+1}^{q} l_j}  H_{l_1}^*  \otimes \dots \otimes H_{l_{i-1}}^*  \otimes K_{l_{i}}^* \otimes H_{l_{i+1}}^*  \otimes
\dots \otimes H_{l_q}^* (\lambda_1 \otimes \dots \otimes \lambda_q)$.
 \end{rem}

 \begin{lemma}
 The axioms $C$ are consistent with the axioms $A$ and $B$, that is if $K$ satisfies axioms $C$ and $\delta$ satisfies axioms $A$,
 then $\delta \circ K + K \circ \delta$ satisfies axioms $B$.
 \end{lemma}
 
 {\bf Proof:}
 
 $$\delta^* \circ K^* =  \sum_{q, l_1 , l_2 \in \N}  \Big( (-1)^{q(\lambda_2) + (q(\lambda_1) + 1) l_2}  \big(
 (-1)^{q(\lambda_2) + l_2} (\delta_q^* \circ K_{l_1}^*)  \otimes  H_{l_2}^*  + (-1)^{(q(\lambda_1) + 1)(q+1)}
 K_{l_1}^* \otimes (\delta_q^* \circ H_{l_2}^*) \big) $$ $$
 + (-1)^{q(\lambda_1) (l_2 + 1)}  \big( (-1)^{q(\lambda_2) + l_2} (\delta_q^* \circ H_{l_1}^*)  \otimes  K_{l_2}^*  + 
 (-1)^{q(\lambda_1) (q + 1)}  H_{l_1}^* \otimes (\delta_q^* \circ K_{l_2}^*) \big) \Big),$$
and
$$K^*  \circ \delta^* =  \sum_{q, l_1 , l_2 \in \N}  \Big( (-1)^{q(\lambda_2)}  \big(
 (-1)^{q(\lambda_2) + q(\lambda_1)  l_2} (K_{l_1}^* \circ  \delta_q^* )  \otimes  H_{l_2}^*  + (-1)^{(q(\lambda_1) + 1)(l_2+1)}
 (H_{l_1}^*  \circ  \delta_q^* ) \otimes K_{l_2}^* \big) $$ $$
 + (-1)^{q(\lambda_1) (q + 1)}  \big( (-1)^{q(\lambda_2) + q + (q(\lambda_1) + 1)l_2} K_{l_1}^*  \otimes  (H_{l_2}^* \circ  \delta_q^* )  + 
 (-1)^{q(\lambda_1) (l_2 + 1)}  H_{l_1}^* \otimes (K_{l_2}^* \circ \delta_q^* ) \big) \Big). \square$$

 Denote by $\text{Hom} ( K^b ({\cal OS}_N (X, \omega)))$ the set of degree zero chain maps between elements of
 $\text{Ob} ( K^b ({\cal OS}_N (X, \omega)))$  satisfying axioms $B$ modulo degree $-1$ homotopies satisfying axioms $C$,
 so that if $H - H' = \delta \circ K + K \circ \delta$, then the morphisms  $H$ and $H'$ are the same in 
 $\text{Hom} ( K^b ({\cal OS}_N (X, \omega)))$. The composition of two morphisms is a morphism so that $K^b ({\cal OS}_N (X, \omega)) =
(\text{Ob} ( K^b ({\cal OS}_N (X, \omega))) , \text{Hom} ( K^b ({\cal OS}_N (X, \omega))))$ is a small preadditive category.
 
 \begin{defi}
 \label{defaugmentation}
 Let $(\Lambda , \delta^\Lambda) \in \text{Ob} ( K^b ({\cal OS}_N (X, \omega)))$. An augmentation of $(\Lambda , \delta^\Lambda)$
 is a morphism $\delta = \sum_{\lambda \in \Lambda^0} \delta_{\emptyset , \lambda} : \emptyset \to  \Lambda^0$ such that:
 
 1) $\delta^0 \circ \delta = 0$.
 
 2) $ \delta_{\emptyset , \lambda_1 \otimes \dots \otimes \lambda_q} = (-1)^{\sum_{i=1}^q (q-i) \mu (\lambda_i)} \delta_{\emptyset , \lambda_1}
\otimes \dots \otimes \delta_{\emptyset , \lambda_q} $, whenever $\lambda_1 , \dots , \lambda_q$ are elementary strings of $\Lambda^0$
such that $\lambda_1 \otimes \dots \otimes \lambda_q \in  \Lambda^0$.

A complex equipped with such an augmentation is called an augmented complex.
 \end{defi}
 
 \begin{lemma}
 Let $H : (\Lambda , \delta^\Lambda) \to (\Lambda' , \delta^{\Lambda'}) \in \text{Hom} ( K^b ({\cal OS}_N (X, \omega)))$ 
 and $\delta : \emptyset \to  \Lambda^0$ be an augmentation of $(\Lambda , \delta^\Lambda) $. Then,
 $H^0 \circ \delta$ in an  augmentation of $(\Lambda' , \delta^{\Lambda'}) $. 
\end{lemma} 

{\bf Proof:}

In the relation $(H_{l_1} \otimes \dots \otimes H_{l_q}) \otimes (\emptyset_1  \otimes \dots \otimes \emptyset_l) =
(-1)^{\sum_{j=1}^{q-1} (l_{j+1} - 1 + \dots + l_q - 1) \sum_{i= l_1 + \dots + l_{j-1} + 1}^{l_1 + \dots + l_j} \mu_i}
(H_{l_1} \otimes \emptyset_1  \otimes \dots \otimes \emptyset_{l_1}) \otimes \dots \otimes
(H_{l_q} \otimes \emptyset_{l_1 + \dots + l_{q-1} + 1}  \otimes \dots \otimes \emptyset_{l})$,
the $j^{\text{th}}$ term is counted with respect to the sign  $(q-j) (l_j - 1) + 
\sum_{i= l_1 + \dots + l_{j-1} + 1}^{l_1 + \dots + l_j} (l-i) \mu_i \mod (2)$ by the left hand side whereas it is counted
with respect to the sign $(q-j)(\sum_{i= l_1 + \dots + l_{j-1} + 1}^{l_1 + \dots + l_j} \mu_i  + l_j -1) +
\sum_{i= l_1 + \dots + l_{j-1} + 1}^{l_1 + \dots + l_j} (l_1 + \dots + l_j - i)\mu_i \mod(2)$ by the right hand side.
These signs coincide. $\square$\\

\section{Lagrangian conductors}

This paragraph is devoted to the construction of the category ${\cal CL}^\pm (X , \omega)$ of Lagrangian conductors
of $(X , \omega)$. Once more, $(X , \omega)$ is assumed to be of dimension at least four, we do not consider the
special case of surfaces.

\subsection{Singular almost-complex structures}
\label{subsectionsingular}

\begin{defi}
A $S$-neck of the manifold $(X , \omega)$  is an embedding $\phi : S \times [- \epsilon , \epsilon ] \to X$
which satisfies $\phi^* \omega = d(e^t \theta)$, where $(S , \theta)$ is a closed contact manifold of dimension
$2n-1$, $\epsilon \in \R_+^*$ and $t \in [- \epsilon , \epsilon ]$.
\end{defi}

\begin{defi}
\label{defssing}
An almost-complex structure $J$ is called $S$-singular if there exists a $S$-neck $\phi : S \times [- \epsilon , \epsilon ] \to X$
such that:

1) The domain of definition of $J$ is the complement $X \setminus \phi (S \times \{ 0 \})$.

2) The almost-complex structure $\phi^* J$ preserves the contact distribution $\ker (\theta)$
of $S \times \{ t \}$ for every $t \in  [- \epsilon , \epsilon ]   \setminus \{ 0 \}$ and its restriction to
$\ker (\theta)$ does not depend on $t \in  [- \epsilon , \epsilon ]   \setminus \{ 0 \}$.

3) $\forall (x,t) \in S \times ([- \epsilon , \epsilon ]   \setminus \{ 0 \})$, $\phi^* J (\frac{\partial}{\partial t}) \vert_{(x,t)} =
\alpha' (t) R_\theta \vert_{(x,t)}$, where $\alpha' : [- \epsilon , \epsilon ]   \setminus \{ 0 \} \to \R_+^*$ is even with infinite
integral and $R_\theta$ denotes the Reeb vector field of $(S , \theta)$.
\end{defi}

\begin{defi}
An almost-complex structure $J$ of $(X , \omega)$ is called singular if it is $S$-singular for some
$(2n-1)$-dimensional contact manifold $(S , \theta)$.
\end{defi}

Denote by $\partial {\cal J}_\omega$ the space of singular almost-complex structures of $X$ compatible
with $\omega$. It is equipped with the following topology. A singular almost-complex structure $J$
is said to be in the $\eta$-neighborhood of $J_0 \in \partial {\cal J}_\omega$, $\eta > 0$, if these structures
are $S$-singular for the same contact manifold $(S , \theta)$ and if there exists pairs $(\phi_0 , \alpha'_0)$
and $(\phi , \alpha')$ given by Definition \ref{defssing} such that:

1) The distance between $\phi$ and $\phi_0$ is less than $\eta$. This distance in the space of embeddings
of finite regularity is induced by some fixed metric on $X$. The regularity of these embeddings is one more
than the regularity of the almost-complex structures which throughout the paper is supposed to be finite
and much bigger than the Sobolev regularity $k$ chosen for spaces of sections.

2) There exists $0 < \delta < \epsilon$ such that $2\eta \int_\delta^\epsilon \alpha'_0 (t) dt > 1$ and
the distance between the restrictions of $J$ and $J_0$ to the complement $X \setminus \phi_0 (S \times ]- \delta , \delta [)$
is less than $\eta$.

\begin{defi}
\label{defjsneck}
An almost-complex structure $J \in {\cal J}_\omega$ is said to have an $S$-neck if $X$
has an $S$-neck $\phi : S \times [- \epsilon , \epsilon ] \to X$ such that

1) The almost-complex structure $\phi^* J$ preserves the contact distribution $\ker (\theta)$
of $S \times \{ t \}$ for every $t \in  [- \epsilon , \epsilon ]$ and its restriction to
$\ker (\theta)$ does not depend on $t \in  [- \epsilon , \epsilon ] $.

2) $\forall (x,t) \in S \times [- \epsilon , \epsilon ]$, $\phi^* J (\frac{\partial}{\partial t}) \vert_{(x,t)} =
\alpha' (t) R_\theta \vert_{(x,t)}$, where $\alpha' : [- \epsilon , \epsilon ]   \to \R_+^*$ is even.

The integral $\int_{- \epsilon}^\epsilon \alpha (t) dt  $ is called the length of the neck.
\end{defi} 

Hence, an $S$-singular almost-complex structure is an almost-complex structure having an $S$-neck
of infinite length. This terminology comes from symplectic field theory \cite{EGH}. Indeed, if $J \in {\cal J}_\omega$
has an $S$-neck and $\alpha$ is the odd primitive of the function $\alpha'$ given by Definition \ref{defjsneck},
then, the diffeomorphism $(x,t) \in S \times [- \epsilon , \epsilon ] \mapsto (x, \alpha (t)) \in S \times [\alpha (-\epsilon) , \alpha (\epsilon) ]$
pushes forward $J$ to an almost-complex structure which preserves the contact distribution and sends the Liouville
vector field $\frac{\partial}{\partial t}$ onto the Reeb vector field $R_\theta$, compare \S $2.2$ of \cite{HWZfinite}.
In the language of symplectic field theory, a symplectic manifold $(X , \omega)$ equipped with a
$S$-singular almost-complex structure $J$ is an almost-complex manifold $(X \setminus \phi (S \times \{ 0 \}) , J)$
with cylindrical end.

Set $\overline{\cal J}_\omega = {\cal J}_\omega \sqcup \partial {\cal J}_\omega$ and equip this space with the following topology.
An almost-complex structure $J \in {\cal J}_\omega$
is said to be in the $\eta$-neighborhood of the $S$-singular almost-complex structure $J_0 \in \partial {\cal J}_\omega$, $\eta > 0$, 
if it has an $S$-neck and there exists pairs $(\phi_0 , \alpha'_0)$,
$(\phi , \alpha')$ given by Definitions \ref{defssing} and \ref{defjsneck} such that:

1) The distance between $\phi$ and $\phi_0$ is less than $\eta$ in the space of embeddings
of our fixed finite regularity.

2) There exists $0 < \delta < \epsilon$ such that $2\eta \int_\delta^\epsilon \alpha'_0 (t) dt > 1$ and
the distance between the restrictions of $J$ and $J_0$ to the complement $X \setminus \phi_0 (S \times ]- \delta , \delta [)$
is less than $\eta$.

 In particular, when $\eta$ is closed to zero, the length of the $S$-neck of $J$ is closed to infinity.
 
 Let us equip now the $n$-dimensional sphere $S^n$ with a metric $g$ having constant curvature. 
 For every $r >0$, denote by $U^*_r S^n = \{ (q,p) \in T^* S^n \, \vert \,  ||p||_g \leq r \}$
 and $S^*_r S^n = \partial U^*_r S^n$. The latter is equipped with the restriction $\theta$
 of the Liouville one-form of $T^* S^n$.
 
 \begin{defi}
 \label{defa1singular}
 An almost-complex structure $J \in \partial {\cal J}_\omega$ is said to be $A_1$-singular if there exists
 $r > 0$ such that it is $S^*_r S^n$-singular  and if the embedding $\phi : S^*_r S^n \times [- \epsilon , 0 ] \to X$
 given by Definition \ref{defssing} extends to a symplectic embedding $\phi : (U^*_r S^n ,d\theta) \to (X , \omega)$.
 \end{defi}
 Let $L$ be a Lagrangian sphere embedded in $(X , \omega)$, we denote by ${\cal J}_\omega^\infty (L) \subset \partial {\cal J}_\omega$
 the space of $A_1$-singular almost-complex structures $J$ for which $L \subset \phi (U^*_r S^n \setminus S^*_r S^n)$
 and $\phi^{-1} (L)$ is Hamiltonian isotopic to the zero section  of $U^*_r S^n$, where $\phi$ is the embedding given
 by Definition \ref{defa1singular}. Note that every Lagrangian sphere of $U^*_r S^n$ is conjectured to be Hamiltonian
  isotopic to the zero section, see \cite{FSS} for a recent work on this conjecture.
 
 {\bf Example:}
 
 Let $X_0 \in \C P^{n+1}$ be a projective hypersurface having a unique singular point $x$ of type $A_1$,
 so that its local model is given by the equation $z_0^2 + \dots + z_n^2 = 0$ in $\C^{n+1}$. Let $(X_\eta)_{\eta > 0}$
 be a smoothing of $X_0$, $L_\eta \subset X_\eta$ be the vanishing cycle of $x$ and $J_\eta$ be the complex
 structure of $X_\eta$. Then, $J_\eta$ lies in the $\eta$-neighborhood of some $A_1$-singular almost-complex structure 
 $J_0 \in {\cal J}_\omega^\infty (L_\eta)$ of $(X_\eta , \omega \vert_{X_\eta})$.
 
 This example justifies the following Definition \ref{defvanishingcycle}.
 
 \begin{defi}
 \label{defvanishingcycle}
 A vanishing cycle of $(X , \omega)$ is a pair $\widetilde{L} = (L , J_L)$ where $L$ is a smooth Lagrangian sphere of $(X , \omega)$
 and $J_L \in {\cal J}_\omega^\infty (L)$.
 \end{defi}
 
 \subsection{Stasheff's associahedron}
 
 \subsubsection{Definition}
 \label{pardefassociahedron}
 
\begin{figure}[ht]
\begin{center}
\includegraphics[scale=0.5]{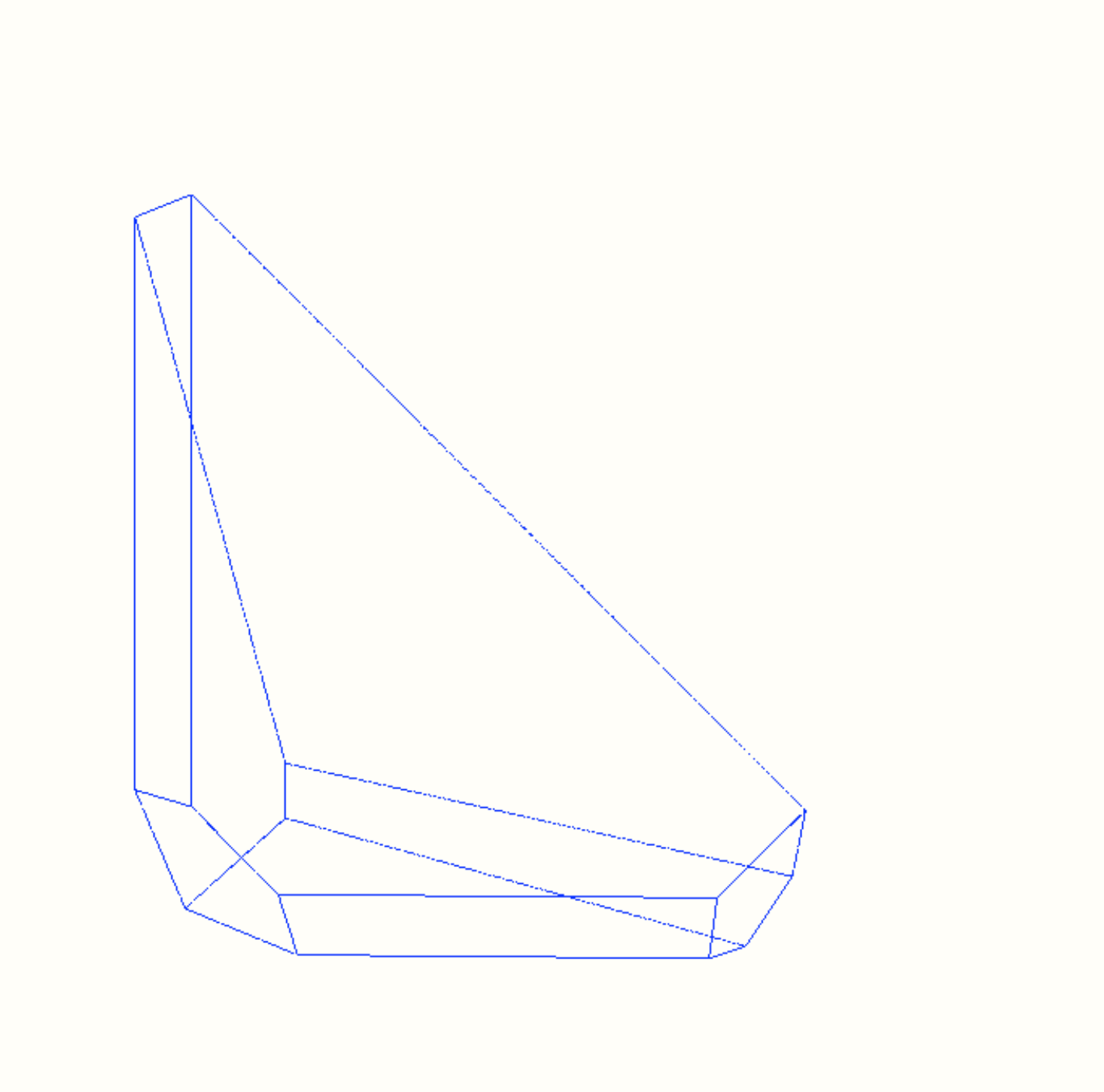}
\end{center}
\caption{The three-dimensional associahedron $K_5$.}
\label{figureassociahedron}
\end{figure}

 Let $D = \{ z \in \C \, \vert \, \vert z \vert \leq 1 \}$ be the complex unit disk equipped with its canonical orientation.
 For every integer $l \geq 2$, denote
 $\stackrel{\circ}{K}_l = \{ (z_0 , \dots , z_l) \in \partial D \, \vert \, \forall 0 \leq i \leq l, \, z_i \in ] z_{i-1} , z_{i+1} [ \}/\text{Aut} (D)$,
 where $z_{-1} = z_l$ and $z_{l+1} = z_0$. Hence, $\stackrel{\circ}{K}_l $ denotes the moduli space of punctured
 holomorphic disks having $l+1$ punctures cyclically ordered on the boundary. For every $\underline{z} \in \stackrel{\circ}{K}_l $,
 let $D_{\underline{z}} = D \setminus \underline{z}$ and for every $0 \leq i \leq l$, denote by $\partial_i D_{\underline{z}}$
 the interval $] z_{i} , z_{i+1} [ \in \partial D$. Let $K_l$ be the stable compactification
 of $\stackrel{\circ}{K}_l $, it has the structure of a $(l-2)$-dimensional convex polytope of the Euclidian space isomorphic
 to  Stasheff's associahedron, see Figure \ref{figureassociahedron}, \cite{Sta}, \cite{Forc}, \cite{FukOh} and references therein.
 We agree that $K_0$ and $K_1$ are points.
 
\begin{figure}[ht]
\begin{center}
\includegraphics{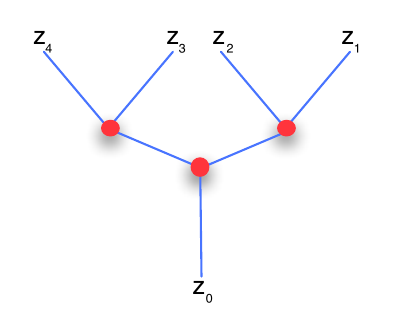}
\end{center}
\caption{Binary tree encoding the bracketing $(43)(21)$.}
\label{figurebracket}
\end{figure}

 This associahedron can also be defined as the space of connected metric binary trees with $l+1$ free
 edges and such that the interior edges have lengths between $0$ and $1$. The vertices encodes 
 bracketings of the ordered set $(l, \dots , 1)$ whereas edges encodes the associativity rule applied to 
 one bracketing, see Figure \ref{figurebracket}.
 
 Let $U_l \to K_l$, $l \geq 2$, be the universal curve, so that its fiber over $\underline{z} \in K_l$
 is isomorphic to the pointed stable disk $D_{\underline{z}}$ encoded by $\underline{z}$.
 For every $\underline{z} \in \partial K_l$ and $0 \leq i \leq l$, denote by $\partial_i D_{\underline{z}}$
 the oriented reducible component of  $\partial D_{\underline{z}}$ going from $z_{i}$ to $z_{i+1}$
 without meeting $\underline{z} \setminus \{ z_{i} , z_{i+1} \}$. This definition extends the one given
 for $\underline{z} \in \stackrel{\circ}{K}_l $. The interior of every $j$-codimensional face of $K_l$,
 $1 \leq j \leq l-2$ encodes a stable disk having $j+1$ irreducible components $D^1 , \dots , D^{j+1}$.
 Such a face is thus canonically isomorphic to a product $K_{l_1} \times \dots \times K_{l_{j+1}}$ of
 associahedra, where $l_1 +  \dots + l_{j+1} = l + j$ and $D^i$ contains $l_i + 1$ punctures, $1 \leq i \leq j+1$.
 Denote by $\sigma_0 , \dots , \sigma_l : K_l \to U_l$ the tautological sections defined by
 $\sigma_i (\underline{z}) = z_i$ for every $\underline{z} \in K_l$ and $0 \leq i \leq l$.
 To every codimension one face $F$ of $K_l$ is associated an additional tautological section $ \sigma_F : F \to U_l$
 which maps $\underline{z} \in \stackrel{\circ}{F}$ to the puncture of $D_{\underline{z}}$ linking its two irreducible components.
 
 \begin{lemma}
 \label{lemmastriplikeend}
 There exists a family $(V_l \subset U_l)_{l \geq 2}$ such that for every $l \geq 2$, $\underline{z} \in K_l$ and
 every irreducible component $\Delta_{\underline{z}}$ of $D_{\underline{z}}$, $(V_l  \cap \Delta_{\underline{z}}) \subset
 \partial \Delta_{\underline{z}}$ is a compact neighborhood of the punctures of $\Delta_{\underline{z}}$. Each connected component
 of this intersection is homeomorphic to an interval having a unique puncture in its interior. Moreover, 
 $V_l  \cap \Delta_{\underline{z}}$ only depends on the conformal structure of $\Delta_{\underline{z}}$ and neither depends
 on the conformal structure of $D_{\underline{z}} \setminus \Delta_{\underline{z}}$, nor on the type of punctures
 of $\Delta_{\underline{z}}$ or on $l \geq 2$.
 \end{lemma}
 
 {\bf Proof:}
 
 The construction of $V_l$ is done by induction on $l \geq 2$. The initial case $l=2$ is no problem,
 the second condition imposes to choose $V_2$ invariant under the automorphisms of the disk
 cyclically permuting the three punctures. Now assume the construction done up to $l-1$ and let us construct
 $V_l$. If $\underline{z} \in \partial K_l$, $D_{\underline{z}}$ has $j+1$ irreducible components $D^1 , \dots , D^{j+1}$,
 $1 \leq j \leq l-2$, and every such component $D^i$ has $l_i +1 < l+1$ punctures, $1 \leq i \leq j+1$. We set
 $V_l \cap D^i_{\underline{z}} = V_{l_i} \cap D^i_{\underline{z}}$ and then extend $V_l \vert_{\partial K_l}$ to a set
 defined over the whole $K_l$ such that it remains invariant under the action of the cyclic group
 permuting cyclically the punctures. $\square$\\
 
 Following \cite{Seidbook}, we set
 \begin{defi}
 \label{defstriplikeend}
 A coherent choice of strip like ends is a system $(V_l \subset U_l)_{l \geq 2}$ given by Lemma \ref{lemmastriplikeend}.
 \end{defi}
 This Definition \ref{defstriplikeend} gets justified by the following fact (compare Lemma $9.3$ of \cite{Seidbook}).
 For every $\underline{z} \in K_l$, every irreducible component $\Delta_{\underline{z}}$ of $D_{\underline{z}}$
 and every puncture $w$ of $\Delta_{\underline{z}}$,  $V_l  \cap \Delta_{\underline{z}}$ provides a neighborhood
 $[w_- , w_+] \subset \partial \Delta_{\underline{z}}$ of $w$ in $\partial \Delta_{\underline{z}}$. 
 There exists a unique injective holomorphic map $\psi_w : \R_+ \times [-1 , 1] \to
 \Delta_{\underline{z}}$ such that $\psi_w (0, -1) = w_-$, $\psi_w (0, 1) = w_+$ and $\lim_{+ \infty} \psi_w = w$.
 These strip like ends are disjoint to each other. We agree that $U_0$ is a point and $U_1 = [-1 , 1]$.
 
 Finally, we denote by $\overline{\cal J}_\omega (U_l , V_l)$, $l \in \N$, the space of maps $J : U_l \to \overline{\cal J}_\omega$
 having the same regularity as the one chosen for our almost-complex structures and which for $l \geq 2$ satisfy the
 following three properties $P_1$, $P_2$, $P_3$.\\
 
 $P_1$: For every $0 \leq i \leq l$, $J(\partial_i D_{\underline{z}})$ is a point $J_i \in \overline{\cal J}_\omega$ which does not
 depend on $\underline{z} \in K_l$.
 
 $P_2$: For every  tautological section $ \sigma : F \to U_l$, every $\underline{z} \in {F}$ and every associated strip like end
 $\psi_{\sigma (\underline{z})}$ given by $V_l$, the composition $J \circ \psi_{\sigma (\underline{z})} :  \R_+ \times [-1 , 1] \to
 \overline{\cal J}_\omega$ does not depend on the first factor $\R_+$. Moreover, the associated path
 $J_\sigma = J \circ \psi_{\sigma (\underline{z})} \vert_{\{ 0 \} \times  [-1 , 1]}$ neither depends on $\underline{z} \in {F}$ nor on
 the choice of the strip like end $\psi_{\sigma (\underline{z})}$. Note indeed that when $F$ has positive codimension, 
 $\sigma$ encodes a puncture linking two irreducible components, both having an associated strip like end.
 
 $P_3$: For every $\underline{z} \in \partial K_l$ and every irreducible component $\Delta_{\underline{z}}$ of $D_{\underline{z}}$,
 the restriction of $J$ to $\Delta_{\underline{z}}$ does not depend on the conformal structure of $D_{\underline{z}} \setminus
 \Delta_{\underline{z}}$.\\
 
 Property $P_3$ is a compatibility property of $J$ with respect to the product structure of the faces of $K_l$.
 
 \subsubsection{Orientation}
 \label{subsubsectionorientationassociahedron}
 
 Let $\underline{z} = (z_0 , \dots , z_l) \in K_l$, $l \geq 2$. Let us identify the unit disk $D$ with the upper half
 complex plane $\H = \{Êz \in \C P^1 \, \vert \, \text{Im} (z) \geq 0 \}$, so that $z_0$ is the point at infinity and
 $z_1 < \dots < z_l \in \R$. For every $1 \leq i \leq l$, we equip the tangent line $T_{z_i} \R$ with its induced orientation.
 The group $\text{Aut} (\H ; z_0)$ of automorphisms of $\H$ fixing $z_0$ is generated by translations by real numbers
 and dilation based at any real points. We equip its Lie algebra $\text{aut} (\H ; z_0) = T_{id} \text{Aut} (\H ; z_0)$
 with the orientation $T^+ \wedge H^+$, where $T^+$ is an infinitesimal translation by positive real numbers and
 $H^+$ an infinitesimal dilation. Note that the interior of $\H$ has the structure of a principal homogenous space
 over $\text{Aut} (\H ; z_0)$ and that our orientation is nothing but the one induced from it. From the
 short exact sequence
 $$0 \to \text{aut} (\H ; z_0) \to T_{z_1} \R \times \dots \times T_{z_l} \R \to T_{\underline{z}} K_l \to 0,$$
 we deduce an orientation on $T_{\underline{z}} K_l $.
 
 Another equivalent definition of this orientation is the following. Let $1 \leq i \leq l$ and let us
  identify this times $D \setminus \{ z_0 , z_i \}$ with the strip $\Theta_i =  \R \times [-1 , 1]$ so that
  $z_0 = -\infty$,  $z_i = + \infty$, $z_1 < \dots < z_{i-1} \in \R \times \{ 0 \}$ and $z_l < \dots < z_{i+1} \in \R \times \{ 1 \}$.
  The group $\text{Aut} (\Theta_i)$ of automorphisms of the strip $\Theta_i$  is reduced to translations by real numbers,
  we equip its Lie algebra $\text{aut} (\Theta_i)$ with the orientation given by an infinitesimal translation $T^-$ by negative
  real numbers, it flows from $z_i$ to $z_0$. This infinitesimal translation $T^-$ induces an orientation on all the tangent
  lines $T_{z_j} \partial \Theta_i$, $1 \leq j \leq l$, $j \neq i$. From the
 short exact sequence
 $$0 \to \text{aut} (\Theta_i) \to T_{z_1} \partial \Theta_i \times \dots  \times \widehat{T_{z_i} \partial \Theta_i} \times \dots 
 \times T_{z_l} \partial \Theta_i \to T_{\underline{z}} K_l \to 0,$$
 we deduce an orientation on $T_{\underline{z}} K_l $. This new orientation coincides with the preceding one. In particular, 
 it does not depend on the choice of $1 \leq i \leq l$.
 
 Every face of the associahedron $K_l$ is canonically a product of lower dimensional associahedra. In particular,
 faces of codimension one are products of two associahedra and thus now inherit two different orientations,
 one from $K_l$ and one from the products of lower dimensional associahedra. The following Lemma \ref{lemmaassociahedron}
 compare these two orientations.
 
 \begin{lemma}
 \label{lemmaassociahedron}
 Let $F = \{ (\underline{z} , \underline{w}) \in K_{l_1} \times K_{l_2} \, \vert \, z_i = w_0 \}$ be a codimension one face
 of $K_l$, where $l, l_1, l_2 \geq 2$, $1 \leq i \leq l_1$ and $l_1 + l_2 = l+1$. Then, the orientations of $F$ induced
 by $\partial K_l$ and $K_{l_1} \times K_{l_2} $ coincide if and only if $l_1 l _2 + i(l_2 - 1)$ is odd.
 \end{lemma}

 {\bf Proof:}
 
 To begin with, assume that $i>1$. We represent $\underline{z} \in K_{l_1}$ by the strip $\Theta_{\underline{z}} =  \R \times [-1 , 1]$,
 where $z_1 < \dots < z_{i-1} \in \R \times \{ 0 \}$, $z_{l_1} < \dots < z_{i+1} \in \R \times \{ 1 \}$ and $z_{i}= +\infty$.
 We represent $\underline{w} \in K_{l_2}$ by the strip $\Theta_{\underline{w}} =  \R \times [-1 , 1]$,
 where $w_1 < \dots < w_{l_2-1} \in \R \times \{ 0 \}$ and $w_{l_2} = +\infty$. For every $\rho \gg 0$, the glueing
 $\Theta_{\underline{z}} \star_\rho \Theta_{\underline{w}}$ is encoded by an interior point of $K_l$. It is a strip
 $\Theta_{\underline{z}  \star \underline{w}}$ where $z_1 - \rho < \dots < z_{i-1} - \rho < w_1+ \rho < \dots < 
 w_{l_2-1} + \rho \in \R \times \{ 0 \}$ and $z_{l_1} - \rho < \dots < z_{i+1} - \rho \in \R \times \{ 1 \}$. The vector
 $\frac{\partial}{\partial \rho} (\Theta_{\underline{z}} \star_\rho \Theta_{\underline{w}}) = - \frac{\partial}{\partial \underline{z}}
 +  \frac{\partial}{\partial \underline{w}}$ is identified with the outward normal vector of $K_l$ at $F$. Let us identify
 $T_{\Theta_{\underline{z}} \star_\rho \Theta_{\underline{w}}} K_l$ with the complement of the orbit of 
 $\text{aut} (\Theta_{\underline{z}} \star_\rho \Theta_{\underline{w}})$ in 
 $T_{z_1} \R \times \dots \times T_{z_{i-1}} \R \times T_{\underline{w}} \R \times T_{z_{i+1}} \R \times \dots \times T_{z_l} \R$
 for which the first vector $\frac{\partial}{\partial z_1}$ vanishes. The orientation of $T_{\Theta_{\underline{z}} \star_\rho \Theta_{\underline{w}}} K_l$ writes:
 
$$(-\frac{\partial}{\partial z_2}) \wedge \dots \wedge (-\frac{\partial}{\partial z_{i-1}}) \wedge (-\frac{\partial}{\partial w_1}) \wedge \dots 
\wedge (-\frac{\partial}{\partial w_{l_2-1}}) \wedge (-\frac{\partial}{\partial z_{i+1}}) \wedge \dots \wedge (-\frac{\partial}{\partial z_{l_1}}) $$
 $$= (-1)^{i-1} (\frac{\partial}{\partial \underline{w}}- \frac{\partial}{\partial \underline{z}}) \wedge (-\frac{\partial}{\partial z_2}) \wedge \dots \wedge (-\frac{\partial}{\partial z_{i-1}}) \wedge (-\frac{\partial}{\partial w_2}) \wedge \dots 
 \wedge (-\frac{\partial}{\partial w_{l_2-1}}) \wedge (-\frac{\partial}{\partial z_{i+1}}) \wedge \dots \wedge (-\frac{\partial}{\partial z_{l_1}}) $$
 $$ =(-1)^{i-1 + l_2(l_1 - i)} (\frac{\partial}{\partial \underline{w}}- \frac{\partial}{\partial \underline{z}}) \wedge (-\frac{\partial}{\partial z_2}) \wedge \dots \wedge \widehat{(-\frac{\partial}{\partial z_{i}})} \wedge \dots \wedge (-\frac{\partial}{\partial z_{l_1}}) \wedge (-\frac{\partial}{\partial w_2}) \wedge \dots \wedge (-\frac{\partial}{\partial w_{l_2-1}})$$
 
 Since $(-\frac{\partial}{\partial z_2}) \wedge \dots \wedge \widehat{(-\frac{\partial}{\partial z_{i}})} \wedge \dots \wedge (-\frac{\partial}{\partial z_{l_1}}) \wedge (-\frac{\partial}{\partial w_2}) \wedge \dots \wedge (-\frac{\partial}{\partial w_{l_2-1}})$ form a direct basis of
 $T_{(\underline{z} , \underline{w})} (K_{l_1} \times K_{l_2})$, we deduce the result when $i>1$. If $i=1$, the orientation of $T_{\Theta_{\underline{z}} \star_\rho \Theta_{\underline{w}}} K_l$ writes:

$\begin{array}{cl}
 & (-\frac{\partial}{\partial w_2}) \wedge \dots 
 \wedge (-\frac{\partial}{\partial w_{l_2-1}}) \wedge (-\frac{\partial}{\partial z_2}) \wedge \dots \wedge (-\frac{\partial}{\partial z_{l_1}}) \\
 =& (-1)^{ l_2(l_1 - 1)} (-\frac{\partial}{\partial z_2}) \wedge \dots \wedge (-\frac{\partial}{\partial z_{l_1}}) \wedge (-\frac{\partial}{\partial w_2}) \wedge \dots \wedge (-\frac{\partial}{\partial w_{l_2-1}})  \\
 =& (-1)^{l_2(l_1 - 1)} (\frac{\partial}{\partial \underline{w}}- \frac{\partial}{\partial \underline{z}}) \wedge (-\frac{\partial}{\partial z_3}) \wedge \dots \wedge (-\frac{\partial}{\partial z_{l_1}}) \wedge (-\frac{\partial}{\partial w_2}) \wedge \dots \wedge (-\frac{\partial}{\partial w_{l_2-1}}).
 \, \square 
 \end{array}$ \\
 
 \begin{rem}
 \label{remdbarz}
 Let $\underline{z} \in K_l$ and $\overline{\partial}_{\underline{z}} : L^{k,p} (D_{\underline{z}} ; T D_{\underline{z}} , 
 T \partial D_{\underline{z}} ) \to L^{k-1,p} (D_{\underline{z}} ; \Lambda^{0,1} D_{\underline{z}} \otimes T D_{\underline{z}})$
 be the associated Cauchy-Riemann operator. When $l \geq 2$, $\overline{\partial}_{\underline{z}} $ is injective and
 its cokernel is canonically isomorphic to $T_{\underline{z}} K_l$. Hence, the orientation fixed on $K_l$ induces
 an orientation on the real line $\det (\overline{\partial}_{\underline{z}})$. When $l=1$,  $\overline{\partial}_{\underline{z}} $ is
 surjective and its kernel is one-dimensional. We agree that an element of the kernel which flows to $z_0$ form
 a direct basis of $\det (\overline{\partial}_{\underline{z}})$.
 \end{rem}

\subsection{Lagrangian conductors}
\label{subseclagrangianconductors}

Denote by ${\cal L}ag^\pm$ the set of pairs $(L, \p^\pm_L)$ such that $L$ is a closed Lagrangian submanifold embedded in
$X$ and $ \p^\pm_L$ is a $\widetilde{GL}_n^\pm (\R)$-structure on $L$.

\begin{defi}
\label{defelemconductor}
An elementary Lagrangian conductor is a pair $(L,J) \in {\cal L}ag^\pm \times \overline{\cal J}_\omega$ such that $J$ is
non-singular along $L$ and for every $\mu \leq 2$, the space of $J$-holomorphic disks of Maslov index $\mu$ with 
boundary on $L$ is of effective dimension strictly less than the expected one $\mu +n - 3$.
\end{defi}
When $J \in \partial \overline{\cal J}_\omega$, by $J$-holomorphic disk we mean split $J$-holomorphic disk in the
sense of \cite{EGH}, see \S \ref{parc1}. In particular, as soon as $\mu \leq 3-n$, the space of $J$-holomorphic disks of Maslov index $\mu$
has to be empty. Monotone Lagrangian submanifolds in the sense of \cite{Oh1} become Lagrangian conductors once equipped
with any $J \in  \overline{\cal J}_\omega$. To these examples we will add in \S \ref{parc1} vanishing cycles in the sense of Definition 
\ref{defvanishingcycle} when $(X , \omega)$ has vanishing first Chern class.

\begin{defi}
\label{defconductor}
A Lagrangian conductor $\tilde{L}$ of  $(X , \omega)$ is an uple $(L_0 , \dots , L_l ; J^{\tilde{L}})$, $l \in \N$, such that:

1) $J^{\tilde{L}} \in \overline{\cal J}_\omega (U_l , V^{\tilde{L}}_l)$, where $V^{\tilde{L}}_l$ is a coherent choice of strip like
ends given by Definition \ref{defstriplikeend}.

2) $\forall 0 \leq i \leq l$, $(L_i , J^{\tilde{L}}_i)$ is an elementary Lagrangian conductor, where $J^{\tilde{L}}_i$ is the almost-complex
structure given by Property $P_1$ of \S \ref{pardefassociahedron}.
\end{defi}

\begin{figure}[ht]
\begin{center}
\includegraphics[scale=0.75]{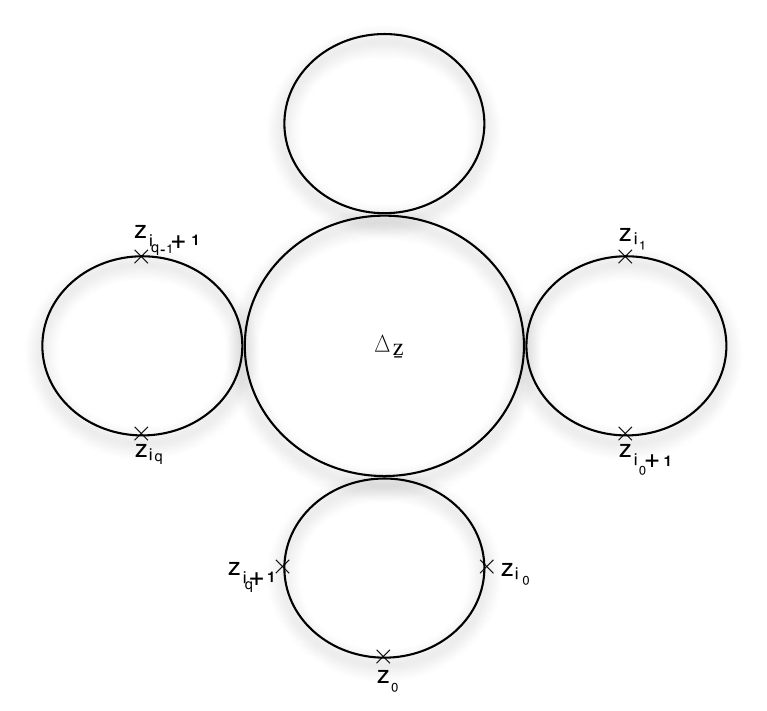}
\end{center}
\caption{Stable disk $D_{\underline{z}} $, $\underline{z} \in F_{i_0  \dots  i_q}$.}
\label{figurefacef}
\end{figure}

Let $2 \leq q \leq l$ and $0 \leq i_0 < \dots < i_q \leq l$. Denote by $F_{i_0  \dots  i_q}$ the maximal face of $K_l$  having
the property that for every $\underline{z} \in F_{i_0  \dots  i_q}$, $D_{\underline{z}} $ contains a sub stable disk 
whose boundary is disjoint from $\cup_{i \notin \{i_0 , \dots ,  i_q \}} \partial_i D_{\underline{z}} $. Denote by
$\Delta_{\underline{z}} $ the maximal sub stable disk of $D_{\underline{z}} $  having
this property. The face $F_{i_0  \dots  i_q}$ canonically decomposes as a product $K_q \times K_{i_1 - i_0} \times
\dots \times K_{i_q - i_{q-1}} \times K_{l - i_q + i_0  +1}$ where the first factor $K_q$ encodes the component
$\Delta_{\underline{z}} $, see Figure \ref{figurefacef}.

Denote by $U_{i_0  \dots  i_q} = \{ \Delta_{\underline{z}} \, \vert \, \underline{z} \in F_{i_0  \dots  i_q} \}$. The projection
$F_{i_0  \dots  i_q} \to K_q$ onto the first factor lifts to a projection $U_{i_0  \dots  i_q} \to U_q$. If now
$J \in \overline{\cal J}_\omega (U_l , V_l)$, Property $P_3$ satisfied by $J$ ensures that $J$ induces on the quotient
a map $J_{i_0  \dots  i_q}  \in \overline{\cal J}_\omega (U_q , V_l \cap U_{i_0  \dots  i_q})$. Properties $P_1$ and
$P_2$ satisfied by $J$ make it possible to extend this construction to $q=0$ and $q=1$ respectively. 
As a result, we get the following Definition \ref{defsubconductor}.

\begin{defi}
\label{defsubconductor}
Let $\tilde{L} = (L_0 , \dots , L_l ; J^{\tilde{L}})$ be a Lagrangian conductor of $(X , \omega)$. A subconductor of
$\tilde{L} $ is a  Lagrangian conductor of the form $\tilde{L}_{i_0  \dots  i_q} = (L_{i_0} , \dots , L_{i_q} ; J^{\tilde{L}}_{i_0  \dots  i_q} )$, where
$0 \leq i_0 < \dots < i_q \leq l$. We then say that $\tilde{L} $ refines $\tilde{L}_{i_0  \dots  i_q}$.
\end{defi}

The particular refinements given by Definition \ref{defextensionlagconductor} play a special r\^ole.

\begin{defi}
\label{defextensionlagconductor}
For every Lagrangian conductor $\tilde{L} = (L_0 , \dots , L_l ; J^{\tilde{L}})$ and every $0 \leq q < l$, we say that
$\tilde{L} $ is a refinement of $\tilde{L}_{0  \dots  q}$ by $\tilde{L}_{q+1 \dots  l}$.
\end{defi}

If $\tilde{L}_1 $ and $\tilde{L}_2 $ are two Lagrangian conductors, there exists a refinement $\tilde{L} $ of $\tilde{L}_1 $ by $\tilde{L}_2 $
if and only if the Lagrangian submanifolds in $\tilde{L}_1 $, $\tilde{L}_2 $ are transversal to each other. We denote by
$\text{Ob} ( {\cal CL}^\pm (X, \omega))$ the set of Lagrangian conductors of the manifold $(X , \omega)$.

\begin{defi}
\label{defeffectivecontinuation}
An effective continuation $H$ from the Lagrangian conductor $\tilde{L} = (L_0 , \dots , L_q ; J^{\tilde{L}})$ to 
the Lagrangian conductor $\tilde{L}' = (L'_0 , \dots , L'_l ; J^{\tilde{L}'})$ is the data of a subset
$I_H \subset \{ 0 , \dots , q \}$, an increasing injection $\phi_H : j \in I_H \to i_j \in \{ 0 , \dots , l \}$ and 
 a path $(L_j^s , J_j^s)_{s \in [0,1]}$ of elementary Lagrangian conductors, for all $ j \in I_H$, such that
 
 1) $(L_j^0 , J_j^0) = (L_j , J_j^{\tilde{L}})$ and $(L_j^1 , J_j^1) = (L'_{i_j} , J_{i_j}^{\tilde{L}'})$.
 
 2) $(L_j^s)_{s \in [0,1]}$ is a family of Hamiltonian isotopic elements of ${\cal L}ag^\pm$ .
\end{defi}

Hence, an effective continuation is given by Hamiltonian isotopies between some Lagrangian submanifolds
of $\tilde{L} $ and $\tilde{L}' $ together with  homotopies between the $\# I_H$ corresponding values of
$J^{\tilde{L}}$ and $J^{\tilde{L}'}$ given by the Property $P_1$ satisfied by these functions. The latter condition
will be essential in \S \ref{parc1}, see Remark \ref{remmonodromy}.

\begin{defi}
\label{defcontinuation}
A continuation between Lagrangian conductors is a homotopy class with fixed extremities of
effective continuations. 
\end{defi}

Let $H : \tilde{L} \to \tilde{L}'$ be a continuation given by Definitions \ref{defeffectivecontinuation}
 and \ref{defcontinuation}. The sub conductor $ \tilde{L}'_{\text{Im} (\phi_H)}$ of  $\tilde{L}'$ is called the image
 of $H$ whereas the subconductor$ \tilde{L}_{I_H}$ is called the cokernel of $H$. A sequence
 $\tilde{L} \stackrel{H}{\to} \tilde{L}' \stackrel{K}{\to} \tilde{L}''$ is called exact if and only if
 $\# (\phi_H (I_H) \cap I_K) \leq 1$. We denote by
$\text{Hom} ( {\cal CL}^\pm (X, \omega))$ the set of continuations between Lagrangian conductors of the manifold $(X , \omega)$.
Denote by ${\cal CL}^\pm (X, \omega)$ the pair $\big( \text{Ob} ( {\cal CL}^\pm (X, \omega)) , \text{Hom} ( {\cal CL}^\pm (X, \omega)) \big)$.

\begin{prop}
\label{propcl}
Let $(X, \omega)$ be a symplectic manifold of dimension at least four. Then, ${\cal CL}^\pm (X, \omega) $
has the structure of a small category with the properties of exact sequence, sub object, refinement, cokernel and image.
\end{prop}

{\bf Proof:}

Objects and morphisms of this pair are sets. The composition of morphisms is given by the composition of injections
given by Definition  \ref{defeffectivecontinuation}. It is associative and has identities for every object.
Properties of sub objects and refinements are given by Definitions \ref{defsubconductor}, \ref{defextensionlagconductor}. 
Properties of cokernel and image have just been defined. $\square$

\begin{rem}
\label{remrefinement}
The property of refinement given by Definition \ref{defextensionlagconductor} provides another structure of category
on ${\cal CL}^\pm (X, \omega)$ with same set of objects but with morphisms between objects given by homotopy classes
of refinements of the first by the second one. Property $P_3$ given in \S \ref{pardefassociahedron} ensures the associativity
of the composition of such morphisms.
\end{rem}
 
\section{Floer functor}

In this paragraph, $(X, \omega)$ stands for a semipositive symplectic manifold of dimension $2n \geq 4$ which is either closed
or convex at infinity.

\subsection{Floer complex}

Let $\tilde{L} = (L_0 , \dots , L_l ; J^{\tilde{L}}) \in \text{Ob} ( {\cal CL}^\pm (X, \omega)) $, we are going to associate
to this Lagrangian conductor a complex ${\cal F} (\tilde{L}) \in \text{Ob} ( K^b ({\cal OS}_1 (X, \omega)))$.
If $l=0$, ${\cal F} (\tilde{L}) = 0$. Otherwise, for every $0 \leq i < j \leq l$, there is a tautological injection
$x \in L_i \cap L_j \mapsto \lambda_x \in \text{Ob} ( {\cal OS}_1 (X, \omega))$. We set
$CF (L_i , L_j) = \oplus_{x \in L_i \cap L_j} \lambda_x$ and then for every $1 \leq q \leq l$,
$$CF_q (\tilde{L}) = \oplus_{0 \leq i_0 < \dots < i_q \leq l} CF (L_{i_0} , L_{i_1}) \otimes \dots \otimes CF (L_{i_{q-1}} , L_{i_q}).$$
Finally, we set $CF (\tilde{L}) = \oplus_{q=1}^l CF_q (\tilde{L})$.

Let $\lambda^+ \in CF (L_{i_0} , L_{i_1}) \otimes \dots \otimes CF (L_{i_{q-1}} , L_{i_q})$ and $\lambda^- \in CF (L_{i_0} , L_{i_q})$ be such
that $\mu (\lambda^-) - 1 = \mu(\lambda^+) - q + 1$.
Denote by $\gamma_{\lambda^+ \lambda^-}$ the sum of elementary trajectories $(u, D_{\underline{z}} , \lambda_{\partial u} , \overline{\partial}_u )$ 
from $\lambda^+$ to $\lambda^-$ which are such that:

1) $D_{\underline{z}}$ has one negative puncture $z_0$ and $q$ positive punctures $z_1, \dots z_q$.

2) For every $0 \leq j \leq q$, $u (\partial_j D_{\underline{z}}) \subset L_{i_j}$ and $\lambda_{\partial u} (\partial_j D_{\underline{z}})
\subset TL_{i_j}$.

3) $u : D_{\underline{z}} \to X$ satisfies the Cauchy-Riemann equation $J^{\tilde{L}} \vert_{D_{\underline{z}}} \circ du = du \circ J_{D_{\underline{z}}}$
and $\overline{\partial}_u$ is the associated Cauchy-Riemann operator. Its orientation is induced from the following exact sequence.

$$
0  \to  L^{k,p} (D_{\underline{z}} ; TD_{\underline{z}} , T \partial D_{\underline{z}})  \stackrel{du}{\to} 
L^{k,p} (D_{\underline{z}} ; u^* TX , \lambda_{\partial u})  \to  
L^{k,p} (D_{\underline{z}} ; u^* TX , \lambda_{\partial u})/\text{Im}(du) \to  0 $$
$$\downarrow \overline{\partial}_{\underline{z}} \hspace{4cm} \downarrow \overline{\partial}_u \hspace{4cm} \downarrow \overline{\partial}_N$$
$$0 \to L^{k-1,p} (D_{\underline{z}}  , \Lambda^{0,1} D_{\underline{z}}  ) \stackrel{du}{\to}  
L^{k-1,p} (D_{\underline{z}}  , \Lambda^{0,1} D_{\underline{z}}  \otimes u^* TX)  \to 
L^{k-1,p} (D_{\underline{z}}  , \Lambda^{0,1} D_{\underline{z}}  \otimes u^* TX)/\text{Im}(du)  \to  0 $$

Indeed, the Fredholm index of $\overline{\partial}_u$ coincides with the index of $\gamma_{\lambda^+ \lambda^-}$,
the Fredholm index of $\overline{\partial}_{\underline{z}} $ is $2-q$ and the genericness of $J^{\tilde{L}}$  ensures
that $\overline{\partial}_N$ is an isomorphism, see Proposition \ref{propindex}. The latter is thus canonically oriented,
$\overline{\partial}_{\underline{z}} $ is oriented from Remark \ref{remdbarz} so that $\overline{\partial}_u$ gets an orientation.

\begin{prop}
\label{propindex}
Let $\lambda^+ \in CF (L_{i_0} , L_{i_1}) \otimes \dots \otimes CF (L_{i_{q-1}} , L_{i_q})$ and $\lambda^- \in CF (L_{i_0} , L_{i_q})$ be such
that $\mu (\lambda^-) - 1 = \mu(\lambda^+) - q + 1$, where $0 \leq i_0 < \dots < i_q \leq l$, $q \geq 1$. Then, $\gamma_{\lambda^+ \lambda^-}$
is a trajectory in the sense of Definition \ref{deftrajectory}.
\end{prop}

{\bf Proof:}

Let ${\cal M}(\lambda^+ ,\lambda^-)$ be the space of triples $(u, D_{\underline{z}} , J)$ from $\lambda^+$ to $\lambda^-$ such that 
$\underline{z} \in \stackrel{\circ}{K}_l$,
$J \in \overline{\cal J}_\omega (U_q , V_q)$ and $u (\partial_j D_{\underline{z}}) \subset L_{i_j}$, $0 \leq j \leq q$.
This space ${\cal M}(\lambda^+ ,\lambda^-)$ is a separable Banach manifold whose regularity is the difference
between the regularity chosen for our almost-complex structures and the regularity $k$ of our maps $u$. This follows from
a standard argument which we do not reproduce here, see Proposition $3.2.1$ of \cite{McDSal}, \cite{Oh1} and \cite{Floer}.
When $q=1$, this argument is however slightly more elaborate and given in \cite{Oh3}. Moreover, the index of the projection
$(u, D_{\underline{z}} , J) \in {\cal M}(\lambda^+ ,\lambda^-) \mapsto (\underline{z} , J) \in K_q \times \overline{\cal J}_\omega (U_q , V_q)$
is the same as the Fredholm index of the operator $\overline{\partial}_N$ appeared in the previous exact sequence.
The index of the operator $\overline{\partial}_u$ may be computed as follows. The construction of the double applied to
$D_{\underline{z}} $ provides a symplectic vector bundle $u^* TX$ over $\C P^1 \setminus \{ \underline{z} \}$ equipped with
an antisymplectic involution whose fixed point set is a Lagrangian sub bundle over $\R P^1 \setminus \{ \underline{z} \}$, see
for example \cite{HLS}. The choice of an extension of $\lambda_{\partial u}$ over the whole $D_{\underline{z}}$ given by the third
property of elementary trajectories provides an extension of the
 Lagrangian sub bundle of $u^* TX$ over $\C P^1 \setminus \{ \underline{z} \}$ and thus a trivialization of $u^* TX$. Denote by
 $\overline{\partial}_{\C P^1}$ the $\Z/2\Z$-equivariant Cauchy-Riemann operator associated to $u^* TX$, so that
 $\ind_\R (\overline{\partial}_{\C P^1}) = 2\ind_\R (\overline{\partial}_u)$, see for example \cite{WelsInvent}. The index of
 $\overline{\partial}_{\C P^1}$ was computed in the thesis \cite{Schw} or in \cite{HWZ1} and is given by the formula
 $\ind_\R (\overline{\partial}_{\C P^1}) = \sum_{j=0}^q \mu_{CZ} (u(z_j)) + n(1-q)$, where $\mu_{CZ} (u(z_j))$ denotes the
 Conley-Zehnder index of the identity of $T_{z_j} X$ computed in our trivialization. Thus, $\mu_{CZ} (u(z_j))$ is twice the Robbin-Salamon
 index $\mu(l_{-1}^j , l^j)$. We deduce that $\ind_\R (\overline{\partial}_u) = \frac{1}{2} \big( \sum_{j=0}^q (2\mu_{CZ} (u(z_j)) - n) + 2n\big)=
 \mu (\lambda^-) - \mu(\lambda^+) =2-q$, so that the index of $\overline{\partial}_N$ vanishes. Finally, the fact that
 $\gamma_{\lambda^+ \lambda^-}$ satisfies the Novikov condition of Definition \ref{deftrajectory} follows from Gromov-Floer compactness
 Theorem  since for generic $J$ there is no bubbled punctured $J$-holomorphic disk from $\lambda^+$ to $\lambda^-$, see
 Proposition $4.1$ of \cite{Oh1}. Note that here the condition on disks of Maslov indices less than three given in Definition
 \ref{defelemconductor} plays a crucial r\^ole. $\square$\\

Following \cite{Fuk2}, denote by $m_q : CF (L_{i_0} , L_{i_1}) \otimes \dots \otimes CF (L_{i_{q-1}} , L_{i_q}) \to CF (L_{i_0} , L_{i_q})$
the sum of trajectories $\gamma_{\lambda^+ \lambda^-}$ of index $2-q$ going from an open string
$\lambda^+ \in CF (L_{i_0} , L_{i_1}) \otimes \dots \otimes CF (L_{i_{q-1}} , L_{i_q})$ to an open string $\lambda^- \in CF (L_{i_0} , L_{i_q})$.
Let then $\delta^{CF} = \oplus_{q=1}^l  \delta^{CF}_q : CF (\tilde{L}) \to CF (\tilde{L})$ be the morphism
defined by
$\delta^{CF}_q = \oplus_{l_2=1}^q (-1)^{ql_2} \oplus_{i=1}^{q_1} (-1)^{i(l_2 - 1)} id_{i-1} \otimes m_{l_2} \otimes id_{q_1 - i} \, ,$
where $q_1 + l_2 - 1 = q$, $1 \leq q \leq l$ and $\delta^{CF}_q : CF_q (\tilde{L}) \to CF (\tilde{L})$.
In other words, $\delta^{CF}$ is the morphism which satisfies axioms $A$ of \S \ref{subsectionchain}
and whose dual takes the opposite values on elementary strings as the coproducts $m_q^*$, $1 \leq q \leq l$.

\begin{theo}
\label{theocomplex}
Let $(X, \omega)$ be a semipositive symplectic manifold of dimension $2n \geq 4$ which is either closed
or convex at infinity and $\tilde{L} = (L_0 , \dots , L_l ; J^{\tilde{L}}) \in \text{Ob} ( {\cal CL}^\pm (X, \omega)) $, $l \in \N$.
Then, $\delta^{CF} \circ \delta^{CF} = 0$, so that $(CF (\tilde{L}) , \delta^{CF}) \in \text{Ob} ( K^b ({\cal OS}_1 (X, \omega)))$.
\end{theo}

{\bf Proof:}

\begin{figure}[ht]
\begin{center}
\includegraphics[scale=0.75]{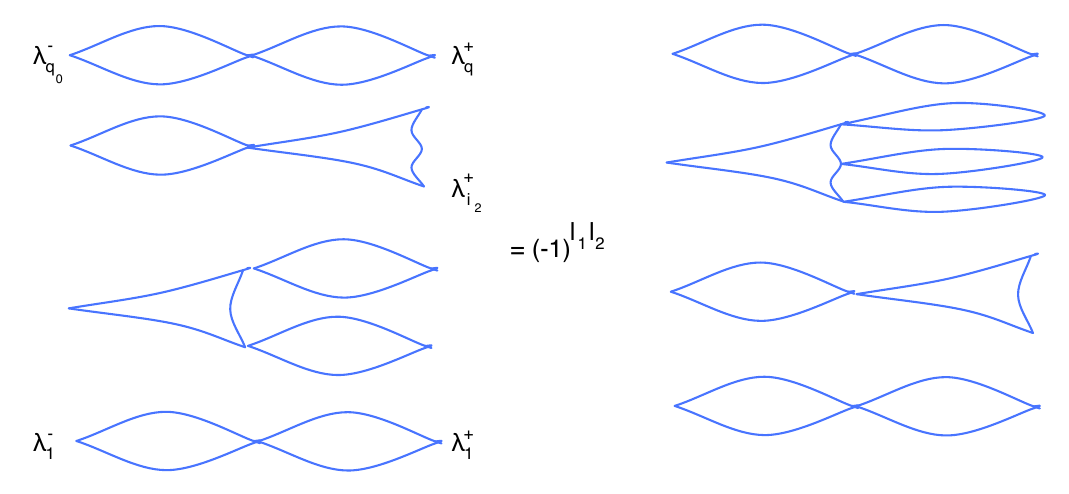}
\end{center}
\caption{$l_2 = 3$, $l_1 = 2$, $q=7$, $q_1 = 5$, $q_0 = 4$, $i_2 = 4$ and $i_1 = 2$.}
\label{figurecommute}
\end{figure}

Let $\lambda^+$ and $\lambda^-$ in $CF (\tilde{L})$ be such that $\mu (\lambda^-) - 1 = \mu(\lambda^+) - q + 2$, so that
the space ${\cal M}(\lambda^+ ,\lambda^- ; J^{\tilde{L}})$ of pairs $(u, D_{\underline{z}} )$ from $\lambda^+$ to $\lambda^-$ such that 
$\underline{z} \in \stackrel{\circ}{K}_l$,  $u (\partial_j D_{\underline{z}}) \subset L_{i_j}$, $0 \leq j \leq q$ and 
$J^{\tilde{L}} \vert_{D_{\underline{z}}} \circ du = du \circ J_{D_{\underline{z}}}$ is one-dimensional. 
It is compact from the condition on disks of Maslov indices less than three given in Definition
 \ref{defelemconductor}, see Proposition $4.3$ of \cite{Oh1}. From Floer glueing
 Theorem \cite{Floer}, \cite{BirCor}, the union of elementary trajectories from $\lambda^+$ to $\lambda^-$ counted by $\delta^{CF} \circ \delta^{CF}$ 
 is in bijection with the boundary of this space. It suffices thus to prove that every such trajectory
 induces the outward normal orientation on ${\cal M}(\lambda^+ ,\lambda^- ; J^{\tilde{L}})$. Two cases have to be considered,
 depending on whether one of the inequalities $i_1 + l_1 \leq i_2$ or $i_2 \leq i_1 - 1$ holds or not. If one of these
 inequalities holds, say the first one, then we have the commuting relation
 $(id_{i_1-1} \otimes m_{l_1} \otimes id_{q_0 - i_1}) \circ (id_{i_2-1} \otimes m_{l_2} \otimes id_{q_1 - i_2})
 = (-1)^{l_1 l_2}  (id_{i_2-l_1} \otimes m_{l_2} \otimes id_{q_1 - i_2})  \circ (id_{i_1-1} \otimes m_{l_1} \otimes id_{q_0 - i_1}) $,
 see Figure \ref{figurecommute}.
 But $\delta^{CF} \circ \delta^{CF}$  counts the left hand side with respect to the sign
 $(-1)^{q_1 l_1 + i_1 (l_1 - 1)} (-1)^{q l_2 + i_2 (l_2 - 1)} $ whereas it counts the operator $(id_{i_2-l_1} \otimes m_{l_2} \otimes id_{q_1 - i_2})  \circ (id_{i_1-1} \otimes m_{l_1} \otimes id_{q_0 - i_1}) $ with respect to the sign
 $$(q+1-l_1)l_2  + (i_2+1-l_1)(l_2-1)+q l_1+i_1(l_1 - 1) = l_1 l_2 -1 + q_1 l_1 + i_1 (l_1 - 1) + q l_2 + i_2 (l_2 - 1) \mod(2),$$
 since $q = q_1 + (l_2 - 1)$. Hence, the contribution to $\delta^{CF} \circ \delta^{CF}$ of trajectories satisfying one of these inequalities vanishes.
 
\begin{figure}[ht]
\begin{center}
\includegraphics{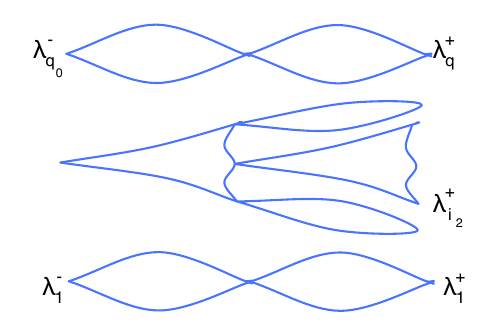}
\end{center}
\caption{$l_1 = l_2 = 3$, $q=7$, $q_1 = 5$, $q_0 = 3$, $i_2 = 3$ and $i_1 = 2$.}
\label{figurenoncommute}
\end{figure}

 Now if none of these inequalities is satisfied, Lemma \ref{lemmaassociahedron} implies the relation
 $(id_{i_1-1} \otimes m_{l_1} \otimes id_{q_0 - i_1}) \circ (id_{i_2-1} \otimes m_{l_2} \otimes id_{q_1 - i_2})
 = (-1)^{l_1 l_2 + (i_2 - i_1 + 1)(l_2 - 1)} id_{i_1-1} \otimes (m_{l_1} \circ m_{l_2})  \otimes id_{q_0 - i_1} $, where trajectories
 counted by $(m_{l_1} \circ m_{l_2})$ are counted with respect to the orientation induced by the associahedron $K_{l_1 + l_2 - 1}$
 on its boundary, see Figure \ref{figurenoncommute}. 
 The morphism $\delta^{CF} \circ \delta^{CF}$  counts this operator
 with respect to the sign 
 $q_1 l_1 + i_1 (l_1 - 1) + q l_2 + i_2 (l_2 - 1) = l_1 l _2 + (i_2 - i_1 + 1)(l_2 - 1) + (q + i_1 + 1)(l_1 + l_2) - 1 \mod(2)$.
 Since the quantities $q$, $l_1 + l_2$ and $i_1$ are constant on every connected component of ${\cal M}(\lambda^+ ,\lambda^- ; J^{\tilde{L}})$,
 and since the orientation of $\partial K_{l_1 + l_2 - 1}$ induces the outward normal orientation on every boundary point
 of ${\cal M}(\lambda^+ ,\lambda^- ; J^{\tilde{L}})$, we deduce the result. $\square$\\

Under the hypothesis of Theorem \ref{theocomplex}, for every Lagrangian conductor $\tilde{L}  \in \text{Ob} ( {\cal CL}^\pm (X, \omega))$,
we set ${\cal F} (\tilde{L}) = (CF (\tilde{L}) , \delta^{CF})  \in \text{Ob} ( K^b ({\cal OS}_1 (X, \omega)))$.

\begin{defi}
\label{deffloercomplex}
The complex ${\cal F} (\tilde{L})$ is called the Floer complex associated to $\tilde{L}$.
\end{defi}
 
\subsection{Stasheff's multiplihedron}
\label{subsectionmultiplihedron}

\subsubsection{Orientation}
 
\begin{figure}[ht]
\begin{center}
\includegraphics[scale=0.5]{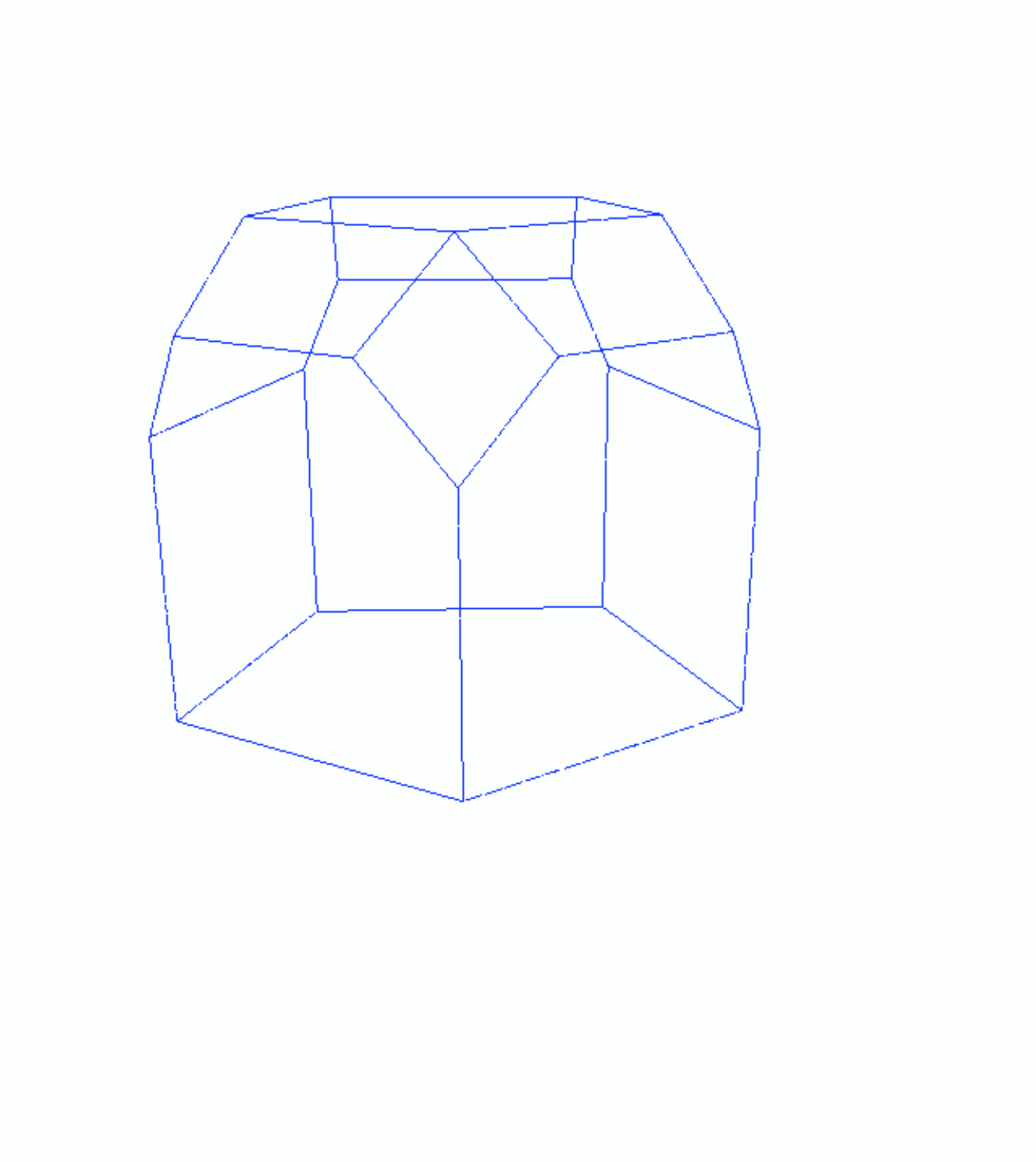}
\end{center}
\caption{The three-dimensional multiplihedron $J_4$.}
\label{figuremultiplihedron}
\end{figure}

\begin{figure}[ht]
\begin{center}
\includegraphics[scale=0.7]{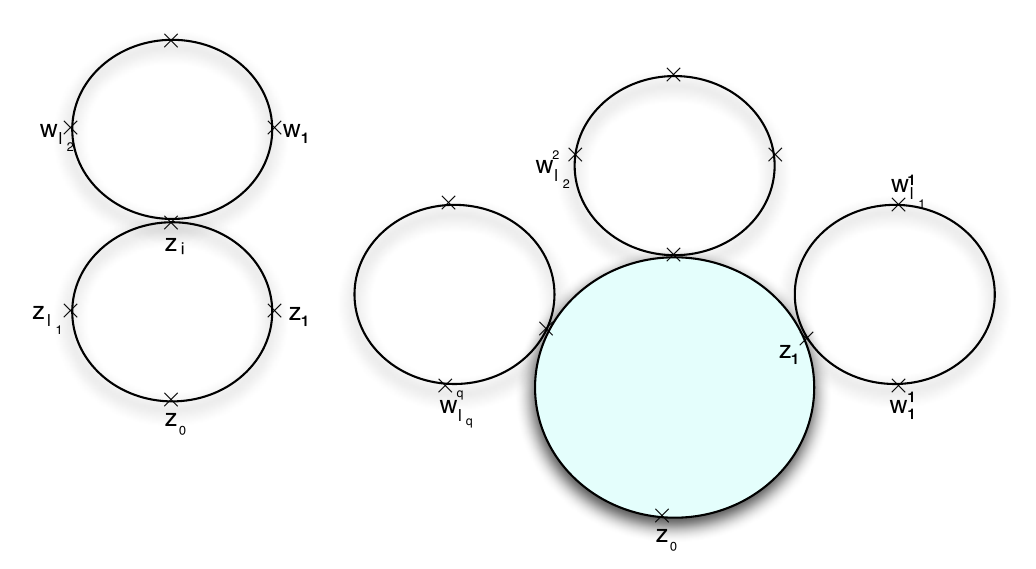}
\end{center}
\caption{Stable disks in a lower and an upper face of the multiplihedron.}
\label{figurelowerupper}
\end{figure}

Stasheff's multiplihedron $J_l$, $l \geq 2$, see \cite{Sta} and Figure \ref{figuremultiplihedron}, may be defined as the space of 
connected painted metric trees with $l+1$ free
 edges and such that the interior edges have lengths between $0$ and $1$, see \cite{Forc} and references therein. 
 The multiplihedron
 $J_l$ is a compactification of $[0,1] \times K_l$ which has the structure of a $(l-1)$-dimensional convex polytope of the Euclidian space. 
 We agree that $J_0 = J_1 = [0,1] $. Note that the propagation of the painting in the tree is non-trivial as soon as this tree
 contains a trivalent vertex, see Figure $3$.
Let us equip $\stackrel{\circ}{J}_l = ] 0,1 [ \times \stackrel{\circ}{K}_l $ with the product orientation. The codimension one faces
of $J_l$ different from $\{Ê0 \} \times K_l$ and $\{Ê1 \} \times K_l$ are of two different natures, see \cite{IwaMim}, \cite{Forc}.
The lower faces are canonically isomorphic to products $J_{l_1} \times K_{l_2}$, $l_1 + l_2 = l + 1$, they encode stable
disks having two irreducible components, one of which is unpainted. The upper faces are canonically isomorphic to products 
$K_q \times J_{l_1} \times  \dots \times  J_{l_q}$, $q \geq 1$, they encode stable
disks having $q$ irreducible components attached to a painted one, see Figure \ref{figurelowerupper}.

\begin{lemma}
\label{lemmalowerface}
Let $F = \{ (\underline{z} , \underline{w}) \in J_{l_1} \times K_{l_2} \, \vert \, z_i = w_0 \}$ be a lower facet
 of  the multiplihedron $J_l$, $l_1 + l_2 = l+1$, $l_1, l_2 \geq 2$, $1 \leq i \leq l_1$. Then, the orientations of $F$ induced
 by $\partial J_l$ and $J_{l_1} \times K_{l_2} $ coincide if and only if $l_1 l _2 + i(l_2 - 1)$ is even.
\end{lemma}

 {\bf Proof:}
 
The proof is analog to the one of Lemma \ref{lemmaassociahedron}. Assume first that $i>1$. Then, using the notations
adopted in the proof of Lemma \ref{lemmaassociahedron}, the orientation of $T_{\Theta_{\underline{z}} \star_\rho \Theta_{\underline{w}}} J_l$ writes:
 
$$\frac{\partial}{\partial s} \wedge (-\frac{\partial}{\partial z_2}) \wedge \dots \wedge (-\frac{\partial}{\partial z_{i-1}}) \wedge (-\frac{\partial}{\partial w_1}) \wedge \dots 
\wedge (-\frac{\partial}{\partial w_{l_2-1}}) \wedge (-\frac{\partial}{\partial z_{i+1}}) \wedge \dots \wedge (-\frac{\partial}{\partial z_{l_1}}) $$
 $$ =(-1)^{l_1 l _2 + i(l_2 - 1)} (\frac{\partial}{\partial \underline{w}}- \frac{\partial}{\partial \underline{z}}) \wedge \frac{\partial}{\partial s} \wedge (-\frac{\partial}{\partial z_2}) \wedge \dots \wedge \widehat{(-\frac{\partial}{\partial z_{i}})} \wedge \dots \wedge (-\frac{\partial}{\partial z_{l_1}}) \wedge (-\frac{\partial}{\partial w_2}) \wedge \dots \wedge (-\frac{\partial}{\partial w_{l_2-1}}).$$
 
 Since $(\frac{\partial}{\partial \underline{w}}- \frac{\partial}{\partial \underline{z}})$ is identified with the outward normal of $J_l$ at
 $F$, $(-\frac{\partial}{\partial z_2}) \wedge \dots \wedge \widehat{(-\frac{\partial}{\partial z_{i}})} \wedge \dots \wedge (-\frac{\partial}{\partial z_{l_1}})$
 form a direct basis of $T_{\underline{z}} J_{l_1}$ and $(-\frac{\partial}{\partial w_2}) \wedge \dots \wedge (-\frac{\partial}{\partial w_{l_2-1}})$ form a direct basis of  $T_{\underline{w}} K_{l_2}$, we deduce the result when $i>1$. The proof goes along the same lines when $i=1$. $\square$

\begin{lemma}
\label{lemmaupperface}
Let $F = \{ (\underline{z} , (s_1 , \underline{w}^1) , \dots , (s_q , \underline{w}^q)) \in K_q \times J_{l_1} \times \dots \times J_{l_q} \, \vert \, z_i = w_0^i
, 1 \leq i \leq q  \}$ be an upper facet
 of  the multiplihedron $J_l$, where $l_1 + \dots + l_q = l$, $q \geq 2$ and $l_i \geq 2$ for every $1 \leq i \leq q$. Then, the orientations of $F$ induced
 by $\partial J_l$ and $K_q \times J_{l_1} \times \dots \times J_{l_q}$ coincide if and only if $\sum_{i=1}^q (q-i)(l_i - 1)$ is even.
\end{lemma}

 {\bf Proof:}
 
Let $(\underline{z} , (s_1 , \underline{\tilde{w}}^1) , \dots , (s_q , \underline{\tilde{w}}^q)) \in F$ and 
$(s , \underline{w}^1 , \dots , \underline{w}^q)$ be an interior point of $J_l$ close to this point.
It follows from Definition $5.1$ of \cite{Forc} that $s = s_i + L_i$ for every $1 \leq i \leq q$, where $L_i$ is an increasing
function of the modulus of the neck which links the component $\Delta_{\underline{z}}$ to the component
$\Delta_{\underline{\tilde{w}}^i}$ of the stable disk encoded by $(\underline{z} , (s_1 , \underline{\tilde{w}}^1) , \dots , (s_q , \underline{\tilde{w}}^q))$.
Let us choose the model of the complex upper half plane, so that
$w_1^1 < \dots < w^1_{l_1} < w_1^2 < \dots < w^2_{l_2} < \dots < w_1^q < \dots < w^q_{l_q}$ and $\underline{w}^i$
is close to the point $z_i \in \R$, $1 \leq i \leq q$. The orientation of $T_{(s , \underline{w}^1 , \dots , \underline{w}^q)} J_l$ writes
$\frac{\partial}{\partial s} \wedge \widehat{(\sum_{i=1}^q \frac{\partial}{\partial \underline{w}^i})} \wedge 
\widehat{(\sum_{i=1}^q \frac{\partial}{\partial \underline{w}^i} - \frac{\partial}{\partial w_1^1})} \wedge \frac{\partial}{\partial w_3^1} \wedge \dots 
\wedge \frac{\partial}{\partial w_{l_1}^1} \wedge \frac{\partial}{\partial \underline{w}^2} \wedge \dots 
\wedge \frac{\partial}{\partial \underline{w}^q}$, where the vectors $(\sum_{i=1}^q \frac{\partial}{\partial \underline{w}^i})$ and 
$(\sum_{i=1}^q \frac{\partial}{\partial \underline{w}^i}) - \frac{\partial}{\partial w_1^1}$ are omitted, they are tangent to the orbit of
$\text{aut} (\H ; z_0)$, see \S \ref{subsubsectionorientationassociahedron}. For every $1 \leq i \leq q$, we identify
$\frac{\partial}{\partial w_1^i} + \dots + \frac{\partial}{\partial w_{l_i}^i}$ with $\frac{\partial}{\partial z_i}$ and
$\frac{\partial}{\partial w_{l_i}^i} - \frac{\partial}{\partial w_1^i} $ with $- \frac{\partial}{\partial L_i}$. This orientation rewrites

$$\frac{\partial}{\partial s} \wedge \widehat{\frac{\partial}{\partial z_1}}  \wedge \widehat{\frac{\partial}{\partial w^1_2}} 
\wedge \frac{\partial}{\partial w_3^1} \wedge \dots \wedge \frac{\partial}{\partial w_{l_1}^1} \wedge 
\big( \frac{\partial}{\partial z_2} \wedge \frac{\partial}{\partial w_2^2} \wedge \dots \wedge \frac{\partial}{\partial w_{l_2 - 1}^2} 
\wedge (- \frac{\partial}{\partial L_2}) \big)$$ $$  \wedge \dots \wedge
\big( \frac{\partial}{\partial z_q} \wedge \frac{\partial}{\partial w_2^q} \wedge \dots \wedge \frac{\partial}{\partial w_{l_q - 1}^q} 
\wedge (- \frac{\partial}{\partial L_q}) \big)$$
$$= (-1)^{\sum_{i=1}^q (q-i)(l_i - 1)} \frac{\partial}{\partial s} \wedge
\big( \widehat{\frac{\partial}{\partial z_1}}  \wedge \widehat{\frac{\partial}{\partial z_2}} \wedge \frac{\partial}{\partial z_3} \wedge \dots \wedge \frac{\partial}{\partial z_q} \big) \wedge
\big( \frac{\partial}{\partial w_2^1} \wedge \dots \wedge \frac{\partial}{\partial w_{l_1 - 1}^1} 
\wedge \frac{\partial}{\partial s_1} \big)$$ $$ \wedge \dots \wedge
\big( \frac{\partial}{\partial w_2^q} \wedge \dots \wedge \frac{\partial}{\partial w_{l_q - 1}^q} 
\wedge \frac{\partial}{\partial s_q} \big),$$
hence the result. $\square$

\subsubsection{Bordism of almost-complex structures}

Let $H$ be an effective continuation from the Lagrangian conductor $\tilde{L} = (L_0 , \dots , L_q ; J^{\tilde{L}})$ to 
$\tilde{L}' = (L'_0 , \dots , L'_l ; J^{\tilde{L}'})$, see Definition \ref{defeffectivecontinuation}. Restricting ourselves to 
the subconductor of $\tilde{L}' $ image of $H$, we may assume that $q=l$.
The continuation $H$ provides in particular, for every $0 \leq i \leq l$,  a homotopy
$s \in [0,1] \mapsto \tilde{L}_i^s =  (L_i^s , J_i^s)$ of elementary Lagrangian conductors
from the $i^{th}$  elementary Lagrangian sub conductor $\tilde{L}_i^0$ of $\tilde{L}$ to
the $i^{th}$  elementary Lagrangian sub conductor $\tilde{L}_i^1$ of $\tilde{L}'$. For every $0 \leq i < j \leq l$, we
extend these homotopies to generic homotopies $s \in [0,1] \mapsto \tilde{L}_{i, j}^s =  ((L_i^s , L_j^s) , J_{i,j}^s)$ 
so that $\tilde{L}_{i, j}^0$ and $\tilde{L}_{i, j}^1$ are subconductors of $\tilde{L} $ and $\tilde{L}'$ respectively.
We asociate to this homotopy a function $(s,t) \in \R \times [-1 , 1] \mapsto J_{i,j}^H (s,t) \in \overline{\cal J}_\omega$
defined by $J_{i,j}^H (s,t) = J_{i,j}^0 (t)$ if $s \leq 0$, $J_{i,j}^H (s,t) = J_{i,j}^s (t)$ if $0 \leq s \leq 1$ and
$J_{i,j}^H (s,t) = J_{i,j}^1 (t)$ if $s \geq 1$. We may assume withtout loss of generality that this function has the 
same regularity as the one chosen for our almost complex structures throughout the paper. Forgetting the painting
of trees provides a map $J_l \to K_l$ and we denote by $U_{J_l}$ the pull back of the universal curve $U_l \to K_l$
by this map. For every point $(0, \underline{z}) \in \{ 0 \} \times K_l \subset J_l$, we attach a copy of $\R \times [-1 , 1]$,
denoted by $\R \times [-1 , 1]_{z_0}$, to $D_{\underline{z}} \subset U_{J_l} \vert_{(0, \underline{z})}$ by identifying the
point $+ \infty$ of $\R \times [-1 , 1]_{z_0}$ with $z_0 \in D_{\underline{z}}$. Likewise, to every point $(1, \underline{z}) \in \{ 1 \} \times K_l \subset J_l$
and every $1 \leq j \leq l$, we attach a copy of $\R \times [-1 , 1]$,
denoted by $\R \times [-1 , 1]_{z_j}$, to $D_{\underline{z}} \subset U_{J_l} \vert_{(1, \underline{z})}$ by identifying the
point $- \infty$ of $\R \times [-1 , 1]_{z_j}$ with $z_j \in D_{\underline{z}}$. Finally, for every point $(s , \underline{z}) \in \partial J_l$
which encodes a painted metric tree having one bivalent vertex with both adjacent edges of length one, one painted,
the other one unpainted, we denote by $\sigma_{i,j} (s , \underline{z})$ the puncture of $D_{\underline{z}} \subset U_{J_l} \vert_{(s, \underline{z})}$
encoded by this vertex, where $i$ and $j$ are the labels of the boundary components of $D_{\underline{z}}$ adjacent to
$\sigma_{i,j} (s , \underline{z})$. We then attach a copy of $\R \times [-1 , 1]$,
denoted by $\R \times [-1 , 1]_{\sigma_{i,j} (s , \underline{z})}$, between the two irreducible components of 
$D_{\underline{z}}$ adjacent to the puncture $\sigma_{i,j} (s , \underline{z})$. Denote by 
$\overline{U}_{J_l}$ the union of $U_{J_l}$ and all these strips that we have just added. We then extend the coherent systems of
strip like ends $V_l$ and $V'_l$ given by $\tilde{L}$ and $\tilde{L}'$ to a system $V_{J_l}$  defined over the whole
$\overline{U}_{J_l}$, see Definition \ref{defstriplikeend}.

\begin{prop}
\label{propbordism}
Let $H$ be an effective continuation from the Lagrangian conductor $\tilde{L} = (L_0 , \dots , L_q ; J^{\tilde{L}})$ to 
$\tilde{L}' = (L'_0 , \dots , L'_l ; J^{\tilde{L}'})$. Using the notations just adopted, there exists a deformation of stable maps
$J^H : \overline{U}_{J_l} \to \overline{\cal J}_\omega$ over $J_l$ satisfying the following four properties.

1) For every $0 \leq i \leq l$ and every $(s ,\underline{z}) \in J_l$, $J^H (\partial_i D_{\underline{z}}) \subset \cup_{s \in [0 , 1]} J_i^s \in 
\overline{\cal J}_\omega$, so that in particular it does not depend on $\underline{z}$.

2) The restriction of $J^H$ to every strip $\R \times [-1 , 1]_{\sigma_{i,j} (s , \underline{z})}$ coincides with $J_{i,j}^H $, while its restriction to
$\R \times [-1 , 1]_{z_0}$ equals $J_{i_0,i_l}^H $ and its restriction to $\R \times [-1 , 1]_{z_i}$ equals $J_{i-1,i}^H $. Moreover, the restriction
of $J^H$ to every strip like end $\R_+ \times [-1 , 1]$ given by $V_{J_l}$ does not depend on the first factor $\R_+$ as soon as
$\tau \in \R_+$ is big enough. This restriction equals $J_{i,j}^0$ or $J_{i,j}^1$ depending on whether the edge of the painted tree
corresponding to this end is painted or not, where $i$ and $j$ are the labels of the boundary components $\R_+ \times \{Ê-1 \}$
and $\R_+ \times \{Ê1 \}$ of $\R_+ \times [-1 , 1]$.

3) For every $(s, \underline{z}) \in \partial J_l$ and every irreducible component $\Delta_{\underline{z}}$ of $D_{\underline{z}}$ corresponding
to a painted (resp. unpainted) metric subtree of $(s, \underline{z}) $, the restriction of $J^H$ to $\Delta_{\underline{z}}$ coincides with the
restriction of $J^{\tilde{L}}$ (resp. $J^{\tilde{L}'}$)  to $\Delta_{\underline{z}}$. In particular, the restriction of $J^H$ to $U_{J_l}$
coincides with $J^{\tilde{L}'}$ (resp. $J^{\tilde{L}}$) over $\{ 0 \} \times K_l$ (resp. $\{ 1 \} \times K_l$).

4) For every $(s , \underline{z}) \in J_l$ and every irreducible component $\Delta_{\underline{z}}$ of $D_{\underline{z}}$,
 the restriction of $J^H$ to $\Delta_{\underline{z}}$ does not depend on the conformal structure of $D_{\underline{z}} \setminus
 \Delta_{\underline{z}}$.\\
\end{prop}

Note that every face $F$ of $\partial J_l$  is canonically isomorphic to a product of lower dimensional associahedra or
multiplihedra. If such a product contains $j \geq 2$ factors, the stable disk $\Delta_{\underline{z}}$ associated to any interior
point $(s , \underline{z}) $ of $F$ cointains $j$ irreducible components in bijection with the factors of the product. Property $4$
given by Proposition \ref{propbordism} is a compatibility property of $J^H$ with this product structure. We denote by
$\overline{\cal J}_\omega (\overline{U}_{J_l} , V_{J_l})$ the space of extensions $J^H$ of $J^{\tilde{L}}$, $J^{\tilde{L}'}$
given by Proposition \ref{propbordism}, it is arc-connected. \\

{\bf Proof of Proposition \ref{propbordism}:}

\begin{figure}[ht]
\begin{center}
\includegraphics[scale=0.65]{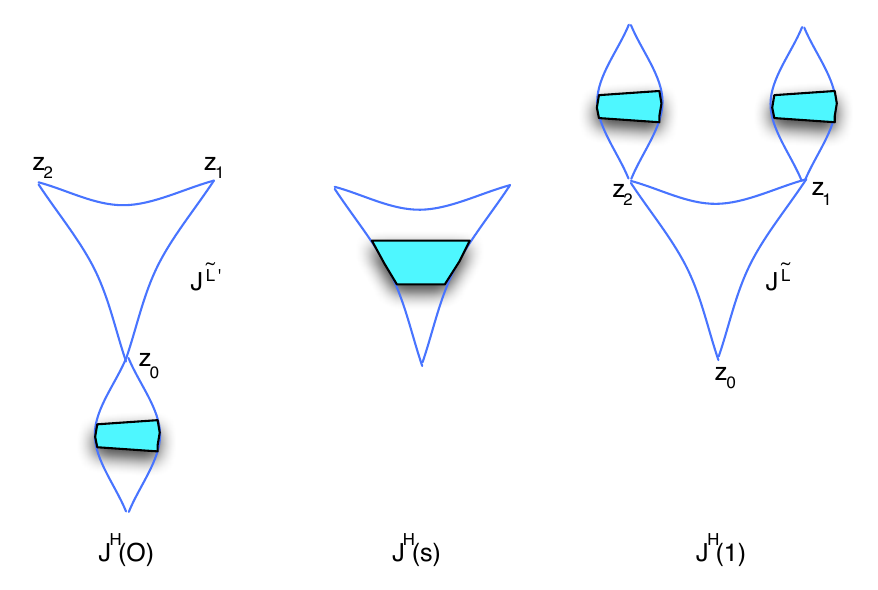}
\end{center}
\caption{Propagation of a homotopy over the multiplihedron.}
\label{figurepropagation}
\end{figure}

For every $0 \leq i_0 < i_1 < i_2 \leq l$, using notations of \S \ref{subseclagrangianconductors}, $U_{i_0 , i_1 , i_2}$ is a 
punctured disk $D_{\underline{z}}$ with three punctures $\underline{z} = \{ z_0 , z_1 , z_2 \}$ equipped with two maps
$J_{i_0 , i_1 , i_2}^{\tilde{L}}$ and $J_{i_0 , i_1 , i_2}^{\tilde{L}'}$ with values in $\overline{\cal J}_\omega$. The existence of
an extension $J^H$ of these maps over $J_2 = [0,1]$ is no problem, it is illustrated by Figure \ref{figurepropagation}.
By induction, assume the existence of the extension $J^H$ over all Lagrangian subconductors of size $q$, $2 \leq q < l$,
and let us prove this existence over a subconductor of size $q+1$. Since every face of $\partial J_q$  is a 
product of lower dimensional associahedra or multiplihedra, our already chosen extensions together with the compatibility
property prescribe $J^H$ over the boundary of $J_q$. There is then no obstruction to extend $J^H$ over the whole
$J_q$. $\square$

\subsection{Floer functor}

The aim of this paragraph is to associate to every continuation $H$ between Lagrangian conductors an element
${\cal F} (H) \in \text{Hom} ( K^b ({\cal OS}_1 (X, \omega)))$ in order to define a contravariant functor
${\cal F} : {\cal CL}^\pm (X, \omega) \to K^b ({\cal OS}_1 (X, \omega))$ called the Floer functor.

\subsubsection{Floer continuations}
\label{subsubsectionfloercontinuations}

Let $H$ be an effective continuation from the Lagrangian conductor $\tilde{L} = (L_0 , \dots , L_{l} ; J^{\tilde{L}})$ to 
$\tilde{L}' = (L'_0 , \dots , L'_{l'} ; J^{\tilde{L}'})$, see Definition \ref{defeffectivecontinuation}. Restricting ourselves to 
the subconductor of $\tilde{L}' $ image of $H$, we may assume that $l'=l$. Let $J^H \in \overline{\cal J}_\omega (\overline{U}_{J_l} , V_{J_l})$ 
be a generic extension of $J^{\tilde{L}}$, $J^{\tilde{L}'}$ given by Proposition \ref{propbordism}.
For every $1 \leq q \leq l$ and $0 \leq i_0 < \dots < i_q \leq l$, denote by $H_q : CF (L'_{i_0} , L'_{i_1}) \otimes \dots \otimes CF (L'_{i_{q-1}} , L'_{i_q})
\to CF (L_{i_0} , L_{i_q})$ the sum of elementary trajectories $\gamma_{\lambda^+ \lambda^-} = (u, D_{\underline{z}} , \lambda_{\partial u} , \overline{\partial}_u )$  going from an open string $\lambda^+ \in CF (L'_{i_0} , L'_{i_1}) \otimes \dots \otimes CF (L'_{i_{q-1}} , L'_{i_q})$ 
to an open string $\lambda^- \in CF (L_{i_0} , L_{i_q})$ satisfying $\mu (\lambda^-) - 1 = \mu(\lambda^+) - q$, such that:

1) $D_{\underline{z}}$ has one negative puncture $z_0$ and $q$ positive punctures $z_1, \dots z_q$.

2) There exists a lift $(s , \underline{z}) $ of $\underline{z}$ in $J_l$ such that 
$u : D_{\underline{z}} \to X$ satisfies the Cauchy-Riemann equation $J^H \vert_{(s , \underline{z})} \circ du = du \circ J_{D_{\underline{z}}}$.
Property $1$ of Proposition \ref{propbordism} provides then a map $h_j : \partial_j D_{\underline{z}} \to [0 , 1]$, $0 \leq j \leq q$,
such that  $\forall w \in \partial_j D_{\underline{z}} , J^H \vert_{(s , \underline{z})} (w) = J^{h_j (w)}_{i_j}$.

3) For every $0 \leq j \leq q$ and every $w \in \partial_j D_{\underline{z}}$, $u(w) \in L_{i_j}^{h_j (w)}$ and $\lambda_{\partial u} (w) = 
(T_{u(w)} L_{i_j}^{h_j (w)} , \p^\pm_{L_{i_j}^{h_j (w)}})$.

4) $\overline{\partial}_u$ is the Cauchy-Riemann operator associated to $u$, namely the first order deformation of its Cauchy-Riemann equation. 
The orientation of $\overline{\partial}_u$ is induced from the following exact sequence.
$$
0  \to  L^{k,p} (D_{\underline{z}} ; TD_{\underline{z}} , T \partial D_{\underline{z}})  \stackrel{du}{\to} 
L^{k,p} (D_{\underline{z}} ; u^* TX , \lambda_{\partial u})  \to  
L^{k,p} (D_{\underline{z}} ; u^* TX , \lambda_{\partial u})/\text{Im}(du) \to  0 $$
$$\downarrow 0 \oplus \overline{\partial}_{\underline{z}} \hspace{4cm} \downarrow \overline{\partial}_u \hspace{4cm} 
\downarrow \overline{\partial}'_N$$
$$0 \to \R \oplus L^{k-1,p} (D_{\underline{z}}  , \Lambda^{0,1} D_{\underline{z}} ) \stackrel{du}{\to}  
L^{k-1,p} (D_{\underline{z}}  , \Lambda^{0,1} D_{\underline{z}}  \otimes u^* TX)  \to 
L^{k-1,p} (D_{\underline{z}}  , \Lambda^{0,1} D_{\underline{z}}  \otimes u^* TX)/\R \oplus \text{Im}(du)  \to  0, $$
since $\overline{\partial}'_N$ is an isomorphism and the chosen orientation of $J_l$ induces an orientation of 
$0 \oplus  \overline{\partial}_{\underline{z}} $, see Remark \ref{remdbarz}. As in the proof of Proposition \ref{propindex}, Gromov-Floer compactness
 Theorem ensures that $H_q$ satisfies the Novikov condition of Definition \ref{deftrajectory} thanks to the condition on disks of Maslov indices less 
 than three given in Definition \ref{defelemconductor}.
 
 Let then ${\cal F} (H) : CF (\tilde{L}') \to CF (\tilde{L}) $ be the morphism defined by
 ${\cal F} (H) = \oplus_{q=1}^l \oplus_{l_1 + \dots l_{q_1} = q} (-1)^{\sum_{j=1}^{q_1} (q_1 - j)(l_j - 1)} H_{l_1} \otimes \dots \otimes H_{l_{q_1}}$.
In other words, ${\cal F} (H)$ is the morphism satisfying axioms $B$ of \S \ref{subsectionchain}
with dual given on elementary strings by $H_q^*$, $1 \leq q \leq l$.
In the case $l' > l$, we extend ${\cal F} (H) $ by $0$ outside of the sub complex of $CF (\tilde{L}') $ which is given by the subconductor
of $\tilde{L}' $ image of $H$. Note that the image of ${\cal F} (H)$ is the subcomplex associated to the cokernel of $H$, see \S
 \ref{subseclagrangianconductors}.
 
 \begin{theo}
\label{theocontinuation}
Let $(X, \omega)$ be a semipositive symplectic manifold of dimension $2n \geq 4$ which is either closed
or convex at infinity and $H$ be an effective continuation from the Lagrangian conductor $\tilde{L}$ to $\tilde{L}'$.
Then, ${\cal F} (H) : CF (\tilde{L}') \to CF (\tilde{L}) $ is a chain map.
\end{theo}

\begin{defi}
The chain map ${\cal F} (H) $ is called the Floer continuation associated to $H$.
\end{defi}
These continuations  ${\cal F} (H) $ have been introduced by Floer when $l=1$ and by Donaldson, Fukaya
when $l > 1$.\\

{\bf Proof:}

Let $\lambda^+ \in CF (\tilde{L}')$ and $\lambda^- \in CF (\tilde{L})$ be such that $\mu (\lambda^-) - 1 = \mu(\lambda^+) - q + 1$, so that
the space ${\cal M}(\lambda^+ ,\lambda^- ; J^H)$ of trajectories satisfying the same four conditions as the ones defining ${\cal F} (H) $ 
is one-dimensional. Genericness of $J^H$, semipositivity of $(X, \omega)$ and the condition on disks of Maslov indices less than three 
given in Definition \ref{defelemconductor} ensure compactness of this space, see Proposition $4.3$ and Theorem $5.1$ of \cite{Oh1}
as well as \S \ref{subsectcompactness}. From Floer glueing Theorem \cite{Floer}, \cite{BirCor}, the union of elementary trajectories 
from $\lambda^+$ to $\lambda^-$ counted by $\delta^{CF} \circ {\cal F} (H) - {\cal F} (H) \circ \delta^{CF'}$ 
 is in bijection with the boundary of this space. It suffices thus to prove that every such trajectory
 induces the outward normal orientation on ${\cal M}(\lambda^+ ,\lambda^- ; J^H)$. 
 
 Lemma \ref{lemmaupperface} provides the relation: \\
 
 $
 \begin{array}{cl}
& (id_{i-1} \otimes m_{l_0} \otimes id_{q_0 - i}) \circ (H_{l_1} \otimes \dots \otimes H_{l_{q_1}}) \\
 =& (-1)^{l_0 ( l_1 + \dots + l_{i-1} + i-1)}  H_{l_1} \otimes \dots \otimes H_{l_{i-1}} \otimes (m_{l_0} \otimes H_{l_i} \otimes \dots \otimes H_{l_{i + l_0 -1}} )
  \otimes H_{l_{i + l_0}} \otimes \dots \otimes H_{l_{q_1}} \\
 = &  (-1)^{l_0 ( l_1 + \dots + l_{i-1} + i-1) +   \sum_{j=1}^{l_0} (l_0 - j)(l_{i+j-1} - 1)} m_{l_0} \circ H, \\
 \end{array}$
 
 where trajectories
 counted by $m_{l_0} \circ H$ are counted with respect to the orientation induced by the multiplihedron on its upper faces.
 The morphism $\delta^{CF} \circ {\cal F} (H)$ counts this operator with respect to the sign 
 $(-1)^{q_1 l_0 + i (l_0 - 1)  +   \sum_{j=1}^{q_1} (q_1 - j)(l_j - 1)}$. Summing up, $m_{l_0} \circ H$ is counted  with respect to the sign 
 given by the following quantity modulo two
\begin{eqnarray}
\label{eqnsign1}
 \sum_{j=i + l_0}^{q_1} (q_1 - j)(l_j - 1) +  \sum_{j=i}^{i + l_0 -1} \epsilon (l_j - 1) + \sum_{j=1}^{i-1} (q_1 - j)(l_j - 1) + l_0 ( l_1 + \dots + l_{i-1} + q_1-1) - i,
\end{eqnarray}
where $\epsilon = q_1 + l_0 +i +1$. Likewise, Lemma \ref{lemmalowerface} provides the relation: \\
 
 $
 \begin{array}{cl}
& (H_{l_1} \otimes \dots \otimes H_{l_{i-1}} \otimes H_{l'} 
  \otimes H_{l_{i + l_0}} \otimes \dots \otimes H_{l_{q_1}} ) \circ (id_{i'-1} \otimes m_{l'_0} \otimes id_{l-l'_0 +1 - i'}) \\
  =& (-1)^{l'_0 ( l_{i+l_0} + \dots + l_{q_1} + q_1 + 1 - i -l_0)}  H_{l_1} \otimes \dots \otimes H_{l_{i-1}} \otimes (H_{l'} \circ m_{l'_0})
  \otimes H_{l_{i + l_0}} \otimes \dots \otimes H_{l_{q_1}} \\
 = &  (-1)^{l' l'_0 + i'' (l'_0 -1) + l'_0 ( l_{i + l_0} + \dots + l_{q_1} + q_1 + 1 - i-l_0) } H \circ m_{l'_0}, \\
 \end{array}$
 
 where $i' =  l_1 + \dots + l_{i-1} + i''$ and trajectories
 counted by $H \circ m_{l'_0}$ are counted with respect to the orientation induced by the multiplihedron on its lower faces.
 The morphism ${\cal F} (H) \circ \delta^{CF}$ counts this operator with respect to the sign 
$\sum_{j=i + l_0}^{q_1} (q_1 - j)(l_j - 1) + \epsilon (l' - 1) + \sum_{j=1}^{i-1} (i - j - \epsilon)(l_j - 1) + l l'_0 + i' (l'_0 - 1)  \mod (2)$,
since $\epsilon$ is odd if $q_1 - l_0 - i$ is even and even otherwise. Summing up, $H \circ m_{l'_0}$ is counted by ${\cal F} (H) \circ \delta^{CF}$ 
with respect to the sign 
\begin{eqnarray}
\label{eqnsign2}
 \sum_{j=i + l_0}^{q_1} (q_1 - j)(l_j - 1) +  \epsilon (l_i + \dots + l_{i + l_0 -1}) + \sum_{j=1}^{i-1} (i - j - \epsilon)(l_j - 1) + l_1 + \dots + l_{i-1}
 \mod (2),
\end{eqnarray}
 since $l' + l'_0 - 1 = l_i + \dots + l_{i + l_0 -1}$.
 
 In order to prove that the difference $\delta^{CF} \circ {\cal F} (H) - {\cal F} (H) \circ \delta^{CF}$ vanishes, it suffices
 to prove that the difference between (\ref{eqnsign1}) and (\ref{eqnsign2}) is odd.The latter writes
 
$
 \begin{array}{cl}
 &\sum_{j=1}^{i-1} (q_1 + i - \epsilon)(l_j - 1) + l_0 ( l_1 + \dots + l_{i-1} +  \epsilon  + q_1 - 1) + l_1 + \dots + l_{i-1} -i\\
 =& l_0 (\epsilon  + q_1 + i) - 1 = 1 \mod (2). \; \square
  \end{array}$ 

\subsubsection{Functoriality}

 \begin{theo}
\label{theofunctor}
Let $(X, \omega)$ be a semipositive symplectic manifold of dimension $2n \geq 4$ which is either closed
or convex at infinity. If $H^0, H^1$ are effective continuations from the Lagrangian conductor $\tilde{L}$ to $\tilde{L}'$
which are homotopic with fixed extremities, then ${\cal F} (H^0), {\cal F} (H^1) : CF (\tilde{L}') \to CF (\tilde{L}) $
are homotopic. Likewise, if $H^0 : \tilde{L} \to \tilde{L}'$ and $H^1 : \tilde{L}' \to \tilde{L}''$ are effective continuations,
then ${\cal F} (H^1 \circ H^0) = {\cal F} (H^0) \circ {\cal F} (H^1)$.
\end{theo}
Theorem \ref{theofunctor} means that ${\cal F} $ quotients out to a map 
$\text{Hom} (  {\cal CL}^\pm (X, \omega) ) \to  \text{Hom} ( K^b ({\cal OS}_1 (X, \omega)))$ which together with the map
given by Definition \ref{deffloercomplex} provides a contravariant functor
${\cal F} : {\cal CL}^\pm (X, \omega) \to K^b ({\cal OS}_1 (X, \omega))$.

\begin{defi}
\label{deffloerfunctor}
The contravariant functor
${\cal F} : {\cal CL}^\pm (X, \omega) \to K^b ({\cal OS}_1 (X, \omega))$ is called the Floer functor.
\end{defi}

{\bf Proof:}

Let $(H^r)_{r \in [0,1]}$ be a generic homotopy between $H^0 $ and $H^1$. Equip $[0,1] \times J_l$ with the product orientation
and the pull-back  of the universal curve $\overline{U}_{J_l}$, denoted by $\overline{U}_{[0,1] \times J_l}$. Let 
$J^H : \overline{U}_{[0,1] \times J_l} \to  \overline{\cal J}_\omega$ be the map whose restriction to every slice
$\{ r \} \times J_l$ is $J^{H^r}$, $r \in [0,1]$. In the same way as we defined the morphisms $H_q$ in the previous 
\S \ref{subsubsectionfloercontinuations}, we define, 
for every $1 \leq q \leq l$, $0 \leq i_0 < \dots < i_q \leq l$ and $r \in [0,1]$, the morphisms
$K^r_q : CF (L'_{i_0} , L'_{i_1}) \otimes \dots \otimes CF (L'_{i_{q-1}} , L'_{i_q})
\to CF (L_{i_0} , L_{i_q})$ as the sum of elementary trajectories $\gamma_{\lambda^+ \lambda^-} = (u, D_{\underline{z}} , \lambda_{\partial u} , \overline{\partial}_u )$  which are $J^{H^r}$-holomorphic going from an open string $\lambda^+ \in CF (L'_{i_0} , L'_{i_1}) \otimes \dots \otimes 
CF (L'_{i_{q-1}} , L'_{i_q})$ 
to an open string $\lambda^- \in CF (L_{i_0} , L_{i_q})$ satisfying $\mu (\lambda^-) - 1 = \mu(\lambda^+) - q - 1$. These trajectories
are oriented from the orientation just fixed on $[0,1] \times J_l$ and can only appear for a discrete set of values of
$r \in [0,1]$ since the homotopy $(H^r)_{r \in [0,1]}$ is generic. Moreover, they satisfy the Novikov condition of Definition \ref{deftrajectory} 
thanks to the condition on disks of Maslov indices less than three given in Definition \ref{defelemconductor}, see Proposition \ref{propindex}.

Then, we denote by $K : CF (\tilde{L}') \to CF (\tilde{L}) $ the morphism defined by
 $K = \oplus_{q=1}^l (-1)^q \oplus_{l_1 + \dots + l_q = q'} (-1)^{\sum_{j=1}^{q} (q - j)(l_j - 1)}  \oplus_{i=1}^{q} (-1)^{\sum_{j=1}^{i-1} (l_j - 1)}
 \int_0^1  H^r_{l_1} \otimes \dots \otimes H^r_{l_{i-1}} \otimes K^r_{l_i} 
  \otimes H^r_{l_{i + 1}} \otimes \dots \otimes H^r_{l_{q}}$, where the integral is taken with respect to the counting measure.
In other words, $K$ is the morphism satisfying axioms $C$ of \S \ref{subsectionchain}
with dual given on elementary strings by $-K_q^*$, $1 \leq q \leq l$. This morphism $K$ satisfies the relation
${\cal F} (H^0) - {\cal F} (H^1) = \delta^{CF} \circ K + K \circ \delta^{CF} $.

Indeed, 
let $\lambda^+ \in CF (\tilde{L}')$ and $\lambda^- \in CF (\tilde{L})$ be such that $\mu (\lambda^-) - 1 = \mu(\lambda^+) - q $, so that
the space ${\cal M}(\lambda^+ ,\lambda^- ; J^H)$ of trajectories satisfying the same four conditions as the ones defining ${\cal F} (H) $ 
and which are $J^{H^r}$-holomorphic for some $r \in [0,1]$
is one-dimensional. Genericness of $J^H$, semipositivity of $(X, \omega)$ and the condition on disks of Maslov indices less than three 
given in Definition \ref{defelemconductor} ensure compactness of this space, see Proposition $4.3$ and Theorem $5.1$ of \cite{Oh1}
as well as \S \ref{subsectcompactness}. From Floer glueing Theorem \cite{Floer}, \cite{BirCor}, the union of elementary trajectories 
from $\lambda^+$ to $\lambda^-$ counted by ${\cal F} (H^1) - {\cal F} (H^0) + \delta^{CF} \circ K + K \circ \delta^{CF} $
 is in bijection with the boundary of this space. It suffices thus to prove that every such trajectory
 induces the outward normal orientation on ${\cal M}(\lambda^+ ,\lambda^- ; J^H)$. 
 
The analog of  Lemma \ref{lemmaupperface} provides the relation: \\
 
 $
 \begin{array}{cl}
& (id_{I-1} \otimes m_{l_0} \otimes id_{q_0 - I}) \circ  
(H_{l_1} \otimes \dots \otimes H_{l_{i-1}} \otimes K_{l_i}  \otimes H_{l_{i + 1}} \otimes \dots \otimes H_{l_{q_1}})\\
 =& (-1)^{l_0 ( l_1 + \dots + l_{I-1} + I-1)}  H_{l_1} \otimes \dots \otimes H_{l_{I-1}} \otimes (m_{l_0} \otimes H_{l_I} \otimes \dots \otimes
K_{l_i} \otimes \dots \otimes H_{l_{I + l_0 -1}} )\\
&  \otimes H_{l_{I + l_0}} \otimes \dots \otimes H_{l_{q_1}} \\
 = &  (-1)^{l_0 ( l_1 + \dots + l_{I-1} + I-1) +   \sum_{j=1}^{l_0} (l_0 - j)(l_{I+j-1} - 1)   + l_0 - 1 + l_I - 1 + \dots + l_{i-1} -1} m_{l_0} \circ H, \\
 \end{array}$

 where trajectories
 counted by $m_{l_0} \circ H$ are counted with respect to the orientation induced by the product $[0,1] \times J_l$ on its boundary.
 Indeed, using the notations introduced in the proof of Lemma \ref{lemmaupperface}, this relation comes from the identity\\
 
 $
 \begin{array}{l}
 \frac{\partial}{\partial r} \wedge \frac{\partial}{\partial s} \wedge \widehat{(\sum_{i=1}^q \frac{\partial}{\partial \underline{w}^i})} \wedge 
\widehat{(\sum_{i=1}^q \frac{\partial}{\partial \underline{w}^i} - \frac{\partial}{\partial w_1^1})} \wedge \frac{\partial}{\partial w_3^1} \wedge \dots 
\wedge \frac{\partial}{\partial w_{l_1}^1} \wedge \frac{\partial}{\partial \underline{w}^2} \wedge \dots 
\wedge \frac{\partial}{\partial \underline{w}^q} \\
=(-1)^{\sum_{j=1}^{q} (q - j)(l_{j} - 1)   + q - 1 + l_1 - 1 + \dots + l_{i-1} -1} 
\frac{\partial}{\partial s} \wedge \frac{\partial}{\partial \underline{z}} \wedge \frac{\partial}{\partial \underline{w}^1} \wedge \dots 
\wedge \frac{\partial}{\partial \underline{w}^{i-1}} \wedge (  \frac{\partial}{\partial r} \wedge \frac{\partial}{\partial \underline{w}^i}) 
 \wedge \dots \wedge \frac{\partial}{\partial \underline{w}^q}, \\
\end{array}$
 and from the fact that the commutation of $\frac{\partial}{\partial r}$ with the operator $H_{l_1} \otimes \dots \otimes H_{l_{I-1}} $ turns out
 to be compensated by the commutation of the outward normal vector $\frac{\partial}{\partial s}$ with this operator.
 
  The morphism $\delta^{CF} \circ K$ counts this operator with respect to the sign 
  
 $(-1)^{q_1 l_0 + I (l_0 - 1)  + q_1 +  \sum_{j=1}^{q_1} (q_1 - j)(l_j - 1) + \sum_{j=1}^{i-1} (l_j - 1)}$. 
 Summing up, $m_{l_0} \circ H$ is counted  with respect to the sign \\
 
$\begin{array}{l}
 l_0 ( l_1 + \dots + l_{I-1} + I-1) + \sum_{j=1}^{l_0} (l_0 - j)(l_{I+j-1} - 1) +  l_0 - 1
 + \sum_{j=1}^{I -1} (l_j - 1) + q_1 l_0 \\
 + I(l_0 -1) + q_1 + \sum_{j=1}^{q_1} (q_1 - j)(l_j - 1) \mod (2) \\
 = \sum_{j=1}^{q_1} (q_1 - j)(l_j - 1) + \sum_{j=I}^{I+ l_0 -1} (I + l_0 -1 - j)(l_{j} - 1) + (l_0 - 1) \sum_{j=1}^{I -1} (l_j - 1) \\
 + (q+1 -I)(l_0 -1) \mod (2) \\
 = \sum_{j=I + l_0}^{q_1} (q_1 - j)(l_j - 1) + (q_1 - I) ( l_I + \dots + l_{I + l_0 -1} -1) + \sum_{j=1}^{I-1} (q_0 - j)(l_j - I) \mod (2), \\
\end{array}$

since $q_0 = q_1 + 1 - l_0$. But it is exactly with respect to this sign that the operator ${\cal F} (H^1)$ counts the elements near
the boundary of ${\cal M}(\lambda^+ ,\lambda^- ; J^H)$, so that $\delta^{CF} \circ K $ induces the outward normal orientation
on the boundary of ${\cal M}(\lambda^+ ,\lambda^- ; J^H)$. 

Likewise, the analog of Lemma \ref{lemmalowerface} provides the relation: \\
 
 $
 \begin{array}{cl}
& (H_{l_1} \otimes \dots \otimes H_{l_{i-1}} \otimes K_{l_i} 
  \otimes H_{l_{i + 1}} \otimes \dots \otimes H_{l_{q_1}} ) \circ (id_{I-1} \otimes m_{l'_0} \otimes id_{l-l'_0 +1 - I}) \\
  =& (-1)^{l'_0 \sum_{j=i+1}^{q_1} (l_j - 1)}  H_{l_1} \otimes \dots \otimes H_{l_{i-1}} \otimes (K_{l_i} \circ m_{l'_0})
  \otimes H_{l_{i + 1}} \otimes \dots \otimes H_{l_{q_1}} \\
 = &  (-1)^{l'_0 \sum_{j=i+1}^{q_1} (l_j - 1)} (-1)^{ l_i l'_0 + i' (l'_0 -1) +  1} H \circ m_{l'_0}, \\
 \end{array}$

 where $I =  l_1 + \dots + l_{i-1} + i'$ and trajectories
 counted by $H \circ m_{l'_0}$ are counted with respect to the orientation induced by the product $[0,1] \times J_l$ on its boundary.
 Indeed, using the notations introduced in the proof of Lemma \ref{lemmalowerface}, this relation comes from the identity\\
 $
 \begin{array}{l}
 \frac{\partial}{\partial r} \wedge  \frac{\partial}{\partial s} \wedge (-\frac{\partial}{\partial z_2}) \wedge \dots \wedge (-\frac{\partial}{\partial z_{i-1}}) \wedge (-\frac{\partial}{\partial w_1}) \wedge \dots 
\wedge (-\frac{\partial}{\partial w_{l_2-1}}) \wedge (-\frac{\partial}{\partial z_{i+1}}) \wedge \dots \wedge (-\frac{\partial}{\partial z_{l_1}}) \\
  =(-1)^{l_1 l _2 + i(l_2 - 1) + 1} (\frac{\partial}{\partial \underline{w}}- \frac{\partial}{\partial \underline{z}}) \wedge 
  ( \frac{\partial}{\partial r} \wedge \frac{\partial}{\partial s} \wedge \frac{\partial}{\partial \underline{z}})  \wedge \frac{\partial}{\partial \underline{w}},\\ 
%   (-\frac{\partial}{\partial z_2}) \wedge \dots \wedge \widehat{(-\frac{\partial}{\partial z_{i}})} \wedge \dots \wedge (-\frac{\partial}{\partial z_{l_1}}) \wedge (-\frac{\partial}{\partial w_2}) \wedge \dots \wedge (-\frac{\partial}{\partial w_{l_2-1}}),\\
 \end{array}$
 
 where once more the commutation of $\frac{\partial}{\partial r}$ with the operator $H_{l_1} \otimes \dots \otimes H_{l_{I-1}} $ gets 
 compensated by the commutation of the outward normal vector 
 $(\frac{\partial}{\partial \underline{w}}- \frac{\partial}{\partial \underline{z}})$ with this operator.
 
  The morphism $K \circ \delta^{CF}$ counts this operator with respect to the sign 
$q_1 + \sum_{j=1}^{q_1} (q_1 - j)(l_j - 1) + \sum_{j=1}^{i-1} (l_j - 1) + l l'_0 + I (l'_0 - 1)  \mod (2)$. 
Summing up, $H \circ m_{l'_0}$ is counted by $K \circ \delta^{CF}$ with respect to the sign \\

$\begin{array}{l}
l'_0 \sum_{j=i+1}^{q_1} (l_j - 1) + l_i l'_0 + i' (l'_0 -1) +  1 + q_1 + \sum_{j=1}^{q_1} (q_1 - j)(l_j - 1) + \sum_{j=1}^{i-1} (l_j - 1) \\
+ l l'_0 + I (l'_0 - 1)  \mod (2) \\
= \sum_{j=1}^{q_1} (q_1 - j)(l_j - 1) + l'_0 \sum_{j=i}^{q_1} (l_j - 1) + l'_0 (l-1) + l'_0 \sum_{j=1}^{i-1} l_j + i + q_1  \mod (2) \\
= \sum_{j=1}^{q_1} (q_1 - j)(l_j - 1) + (q_1 - i)(l'_0 - 1)  \mod (2). \\
\end{array}$
 
 But it is again exactly with respect to this sign that the operator ${\cal F} (H^0)$ counts the elements near
the boundary of ${\cal M}(\lambda^+ ,\lambda^- ; J^H)$, so that $K   \circ \delta^{CF} $ induces the outward normal orientation
on the boundary of ${\cal M}(\lambda^+ ,\lambda^- ; J^H)$. Since the outward normal orientations at the two
boundary points of a connected component of ${\cal M}(\lambda^+ ,\lambda^- ; J^H)$ induce opposite orientations on this component,
we deduce the vanishing ${\cal F} (H^1) - {\cal F} (H^0) + \delta^{CF} \circ K + K \circ \delta^{CF} =0$, hence the first part of
Theorem \ref{theofunctor}.

Let us now write
$(H^0_{l'_1} \otimes \dots \otimes H^0_{l'_{q_0}} ) \otimes (H^1_{l_1} \otimes \dots \otimes H^1_{l_{q_1}} ) =  (-1)^s
(H^0_{l'_1} \otimes H^1_{l_1} \otimes \dots \otimes H^1_{l_{l'_1}} ) \otimes \dots \otimes
(H^0_{l'_{q_0}} \otimes H^1_{l_{l'_1 + \dots + l'_{q_0} + 1}} \otimes \dots \otimes H^1_{l_{q_1}} ) $, with
$s = \big( \sum_{i = l'_1 + \dots + l'_{q_0 -2} + 1}^{l'_1 + \dots + l'_{q_0-1}} (l_i - 1) \big) (l'_{q_0} - 1) +
\big( \sum_{i = l'_1 + \dots + l'_{q_0 -3} + 1}^{l'_1 + \dots + l'_{q_0-2}} (l_i - 1) \big) (l'_{q_0} + l'_{q_0-1}) + \dots +
\big( \sum_{i =1}^{l'_1} (l_i - 1) \big) (l'_2 + \dots + l'_{q_0} - q_0 + 1) \mod (2)$.
In this product, the $j^{th}$ factor $(H^0_{l'_j} \otimes H^1_{l_{l'_1 + \dots + l'_{j-1} + 1}} \otimes \dots \otimes H^1_{l_{l'_1 + \dots + l'_j}} )$
comes with the sign
$\big( \sum_{i = l'_1 + \dots + l'_{j-1} + 1}^{l'_1 + \dots + l'_{j}} (l_i - 1) \big) (\epsilon_j - q_0 + j)  $, where
$\epsilon_j = l'_{j+1} + \dots + l'_{q_0}$, $1 \leq j \leq q$, and where this $j^{th}$ factor is equipped with the orientation induced
from the product of $l'_j+1$ multiplihedra. In the same way as in \S \ref{subsectionmultiplihedron}, this product of multiplihedra
is a face of a polytope $J^2_l$ which encodes painted connected metric trees with $l+1$  free edges among which one root
and with interior edges having lengths in $[0,1]$. But this time, these trees are painted with two different colors,
one which encodes the propagation of the homotopy $H^0$ and the other one which encodes the propagation of the homotopy $H^1$.
We deduce as in Lemma \ref{lemmaupperface} that the orientation of this face induced by $J^2_l$ and the one induced
by the product of multiplihedra differ by the sign $(-1)^{Ê\sum_{i = l'_1 + \dots + l'_{j-1} + 1}^{l'_1 + \dots + l'_{j}} (l'_1 + \dots + l'_{j} - i)(l_i - 1)}$.
Indeed, using the notations of the proof of Lemma \ref{lemmaupperface}, this difference comes from the relation\\

$
 \begin{array}{l}
 (\frac{\partial}{\partial r} - \frac{\partial}{\partial s}) \wedge (\frac{\partial}{\partial r} + \frac{\partial}{\partial s}) \wedge 
 \widehat{(\sum_{i=1}^q \frac{\partial}{\partial \underline{w}^i})} \wedge 
\widehat{(\sum_{i=1}^q \frac{\partial}{\partial \underline{w}^i} - \frac{\partial}{\partial w_1^1})} \wedge \frac{\partial}{\partial w_3^1} \wedge \dots 
\wedge \frac{\partial}{\partial w_{l_1}^1} \wedge \frac{\partial}{\partial \underline{w}^2} \wedge \dots 
\wedge \frac{\partial}{\partial \underline{w}^q} \\
=(-1)^{\sum_{i=1}^{q} (q - i)(l_i - 1) } (\frac{\partial}{\partial r} - \frac{\partial}{\partial s}) \wedge
(\frac{\partial}{\partial r} \wedge \frac{\partial}{\partial \underline{z}} ) \wedge ( \frac{\partial}{\partial \underline{w}_2^1} \wedge 
\frac{\partial}{\partial s_1}) \wedge \dots  \wedge ( \frac{\partial}{\partial \underline{w}^q} \wedge 
\frac{\partial}{\partial s_q}), \\
\end{array}$

where $\frac{\partial}{\partial r} , \frac{\partial}{\partial s}$ represent the propagation of the first and second paints in the tree, so that
$\frac{\partial}{\partial r} + \frac{\partial}{\partial s}$ is tangent to every level $\{ \alpha \} \times \stackrel{\circ}{J}_l \subset J^2_l$,
$\alpha \in ] 0, 1 [$ and $\frac{\partial}{\partial r} - \frac{\partial}{\partial s}$ is transversal to these level sets and induces
the orientation of $[ 0 , 1]$.

But the operator  ${\cal F} (H^1) $ counts the operator 
$H^1_{l_{l'_1 + \dots + l'_{j-1} + 1}} \otimes \dots \otimes H^1_{l_{l'_1 + \dots + l'_j}} $ with respect to the sign
$(-1)^{\sum_{i = l'_1 + \dots + l'_{j-1} + 1}^{l'_1 + \dots + l'_{j}} (l'_1 + \dots + l'_j - i + \epsilon_j ) (l_i - 1)} $ whereas the operator
 ${\cal F} (H^0)$ counts the term $H^0_{l'_j}$  with respect to the sign $(-1)^{(q_0 - j) (l'_j - 1)}$. Summing up, the glued operator
 $H^0_{l'_j} \otimes H^1_{l_{l'_1 + \dots + l'_{j-1} + 1}} \otimes \dots \otimes H^1_{l_{l'_1 + \dots + l'_j}} $ is counted with
 respect to the sign 
 $(q_0 - j) \sum_{i = l'_1 + \dots + l'_{j-1} + 1}^{l'_1 + \dots + l'_{j}} (l_i - 1) + (q_0 - j) (l'_j - 1) = 
 (q_0 - j) (\sum_{i = l'_1 + \dots + l'_{j-1} + 1}^{l'_1 + \dots + l'_{j}} l_i - 1) \mod (2)$.
 The result follows from the fact that it is precisely with respect to the latter sign that the glued operator
 $H^0_{l'_j} \otimes H^1_{l_{l'_1 + \dots + l'_{j-1} + 1}} \otimes \dots \otimes H^1_{l_{l'_1 + \dots + l'_j}} $ is counted by
 ${\cal F} (H^1 \circ H^0) $. $\square$\\

 The image ${\cal F} ({\cal CL}^\pm (X, \omega))$ is a subcategory of  $K^b ({\cal OS}_1 (X, \omega))$
 whose objects are Floer complexes $CF (\tilde{L}) \stackrel{\delta^{CF}}{\to} CF (\tilde{L}) \stackrel{\delta^{CF}}{\to} CF (\tilde{L})$
 and morphisms are given by Floer continuations.  Denote by $\text{Ob} ({\cal FK} (X, \omega))$ the set of augmented
 Floer complexes $\emptyset \to CF (\tilde{L}) \stackrel{\delta^{CF}}{\to} CF (\tilde{L}) \stackrel{\delta^{CF}}{\to} CF (\tilde{L})$
 in the sense of Definition \ref{defaugmentation}, and by $\text{Hom} ({\cal FK} (X, \omega))$ the set of 
 Floer continuations between augmented
 Floer complexes. The pair ${\cal FK} (X, \omega) = (\text{Ob} ({\cal FK} (X, \omega)) , \text{Hom} ({\cal FK} (X, \omega)))$
 is a subcategory of $K^b ({\cal OS}_1 (X, \omega))$, it inherits the notion of refinements given by Definition \ref{defextensionlagconductor},
 see also Remark \ref{remrefinement}. We may then propose the following Definition \ref{defFloerKontsevich}.
 
 \begin{defi}
 \label{defFloerKontsevich}
 The category ${\cal FK} (X, \omega) $ is called the Floer-Kontsevich category.
 \end{defi}

Note that augmentations of Floer complexes are exactly the twisted complexes introduced in \cite{Kont}.
The latter may be thought of as triangular matrices acting on formal sums $L_0 \oplus \dots \oplus L_l$ in the
language of Donaldson-Fukaya. These twisted complexes are the objects of the category $D^b {\cal F} (X, \omega)$,
the derived category of the Fukaya category, introduced by Kontsevich in \cite{Kont}. Morphisms of this category are given
by refinements of augmented complexes, see Remark \ref{remrefinement}.

\subsection{Lifts of Floer functor}

We present in this paragraph two situations where Floer functor lifts to a functor 
${\cal F}_N : {\cal CL}^\pm (X, \omega) \to K^b ({\cal OS}_N (X, \omega))$, $N \in \N$. From Remark
\ref{remgraduation}, Floer complexes then inherit a non-trivial graduation modulo $N$ and the same holds for
their cohomology after application of the functor coefficients ${\cal C}$.

\subsubsection{Graded Lagrangian conductors}
\label{subsubsectgraded}

\begin{defi}
Let $N \in \N$. A ${\cal L}^\pm_N$-structure of $(X, \omega)$ is  a bundle
${\cal L}^\pm_N \to {\cal L}^\pm$ whose restriction over every fiber ${\cal L}^\pm_x$, 
$x \in X$, is the cyclic covering of order $N$.
\end{defi}
Here, $N=0$ corresponds to the infinite cyclic cover. As soon as twice the first Chern class of $(X, \omega)$ vanishes in 
$H^2 (X ; \Z / N\Z)$, such a ${\cal L}^\pm_N$-structure 
exist on $(X, \omega)$. The set of such structures is then a principal space over the group $H^1 (X ; \Z / N\Z)$, see
\cite{Seid}.

\begin{defi}
A graduation of a Lagrangian submanifold $L$ of $(X, \omega , {\cal L}^\pm_N)$  is a section
$\text{gr}_L : L \to {\cal L}^\pm_N \vert_L$ which lifts the tautological section $x \in L \mapsto T_x L \in {\cal L}^\pm_x$.
\end{defi}
This notion has been introduced by Kontsevich in \cite{Kont} and studied in detail in \cite{Seid}. In particular,
when $L$ is a sphere, such a graduation always exist, compare \S \ref{parc1}. When  $(X, \omega)$ is equipped
with a ${\cal L}^\pm_N$-structure , we denote by ${\cal CL}^\pm (X, \omega , {\cal L}^\pm_N)$ the category
of graded Lagrangian conductors of $(X, \omega , {\cal L}^\pm_N)$. Its objects are
elements $\tilde{L} = ((L_0 , \text{gr}_{L_0} ) , \dots , (L_l , \text{gr}_{L_l} ) ; J^{\tilde{L}})$  where  $(L_0 , \dots , L_l ; J^{\tilde{L}})$
is a Lagrangian conductor of $(X, \omega)$ and $ \text{gr}_{L_i} $ a graduation of $L_i$, $0 \leq i \leq l$, while its morphisms
are continuations which preserve graduations. For every $\tilde{L} = ((L_0 , \text{gr}_{L_0} ) , \dots , (L_l , \text{gr}_{L_l} ) ; J^{\tilde{L}})
\in \text{Ob} ( {\cal CL}^\pm (X, \omega , {\cal L}^\pm_N))$, $l>0$, graduations provide a tautological injection
$x \in L_i \cap L_j \mapsto \lambda_x \in \text{Ob} ( {\cal OS}_N (X, \omega))$, $0 \leq i < j \leq l$.
The construction of Floer functor then provides a functor ${\cal F}_N : {\cal CL}^\pm (X, \omega , {\cal L}^\pm_N)
\to K^b ({\cal OS}_N (X, \omega))$ for which the following diagram is commutative

$$\begin{array}{ccc}
{\cal CL}^\pm (X, \omega , {\cal L}^\pm_N) & \stackrel{{\cal F}_N}{\to} & K^b ({\cal OS}_N (X, \omega)) \\
\downarrow && \downarrow \\
{\cal CL}^\pm (X, \omega) & \stackrel{{\cal F}}{\to} & K^b ({\cal OS}_1 (X, \omega)), \\
\end{array}$$

where ${\cal CL}^\pm (X, \omega , {\cal L}^\pm_N) \to {\cal CL}^\pm (X, \omega) $ denotes the functor which forgets
graduations. From Remark \ref{remgraduation}, Floer complexes ${\cal F}_N (\tilde{L})$ inherit a non-trivial graduation
modulo $N$.

{\bf Example: $N=2$}

Every symplectic manifold $(X, \omega)$ carries a ${\cal L}^\pm_2$-structure, namely ${\cal L}^\pm_2$ denotes the
bundle of oriented Lagrangian subspaces of $TX$ equipped with a $\widetilde{GL}_n^\pm (\R)$-structure. A graduation
of a Lagrangian submanifold $L$ is then an orientation of $L$. When $l=1$, the Euler characteristic of the complex
${\cal F}_2 (\tilde{L})$  with respect to the graduation given by the function $\mu + q$, see Remark
\ref{remgraduation}, is the opposite of the intersection index $L_0 \circ L_1\in \Z$ where $L_0$, $ L_1$ are the two
oriented Lagrangian submanifolds in $\tilde{L}$.

\subsubsection{Based Lagrangian conductors}
\label{subsubsectbased}

Let us fix a base point $x_0 \in X$ together with a Lagrangian subspace $l_0 \in {\cal L}^\pm_{x_0}$ of $T_{x_0} X$.

\begin{defi}
A based Lagrangian submanifold of $(X, \omega)$  is a pair $(L , \lambda_L)$ where $L \in {\cal L}ag^\pm$ and
$\lambda_L : [0,1] \to {\cal L}^\pm$ satisfies $\lambda_L (0) = l_0$ and $\pi_{\cal L} \circ \lambda_L (1) \in L$.
\end{defi}

\begin{defi}
\label{defbasedcontinuation}
An effective continuation of based Lagrangian submanifolds of $(X, \omega)$  is a path $(L^s , \lambda_{L^s})_{s \in [0,1]}$ of
based Lagrangian submanifolds such that $(L^s)_{s \in [0,1]}$ is a Hamiltonian isotopy of Lagrangian submanifolds.
\end{defi}
Note that there always exist a family $\Phi_s (t)$ of Hamiltonian diffeomorphisms of $(X, \omega)$, $(s,t) \in [0,1]^2$,
such that

1) $\forall s \in [0,1]$, $\Phi_s (0) = Id$.

2) $\forall t \in [0,1]$, $\Phi_0 (t) = Id$.

3) $\forall s \in [0,1]$, $\Phi_s (1) (L) = L^s$ and $\lambda_{L^s} = \Phi_s \circ \lambda_L $.\\

When $(X, \omega)$ comes equipped with a ${\cal L}^\pm_N$-structure, we fix a lift $\tilde{l}_0$ of $l_0$ in
${\cal L}^\pm_N$. We then define a based graded Lagrangian submanifold to be a triple $(L , \text{gr}_L  ,\lambda_L)$
where $(L , \lambda_L)$ is based, $(L , \text{gr}_L )$ graded and $\lambda_L$ lifts to a path
$\tilde{\lambda}_L : [0,1] \to {\cal L}^\pm_N$ such that $\tilde{\lambda}_L (0) = \tilde{l}_0$ and
$\tilde{\lambda}_L (1) \in \text{Im} ( \text{gr}_L )$. Denote by ${\cal CL}^\pm (X, \omega , l_0)$ the category of
based Lagrangian conductors of $(X, \omega , l_0)$. Its objects are elements 
$\tilde{L} = ((L_0 , \lambda_{L_0} ) , \dots , (L_l , \lambda_{L_l} ) ; J^{\tilde{L}})$  where  $(L_0 , \dots , L_l ; J^{\tilde{L}})$
is a Lagrangian conductor of $(X, \omega)$ and $(L_i , \lambda_{L_i})$ is a based Lagrangian submanifold of $(X, \omega , l_0)$,
$0 \leq i \leq l$, while its morphisms are continuations given by Definitions \ref{defcontinuation} and \ref{defbasedcontinuation}.

Let $\tilde{L} = ((L_0 , \lambda_{L_0} ) , \dots , (L_l , \lambda_{L_l} ) ; J^{\tilde{L}})$ be a based Lagrangian conductor.
Let $0 \leq i < j \leq l$, $x \in L_i \cap L_j $ and $u : \R_+ \times [-1,1] \to X$ be such that:

$\star 1)$  $u (\R_+ \times \{ -1 \}) \subset L_i $ and $u (\R_+ \times \{ 1 \}) \subset L_j$.

$\star 2)$ $u(0,t) = 
\left\{
\begin{array}{l}
\lambda_{L_j} (t) \text{ if } t \in [0,1],\\
\lambda_{L_i} (- t) \text{ if } t \in [-1,0].\\
\end{array} 
\right. $

$\star 3)$ $\lim_{+ \infty} u = x$.

Such a map $u$ does not exist in general for every intersection point  $x \in L_i \cap L_j $, only for those which are
homotopic to the concatenation $\lambda^{-1}_{L_i} \star \lambda_{L_j}$ is the space of paths
$\Omega (L_i , L_j)$ from $L_i$ to $L_j$. There is a tautological injection
$(x , [u]) \mapsto \lambda_{(x , [u])} \in \text{Ob} ( {\cal OS} (X, \omega))$, where $ [u]$ denotes a
homotopy class of such maps $u$ satisfying $\star 1$, $\star 2$, $\star 3$.
We then set $\widetilde{CF} (L_i , L_j) = \oplus_{(x , [u])} \lambda_{(x , [u])}$, the sum being infinite in general.
The construction of Floer functor then provides a functor $\widetilde{\cal F} : {\cal CL}^\pm (X, \omega , l_0)
\to K^b ({\cal OS} (X, \omega))$ where the chain complexes in the image of $\widetilde{\cal F}$ are infinite sums
of open strings in general. This functor induces on the quotient a functor
${\cal F} : {\cal CL}^\pm (X, \omega , l_0) \to K^b ({\cal OS}_1 (X, \omega))$ such that the image
${\cal F} (\tilde{L})$ coincides with the subcomplex of $(\widetilde{CF} (\tilde{L}) , \delta^{CF})$ which
contains open strings homotopic to the path $\lambda^{-1}_{L_i} \star \lambda_{L_j}$. This new functor
refines the preceding Floer functor ${\cal F}$.

Let us finally consider the special case of simply connected Lagrangian submanifolds. In this case, the second
homotopy group of $X$ transitively acts by concatenation on  homotopy classes $[u]$ of maps converging to a fixed point 
$x \in L_i \cap L_j $. This action writes $\lambda_{(x , \gamma.[u])} = \lambda_{(x , [u])}^{2c_1 (X , \omega).\gamma}$
where $\gamma \in \pi_2 (X)$. Floer differentials are equivariant for this action so that denoting by
$N = 2 \min_{\gamma \in \pi_2 (X) } \{ \vert c_1 (X , \omega).\gamma \vert \} \setminus \{ 0 \} \in \N^*$, the functor
$\tilde{\cal F} $ mods out to a functor $\tilde{\cal F}_N : {\cal CL}^\pm (X, \omega , l_0)
\to K^b ({\cal OS}_N (X, \omega))$, where  chain complexes in the image of $\tilde{\cal F}_N$ are now finite sums
of open strings. In this case, we denote by
$\A = \{Ê\int_\gamma \omega \in \R \, \vert \, \gamma \in \pi_2 (X) \}$ and by
$\Z ((t^\A)) = \{ \sum_{a \in \A} n_a t^a \, \vert \, \forall C \in \R, \, \# \{ a < C \, \vert \, n_a \neq 0 \in \Z \} < \infty \}$.
The composition of $\tilde{\cal F}_N$ with the functor coefficients ${\cal C}$ provides a functor
${\cal C} \circ \tilde{\cal F}_N : {\cal CL}^\pm (X, \omega , l_0 , \pi_1 = 0) \to {\cal M}od_{\Z ((t^\R))}$ with value
in the category ${\cal M}od_{\Z ((t^\R))}$  of free modules of finite type over the Novikov ring $\Z ((t^\R))$.
This functor turns out to be refined here to a functor with value in the category ${\cal M}od_{\Z ((t^\A))}$ of free 
modules of finite type over the Novikov ring $\Z ((t^\A))$. Indeed, for every $x \in L_i \cap L_j $ ,
${\cal C} (\lambda_x) = \Z_{\lambda_x} ((t^{ \{ -\int_{\R_+ \times [-1,1] } u^* \omega \in \R \, \vert \, \lambda_{(x , [u])} \in
\widetilde{CF} (L_i , L_j)  \} }))$ is a free module of rank one over $\Z ((t^\A))$.

When twice the first Chern class of $(X, \omega)$ vanishes in $H^2 (X ; \Z / N\Z)$, the functor $\tilde{\cal F}_N$ extends
to a true lift of Floer functor ${\cal F}$ as in the preceding \S \ref{subsubsectgraded}, such that ${\cal C} \circ \tilde{\cal F}_N$
remains with values in ${\cal M}od_{\Z ((t^\A))}$. We end this paragraph by discussing this phenomenon.

\begin{prop}
\label{propax}
Let $(X, \omega)$ be a symplectic manifold such that $2 c_1 (X, \omega) = 0 \in H^2 (X ; \Z / N\Z)$ and
${\cal L}^\pm_N \to {\cal L}^\pm$ be a ${\cal L}^\pm_N$-structure on $(X, \omega)$. Let
$(L_0 , \text{gr}_{L_0} ) , \dots , (L_l , \text{gr}_{L_l} )$ be $l+1$ graded simply-connected Lagrangian submanifolds
of $(X, \omega , {\cal L}^\pm_N)$ transversal to each other, $ l \geq 1$. There exists a function
$ x \in \sqcup_{0 \leq i < j \leq l} L_i \cap L_j \mapsto \A_x \in \R / \A$ such that for every map
$u : D \setminus \{ z_0 , \dots , z_q \} \to X$ satisfying

1) $\lim_{z_j} u = x_j \in L_{i_{j-1}} \cap  L_{i_j}$ for every $1 \leq q \leq l$, $0 \leq i_0 < \dots < i_q \leq l$ and $0 \leq j \leq q$.

2) $u(\partial_j D) \subset L_{i_j}$ for every $0 \leq j \leq q$,

the relation $\A_{x_0} = \int_D u^* \omega + \sum_{j=1}^q \A_{x_j} \in \R / \A$ holds, where
$\A = \{Ê\int_\gamma \omega \in \R \, \vert \, \gamma \in \pi_2 (X) \}$.
\end{prop}

{\bf Proof:}

Let us fix a base point $x_0 \in X$ and a Lagrangian subspace $\tilde{l}_0 \in {\cal L}^\pm_N \vert_{x_0}$ of $T_{x_0} X$,
together with a path $\lambda_i : [0,1] \to {\cal L}^\pm$ which satisfies $\lambda_i (0) = \tilde{l}_0$ 
and $\lambda_i (1) \in \text{Im} (gr_{L_i})$,
so that for every $0 \leq i \leq l$, $(L_i , \text{gr}_{L_i}  ,\lambda_i)$ is based graded. Let us choose a finite system
of loops $\underline{\gamma}$ based at $x_0$ which generate the fundamental group $\pi_1 (X ; x_0)$.
This choice provides a surjective morphism from a free group $F$ of finite type onto $\pi_1 (X ; x_0)$.
The kernel of this morphism is the subgroup $R$ of relations, it is equipped with a morphism $a : R \to  \R / \A$
defined as follows. Every relation $r \in R$ is represented by a combination of loops in $\underline{\gamma}$
which bounds a disk $D_r$ of $X$.  We set $a(r) = \int_{D_r} \omega \in  \R / \A$. Since $ \R / \A$ is Abelian,
this morphism quotients out to a morphism $R^{ab} \to  \R / \A$ defined on the Abelianization of $R$.
Since $R^{ab}$ is a subgroup of the free Abelian group $F^{ab}$ of finite type and since $ \R / \A$  is divisible,
the morphism $R^{ab} \to  \R / \A$ extends to a morphism $F^{ab} \to  \R / \A$ still denoted by $a$.
Now, for every $x \in L_i \cap L_j$, there exists a word $\gamma_x$ in the alphabet $\underline{\gamma}$
such that the concatenation $\lambda^{-1}_{L_i} \star  \gamma_x \star \lambda_{L_j}$ is homotopic
to $x$ in the space of paths $\Omega (L_i , L_j)$ from $L_i$ to $L_j$. Let us choose such a function
$x \in L_i \cap L_j \mapsto  \gamma_x$  which is constant on every connected component of  $ \Omega (L_i , L_j)$.
We then set, for every $x \in L_i \cap L_j$, $ \A_x = - \int_{\R_+ \times [-1 , 1]} u_{i, j}^* \omega - a(\gamma_x) \in \R / \A$, 
where $u_{i, j}$ is  homotopy between $\lambda^{-1}_{L_i} \star  \gamma_x \star \lambda_{L_j}$
and $x$ satisfying conditions $\star 1$, $\star 2$ and $\star 3$. This function satisfies the properties
of Proposition \ref{propax}. $\square$

\begin{rem}
Under the hypothesis of Proposition \ref{propax}, consider the two-dimensional $CW$-complex
whose vertices are labeled by $0 , \dots , l$, whose edges between two vertices $i$ and $j$, $0 \leq i < j \leq l$,
are labeled by homotopy classes of paths from $L_i$ to $L_j$ and whose faces are labeled by homotopy
classes of polygons with $q+1$ vertices $0 \leq i_0 < \dots < i_q \leq l$, $1 \leq q \leq l$, which are in cyclic
order on the boundary and such that the concatenation of the $q+1$ paths given by the edges of the 
polygon is homotopic to the constant path in $X$. The difference between two functions given by
Proposition \ref{propax}  reads as a one-cocycle on this $CW$-complex taking value in $\R / \A$.
Moreover, the function we have constructed in the proof of this Proposition \ref{propax} does not depend on the
choice of the base point $x_0$ or on the system $\underline{\gamma}$ of loops, it only depends on the choice
of paths  $\lambda_i $, $0 \leq i \leq l$, and on the choice of the extension $a$ to $F^{ab}$ when
$b_1 (X) = \dim H_1 (X ; \Q) \neq 0$ (compare Corollary \ref{corb1}). The difference between two functions
constructed in the proof of this Proposition \ref{propax} reads as a one-coboundary on our $CW$-complex 
taking value in $\R / \A$.
\end{rem}

\section{Manifolds with vanishing first Chern class}
\label{parc1}

This paragraph is devoted to closed or convex at infinity symplectic manifold $(X, \omega) $ 
of dimension $2n \geq 4$ for which the first Chern class vanishes.
Its aim is to prove that their category ${\cal CL}^\pm (X, \omega) $ of Lagrangian conductors
contains vanishing cycles in the sense of Definition \ref{defvanishingcycle}. In particular, integral Floer cohomology
of Lagrangian spheres in such manifolds happens to be well defined. As was pointed out in the note added in
proof of \cite{Kont}, Floer cohomology in such manifolds is obstructed. A systematic study of this obstruction
is carried out in \cite{FOOO}. Our approach to overcome it relies on a phenomenon of localization of pseudo-holomorphic
disks with boundary on Lagrangian spheres. To observe this phenomenon, a key use is made of symplectic field theory.
This phenomenon was observed in \cite{WelsSFT} and is a particular case of sharpness results obtained there, we reproduce it
in the first paragraph. The second paragraph is devoted to the compactness theorem which is needed here and
the third one to some study of the Floer cohomology obtained.

\subsection{Localization of pseudo-holomorphic membranes}

We recall and adapt to our need here Theorem $1.6$ of  \cite{WelsSFT}. Let us fix integers $k,p$ such that
$2 < p < + \infty$ and $1 \ll k$, where $k$ is also supposed to be much less than the regularity of our
almost complex structures. Let $L$ be a smooth Lagrangian sphere in a closed or convex at infinity symplectic 
manifold $(X, \omega) $ with vanishing first Chern class. We denote by ${\cal M}_{g,v} (X , L)$ the space of
triple $(u , C , J)$ where $J \in {\cal J}_\omega^\infty (L)$, $C$ is a punctured Riemann surface of genus $g$
with $v$ punctures and $u : C \to X \setminus U^*L$  is a proper simple $J$-holomorphic map having $k$ derivatives in $L^p$
and which has finite Hofer energy. Here, $U^*L = \phi (U^* S^n)$ for an $A_1$-neck $\phi$ associated to $J$,
see Definition \ref{defa1singular}, and triples $(u , C , J)$ are considered only up to reparametrization by an automorphism of $C$.
Hofer energy is introduced in \cite{HWZ1}, its finiteness  prescribes the behavior of the $J$-holomorphic map
$u$ near the punctures. Namely, near each puncture, $u$ has to be asymptotic to  a cylinder $u(\gamma \times ]0 , \epsilon])$
over a closed Reeb orbit of $R_\lambda$, see Theorem $2.8$ of \cite{HWZ1} and its adaptation to the Morse-Bott case
in \cite{Bour}. Let ${\cal M}^\infty (X , L) = \sqcup_{v \geq 1} {\cal M}_{0,v} (X , L)$ if the dimension of $X$ is four and
${\cal M}^\infty (X , L) = \sqcup_{v \geq 1} \sqcup_{g \geq 0} {\cal M}_{g,v} (X , L)$ if its dimension is greater than four.
Let $\pi^\infty : (u , C , J) \in {\cal M}^\infty (X , L) \mapsto J \in {\cal J}_\omega^\infty (L)$. The space ${\cal M}^\infty (X , L)$
is a separable Banach manifold of finite regularity, namely the difference between $k$ and the regularity of our almost complex
structures, and $\pi^\infty$ is Fredholm, see Proposition $3.2.1$ of \cite{McDSal}, Theorem $2.8$ of \cite{HWZfred} and its adaptation
to the Morse-Bott case in \cite{Bour}. Let $\mu_T^L \in H^2 (X,L ; \Z)$ be the Maslov class of $L$ and $H \subset X$
be a possibly empty codimension two closed symplectic submanifold of $X$ disjoint from $L$. We denote by 
${\cal J}_\omega^\infty (L , H)$ the space of $J \in {\cal J}_\omega^\infty (L)$ such that $H$ is $J$-holomorphic and
does not meet the singular locus of $J$. We denote by ${\cal M}^\infty (X , L , H) = (\pi^\infty)^{-1} ({\cal J}_\omega^\infty (L , H))$
and by $\pi^\infty_H : {\cal M}^\infty (X , L , H) \to {\cal J}_\omega^\infty (L , H)$ the restriction of $\pi^\infty$.
Let $H^* : H_2 (X , L ; \Z) \to \Z$ be the morphism of intersection with $H$.

\begin{theo}
\label{theolocalization}
Let $L$ be a smooth Lagrangian sphere in a closed or convex at infinity symplectic 
manifold $(X, \omega) $ of dimension $2n \geq 4$.
Let $H \subset X$ be a possibly empty codimension two closed symplectic submanifold disjoint from $L$.
Assume that there exist non-negative real numbers $a,b$ such that
$\mu_T^L + a \omega + b H^* = 0 \in \text{Hom} (H_2 (X , L ; \Z) , \Z)$.
Then, the Fredholm index of $\pi^\infty_H$ is everywhere bounded from above by $-2$.
\end{theo}
Note that only the case where $a=g=0$ and $H = \emptyset$ will be used in this paper.\\

{\bf Proof:}

Let $(u , C , J) \in {\cal M}^\infty (X , L , H) $. In the Morse-Bott set-up used here, the Fredholm
index of the restriction of $\pi^\infty_H$ to the connected component of $(u , C , J)$ has been
computed by F. Bourgeois in \cite{Bour}. This index is equal to
$ \mu_T^{S^* L} \circ u(C) + (n-3)(2-2g) + 2v$, see Proposition $1.12$ of \cite{WelsSFT},
where $ \mu_T^{S^* L} \circ u(C) $ is twice the obstruction to extend over the whole
$u(C)$ the canonical trivialization of $TX$ given near the punctures by the Reeb flow. Now,
for every orbit $\gamma$ of $R_\lambda$, the open subset $U^* L$ contains a symplectic
plane $P_\gamma$ of finite energy converging to $\gamma$ and with index
$\mu_T^{S^* L} (P_\gamma )  = 2(n-1)$. If we equip  the interior of  $U^* L$
with the complex structure coming from the affine ellipsoid $\{ x_0^2 + \dots + x_n^2 = 1 \} \subset \C^{n+1}$,
then this plane $P_\gamma$ can be chosen to be a complex line of this ellipsoid.
At every puncture, the map $u$ travels around an integral number of times the orbit $\gamma$
in the image, this integer is called the multiplicity of the orbit. Let $m_1 , \dots , m_v$ be the
multiplicities of the orbits $\gamma_1 , \dots , \gamma_v$ in the image of the punctures of $C$ and
$S \subset X$ be the symplectic surface obtained as the union of $u(C)$ with coverings of
degree $m_i$ of $P_{\gamma_i}$, $1 \leq i \leq v$. From the hypothesis follows that
$\mu_T^L (S) = - a \int_S \omega - b H^* (S) \leq 0$ so that
$ \mu_T^{S^* L} \circ u(C) \leq -2(n-1) \sum_{i=1}^v m_i \leq -2(n-1) v$.
Hence,

$\begin{array}{rcl}
\ind (\pi^\infty_H) & = &  \mu_T^{S^* L} \circ u(C) + (n-3)(2-2g) + 2v\\
 & \leq & (n-3)(2-2g - 2v) - 2v\\
& \leq & -2 \text{ since } n \geq 3 \text{ or } g=0. \quad \square
\end{array}$ 

\begin{rem}
In the framework of complex algebraic geometry, Theorem \ref{theolocalization} can be interpreted 
in the following way. Let $X$ be a $n$-dimensional projective variety having a singular point
$x$ of type $A_1$ and an effective canonical class $K_X$. Let $\pi : Y \to X$ be the blow-up
of the singular point $x$ so that $Y$ is smooth and $K_Y = \pi^* K_X  + (n-2) {\cal O}_Y (E)$,
where $E$ is the exceptional divisor of the blow-up. The Riemann-Roch index of a curve $C$
in $Y$ thus writes 

$\begin{array}{rcl}
\ind_\C (\pi) & = &  -K_Y . C + (n-3)(1-g) \\
 & \leq & -(n-2) E. C + (n-3)(1-g) \\
& \leq & -1 \text{ if } E . C > 0 \text{ and either } n \geq 3 \text{ or } g=0. 
\end{array}$ 

Hence, the real index of a curve $C$ of $X$ passing through $x$ is bounded from above by $-2$.
\end{rem}

\begin{cor}
\label{cor<E}
Under the hypothesis of Theorem \ref{theolocalization}, for every $E>0$, there exists a
dense Baire subset ${\cal B}_E$ of the second category of $C^l ([0,1] , {\cal J}_\omega^\infty (L , H))$
with the following property. For every $\theta \in {\cal B}_E$, there exists $\eta > 0$ such that
for every almost-complex structure $J \in \overline{\cal J}_\omega$ in the $\eta$-neighborhood of
$\text{Im} (\theta)$, there is no compact $J$-holomorphic curve $S$ in $X$ such that
$\int_S \omega \leq E$ and either $\partial S \subset L$ or $S \cap L \neq \emptyset$ if
$S$ is closed.
\end{cor}

{\bf Proof:}

This follows from Theorem \ref{theolocalization} and the compactness Theorem \cite{BEHWZ} 
in symplectic field theory. $ \square$

\begin{rem}
\label{remcodim2}
As soon as $n \geq 3$, the upper bound given in Theorem \ref{theolocalization} is reached only
when $v=1$, $a,b = 0$ and the multiplicity of the orbit $\gamma$ in the image of the puncture is
one. Let us equip  the interior of  $U^* L$
with the complex structure coming from the affine ellipsoid $Q^n = \{ x_0^2 + \dots + x_n^2 = 1 \} \subset \C^{n+1}$
in such a way that $L$ coincides with the real locus $\R Q^n$. We then see that the boundaries of 
once punctured holomorphic disks sitting on $L$ with multiplicity one at the image $\gamma$ of the puncture foliate the complement
of a codimension two equator $S^{n-2} \subset L$. This space of disks can be compactified by adding
complex lines passing through the point $\gamma \in Q^n$ and meeting $L$.
Denote by ${\cal M}_1 (Q^n , \gamma ; J)$ the space of such holomorphic disks having one marked point on the boundary.
The evaluation map at the marked point from this space to $L$ is a pseudo-cycle of degree $\pm 1$,
see \S $1.3$ of \cite{McDSal}. Here $J$ denotes the complex structure of $Q^n$, but the same result holds
for any almost-complex structure of $Q^n$ cylindrical at infinity and hence for any $A_1$-singular almost-complex
structure of $X$. If we replace the one-parameter family of Corollary \ref{cor<E} by a two-parameter family
$\theta : [0,1]^2 \to  {\cal J}_\omega^\infty (L , H)$ and denote by ${\cal M}_1 (X , \theta)$ the space of split $J$-holomorphic
disks of minimal energy with one marked point on the boundary and $J \in \text{Im} (\theta)$, then
${\cal M}_1 (X , \theta)$ may this time be non-empty and localized at the image of finitely many parameters
$J_i = \theta (t_i)$. For every such parameter, the evaluation map is then a pseudo-cycle of degree $\pm 1$.
Assume that this parameter is unique and perturb $\theta$ to a map with image in ${\cal J}_\omega \subset \overline{\cal J}_\omega$.
The evaluation map remains then a pseudo-cycle of degree $\pm 1$, but ${\cal M}_1 (X , \theta)$ is no more
localized at one almost-complex structure.
\end{rem}

\begin{cor}
\label{corvanishingcycle}
Let $(X, \omega) $ be a closed or convex at infinity symplectic 
manifold of dimension $2n \geq 4$ which has vanishing first Chern class.
Then every vanishing cycle $(L , J_L)$ where $L$ is a smooth Lagrangian
sphere of $(X, \omega) $ and $J_L \in {\cal J}_\omega^\infty (L) \setminus \text{Im} (\pi^\infty)$
is an elementary Lagrangian conductor of $(X, \omega) $. $\square$
\end{cor}
Note that since for every $n \geq 2$, spheres of dimension $n$ have vanishing fundamental
group and second Stiefel-Whitney class, they always carry a unique $\widetilde{GL}_n^\pm (\R)$-structure
so that mention of this structure can be omitted. Corollary \ref{corvanishingcycle} means that as
soon as $(X, \omega) $ contains Lagrangian spheres, the category ${\cal CL}^\pm (X, \omega) $
is non-empty. Note that from the work of Paul Seidel, Lagrangian spheres sometimes split generate
the whole Fukaya category, see \cite{SeidICM}, \cite{Seidquar}, \cite{Seidbook}.

\subsection{Compactness Theorem}
\label{subsectcompactness}

The combinatorial type of a $(l+1)$-punctured stable disk of $K_l$ is encoded by a connected
tree having $l+1$ free edges labeled by $z_0 , \dots , z_l$ and whose vertices are of valence at least
three.

\begin{defi}
A $(l+1)$-punctured prestable disk is a $(l+1)$-punctured nodal disk whose combinatorial
type is encoded by a connected
tree having $l+1$ free edges labeled by $z_0 , \dots , z_l$ and whose vertices are of valence at least
two. In addition, such a prestable disk may have trees of Riemann spheres attached
to its interior points. The associated stable disk is obtained by contraction of the unstable components of the disk.
\end{defi}

The combinatorial type of the stable disk associated to an unstable one is thus obtained by
contraction of the bivalent vertices of the combinatorial type of the unstable one.
Let $D$ be such a $(l+1)$-punctured prestable disk and 
$J \in \overline{\cal J}_\omega (U_l , V_l)$, where $V^{\tilde{L}}_l$ is a coherent choice of strip like
ends given by Definition \ref{defstriplikeend}. Then, $J$ induces a map $D \to \overline{\cal J}_\omega$.
Indeed, every unstable component of $D$ is isomorphic to the strip $\R \times [ -1 , 1 ]$ or to a Riemann sphere.
Every maximal connected chain of strips is attached to a puncture $\sigma (\underline{z})$ of $D_{\underline{z}}$,
where $D_{\underline{z}}$ is the stable disk associated to $D$, $\underline{z} \in K_l$. Property $P_2$ satisfied by $J$,
see \S \ref{pardefassociahedron}, ensures that in the  strip like end
 $\psi_{\sigma (\underline{z})}$ given by $V_l$, the composition $J \circ \psi_{\sigma (\underline{z})} :  \R_+ \times [-1 , 1] \to
 \overline{\cal J}_\omega$ does not depend on the first factor $\R_+$. Hence it extends to the chain of strips
 by $ (\tau , t) \in \R \times [ -1 , 1 ] \mapsto J_\sigma (t) \in  \overline{\cal J}_\omega$, where 
  $J_\sigma = J \circ \psi_{\sigma (\underline{z})} \vert_{\{ 0 \} \times  [-1 , 1]}$ is the associated path.
  Since by Property $P_2$ the latter does not depend on
 the choice of the strip like end $\psi_{\sigma (\underline{z})}$, this definition is consistent. 
 On every maximal connected tree of Riemann spheres, $J$ is constant with value given by the point of the disk
 where it is attached.
 
 \begin{theo}
 \label{theocompactness}
 Let $(X, \omega) $ be a closed or convex at infinity symplectic manifold of dimension $2n \geq 4$ such that 
$c_1 (X, \omega) + a \omega = 0$, $a \in \R_+$.
Let $H$ be an effective continuation from the Lagrangian conductor $\tilde{L} = (L_0 , \dots , L_{l} ; J^{\tilde{L}})$ to 
$\tilde{L}' = (L'_0 , \dots , L'_l ; J^{\tilde{L}'})$. Let $J^H \in \overline{\cal J}_\omega (\overline{U}_{J_l} , V_{J_l})$ 
be a generic extension of $J^{\tilde{L}}$, $J^{\tilde{L}'}$ given by Proposition \ref{propbordism}.
Let $(u_i , \underline{z}^i , s_i)_{i \in \N}$ be a sequence such that $(s_i , \underline{z}^i) \in J_l$,
$u_i : D_{\underline{z}^i} \to X$ satisfies the Cauchy-Riemann equation 
$J^H \vert_{(s_i , \underline{z}^i)} \circ du_i = du_i \circ J_{D_{\underline{z}^i}}$ 
and such that  for every $ w \in \partial_j D_{\underline{z}^i}$, $u(w) \in L_{j}^{h_j (w)}$
where $h_j$ is the map given by Property $1$ of Proposition \ref{propbordism}, compare \S \ref{subsubsectionfloercontinuations}.
Assume that the Lagrangian submanifolds $L_0 , \dots , L_l $ are spheres. Then, after extracting a subsequence of
$i \in \N$ if necessary, the following two properties hold

1) $(s_i , \underline{z}^i)_{i \in \N}$ converges to a point $(s_\infty , \underline{z}^\infty)$ of $J_l$.

2) $(u_i )_{i \in \N}$ converges to a stable map $u_\infty : D \to X$ in Gromov-Floer topology, where
$D$ is a prestable disk with associated stable disk $D_{\underline{z}^\infty}$.                               
\end{theo}

The key point in Theorem \ref{theocompactness} is the absence of disk bubbles attached to the boundary of
the curve in the limit. This absence results from Theorem \ref{theolocalization}.
Gromov-Floer topology is defined in \cite{McDSal}.\\

{\bf Proof:}

Assume for simplicity that $(J^H)^{-1} (\partial \overline{\cal J}_\omega ) $ is contained in the boundaries
of the fibers of the projection $\overline{U}_{J_l} \to J_l$, since only this case is used here.
Since for every $0 \leq j \leq l$, the Hamiltonian isotopy class of $L_j^{s_i}$ does not depend on $i \in \N$,
the energy $\int_{D_{\underline{z}^i}} u_i^* \omega$ is bounded, see \cite{Oh1}. The restriction of
$u_i $ to $\partial D_{\underline{z}^i}$ has then bounded derivative. Indeed, if there would exist a point
$w \in \partial D_{\underline{z}^i}$ such that the derivative of $u_i $ at $w$ diverges, then Lemma $10.7$
of \cite{BEHWZ} would provide a  bubble 
attached to the boundary of the Gromov-Floer limit of our curve. Such a bubble would be 
$J^H \vert_{ (s_\infty , \underline{z}^\infty)}$-holomorphic and have several levels in the sense of Theorem
$10.3$ of \cite{BEHWZ}. Theorem \ref{theolocalization} prevents the existence of such a bubble. 
Theorem \ref{theocompactness} now follows from the classical compactness Theorem in Floer theory,
see \cite{Floer}, \cite{Oh1}, \cite{Fuk2}, \cite{Frauen}. $\square$ 

\begin{cor}
Let $(X, \omega) $ be a closed or convex at infinity symplectic manifold of dimension $2n \geq 4$ 
with vanishing first Chern class. Let $\tilde{L} = (L_0 , \dots , L_{l} ; J^{\tilde{L}})$ be a 
Lagrangian conductor made of vanishing cycles given by Corollary \ref{corvanishingcycle}. Let 
$\lambda^+ , \lambda^- \in CF (\tilde{L})$ be open strings such that
$\mu (\lambda^-) - 1 -\mu (\lambda^+) + q^+ \leq 1$, where $q^+$ is the cardinality of $\lambda^+$
and $\lambda^-$ is elementary. Let 
$(\gamma_i )_{i \in \N}  $ be a sequence of trajectories from $\lambda^+ $ to $\lambda^-$.
Then any limit curve given by Theorem \ref{theocompactness} at most has two irreducible components,
none of which is spherical.
\end{cor}

{\bf Proof:}

Let $u_\infty : D \to X$ be a limit curve given by Theorem \ref{theocompactness}. If the restriction
of $u_\infty$ to a spherical component of $D$ is multiple, let us replace it by the underlying simple map.
Then, from the index computation in Proposition \ref{propindex} or \cite{Fuk2},  the index of
$u_\infty$ is negative unless it has no spherical component and at most two irreducible components. $\square$

\subsection{Floer cohomology and Donaldson category}

Let $(X, \omega) $ be a closed or convex at infinity symplectic manifold of dimension $2n \geq 4$ 
with vanishing first Chern class. Let $\tilde{L} = (L_0 , \dots , L_{l} ; J^{\tilde{L}})$ be a 
Lagrangian conductor made of vanishing cycles given by Corollary \ref{corvanishingcycle}.
Then ${\cal C} \circ {\cal F} (\tilde{L} )$ is a complex of free modules of finite type over the
Novikov ring $\Z ((t^\R))$. 

\begin{defi}
\label{defFloercohomology}
The cohomology of the complex ${\cal C} \circ {\cal F} (\tilde{L} )$  is called Floer cohomology and denoted
by $HF (\tilde{L})$.
\end{defi}

From \S \ref{subsubsectgraded}, when $(X, \omega)$ is equipped with a ${\cal L}^\pm_N$-structure
and the Lagrangian spheres are graded, the Floer cohomology $HF (\tilde{L})$ inherits a non-trivial grading
modulo $N \in \N$. From \S \ref{subsubsectbased},  this complex can be refined to 
a complex of free modules of finite type over the Novikov ring $\Z ((t^\A))$, where 
$\A = \{Ê\int_\gamma \omega \in \R \, \vert \, \gamma \in \pi_2 (X) \}$, it suffices to equip the Lagrangian
spheres with based paths. This Definition \ref{defFloercohomology} extends to vanishing cycles which 
are not necessarily transversal to each other. Indeed, if $\tilde{L} = (L_0 , \dots , L_{l} ; J^{\tilde{L}})$ is such that
the spheres $L_i$ are not transversal to each other, a small Hamiltonian perturbation makes them
transversal to each other thus defining a Lagrangian conductor $\tilde{L}^\epsilon = (L_0^\epsilon , \dots , L_{l}^\epsilon ; J^{\tilde{L}})$.
Small Floer continuations then provide canonical isomorphisms between the Floer cohomologies $HF (\tilde{L}^\epsilon)$
 obtained from different small values of $\epsilon > 0$ or different choices of Hamiltonian perturbations.
 We define $HF (\tilde{L})$ to be this class of modules $HF (\tilde{L}^\epsilon)$ up to canonical isomorphisms.

\begin{rem}
\label{remmonodromy}
For arbitrary Hamiltonian perturbations, that is Hamiltonian perturbations $(L_i^\epsilon)_{\epsilon >0}$ such
that $L_i^\epsilon$ does not stay in the neighborhood $U^* L_i$ given by the $A_1$-singular almost complex structure
$J^{\tilde{L}}_i$, see \S \ref{subsectionsingular}, it is no more possible to keep the map $J^{\tilde{L}}$ to
define a Lagrangian conductor $\tilde{L}^\epsilon$. 
Floer continuations then still provide isomorphisms between the Floer cohomologies $HF (\tilde{L}^\epsilon)$,
but these isomorphisms are no more canonical a priori. Theorem \ref{theolocalization} instead suggests that
the space ${\cal J}_\omega^\infty (L) \setminus \text{Im} (\pi^\infty)$ is not simply connected.
Remark \ref{remcodim2} would then show that the fundamental group of ${\cal J}_\omega^\infty (L) \setminus \text{Im} (\pi^\infty)$
acts non-trivialy by monodromy on the automorphisms of $HF (\tilde{L})$. It is the local simply-connectedness of
${\cal J}_\omega^\infty (L) \setminus \text{Im} (\pi^\infty)$ which implies that this phenomenon does not occur 
locally so that $HF (\tilde{L})$ can be defined even for non-transversal vanishing cycles.
\end{rem}

From Theorem \ref{theocomplex}, the product $m_2 : HF( L_0 , L_1 ; J^{01}) \otimes HF( L_1 , L_2 ; J^{12}) \to
HF( L_0 , L_2 ; J^{02}) $, introduced by Donaldson in early $1990$'s,  is associative.

\begin{prop}
\label{propidentity}
Let $(X, \omega) $ be a closed or convex at infinity symplectic manifold of dimension $2n \geq 4$ 
with vanishing first Chern class. Let $(L , J_L )$ be a vanishing cycle given by Corollary \ref{corvanishingcycle}. 
Then, there exists an element $e \in HF (L , L ; J_L )$ such that for every Lagrangian conductor 
$\tilde{L}' = (L , L' ; J^{\tilde{L}'})$ (resp. $\tilde{L}' = (L' , L ; J^{\tilde{L}'})$) such that $J^{\tilde{L}'}_0 = J_L$
(resp. $J^{\tilde{L}'}_1 = J_L$),
$m_2 ( e , *) = id_{HF (\tilde{L}')}$ (resp. $m_2 ( * , e ) = id_{HF (\tilde{L}')}$).
\end{prop}

{\bf Proof:}

Let $f : L \to \R$ be a Morse function having only two critical points outside the
intersection $L \cap L'$. For every small positive $\epsilon$, let $L_\epsilon$ be the graph
of the derivative of $\epsilon f$ which lies in the neighborhood $U^* L$ of $L$ given 
by the $A_1$-singular almost complex structure $J_L$. Since $\epsilon$ is small, 
$L_\epsilon$ is transversal to $L'$ and there is an obvious bijection between 
$L_\epsilon \cap L'$ and $L \cap L'$. Let $x_{\max} \in L \cap L_\epsilon$ be the maximum of
the function $f$ and $\lambda_{\max}^\epsilon$ be the corresponding elementary open string.
The index of this string vanishes from Example $2$ of \S \ref{subsectionopenstring}
and we will prove that its class in $HF (L , L ; J_L )$ has the required properties.
Let $\lambda^+ \in CF (L' , L ) $ and $\lambda^-_\epsilon \in CF  (L'  , L_\epsilon )$
be elementary strings having same index. Let $(\gamma_\epsilon)_{\epsilon >0}$ be a continuous
family of trajectories $\lambda^+  \otimes \lambda_{\max}  \to \lambda^-_\epsilon$ counted by $m_2$.
We may assume that $\lambda^-_\epsilon$ converges to a string $\lambda^-_0$ as
$\epsilon$ converges to zero. Since $\mu (\lambda^-_0) = \mu (\lambda^+)$, there cannot be
any Floer trajectory from $\lambda^+$ to $\lambda^-_0$ and thus from the compactness Theorem
\ref{theocompactness}, $\lambda^-_0$ and $\lambda^+$ have to coincide. Likewise, the
energy of the trajectory $\gamma_\epsilon$ has to converge to zero as $\epsilon$ converges to zero,
since otherwise the limit curve would contain a $J_L$ holomorphic disk with boundary in $L$,
from the removal of singularities for disks \cite{Ohremov}, which does not exist from 
Theorem \ref{theolocalization}. The trajectory $\gamma_\epsilon$ is thus contained in the
neighborhood $U^* L$ of $L$ as soon as $\epsilon$ is small enough. As a result, 
we may assume that $X = U^* L$ and that $L'$ is a fiber of this bundle and we have to find
a $J^{\tilde{L}'}$ such that only one trajectory goes from $\lambda^+ \otimes \lambda_{\max}$ to $\lambda^-_\epsilon$.
This  follows from \cite{FukOh}, but let us propose here a direct proof in our special case. 
So far, we didn't use any property of $J^{\tilde{L}'}_1$, we may thus assume that $J^{\tilde{L}'}$ takes values
in the set of almost complex structures of $X = U^* L$ cylindrical at infinity.
Let us identify $L$ with the unit sphere $S^n \subset \R^{n+1}$,  $f$ with the height function
$(x_1 , \dots , x_{n+1}) \in S^n \mapsto  x_{n+1} \in \R$ and $L'$ with $T^*_{(1,0, \dots , 0)} S^n$.
Let $\text{rot} : LÊ\to L$ be the isometry $(x_1 , \dots , x_{n+1}) \in S^n \mapsto (x_1 , -x_2 , \dots , -x_{n+1}) \in S^n$.
It fixes the point $(1,0, \dots , 0)$ and exchanges the minimum and maximum of $f$.
Let us identify now $D_{\underline{z}}$ with $(\R_+ \times [-1 , 1]) \setminus \{ (0 , -1) , (0 , 1) \}$ where 
$z_1 = (0 , -1) $, $z_2 = + \infty$ and $z_0 = (0 , 1)$. Let $J^{\tilde{L}'} : (\tau , t) \in \R_+ \times [-1 , 1] 
\mapsto J_t \in {\cal J}_\omega^\infty (L) $, where $J_t $ is the almost-complex structure given by Floer
in \cite{FloerWitten}, \S $5$, $p218$. A trajectory $\gamma_\epsilon : \lambda^+ \otimes  \lambda_{\max}  \to \lambda^-_\epsilon$
then extends to a $J_t$-holomorphic map $u_ \epsilon : (\R \times [-1 , 1]) \setminus \{ (0 , -1) , (0 , 1) \} \to X$
by the formula $u_ \epsilon (\tau , t) = c_L \circ \text{Rot} \circ u (-\tau , t)$, where $\text{Rot} : U^* L \to U^* L$
is induced by $\text{rot} : LÊ\to L$ and $c_L : (q,p) \in U^* L \mapsto (q,-p) \in U^* L$. 
From the removal of singularities for disks \cite{Ohremov}, it finally extends to a trajectory 
$\gamma_\epsilon : \lambda_{\max} \to \lambda_{\min}$. Since there is only one geodesic of $S^n$
going from the maximum of $f$ to its minimum passing through $(1,0, \dots , 0)$, we deduce from 
\S $5$ of  \cite{FloerWitten} that we get only one such trajectory. Hence the result. $\square$

\begin{cor}
Under the hypothesis of Proposition \ref{propidentity}, $HF (L , L ; J_L )$ is canonically isomorphic
to $H^* (S^n , \Z) \otimes  \Z ((t^\R))$ as a $\Z ((t^\R))$-algebra.
\end{cor}

{\bf Proof:}

As a module over $\Z ((t^\R))$, $HF (L , L ; J_L )$ is of rank two, with one generator of index zero
and one of index $n$. The latter can be chosen to be given by the maximum and minimum of
a Morse function having only two critical points, as in the proof of Proposition \ref{propidentity}.
The product on $HF (L , L ; J_L )$ is deduced from Proposition \ref{propidentity}, namely
$m_2 (\lambda_{\max} , \lambda_{\max}) = \lambda_{\max}$, $m_2 (\lambda_{\max} , \lambda_{\min}) = \lambda_{\min}$
and $m_2 (\lambda_{\min} , \lambda_{\min}) = 0$ since there is no string of index $2n$ in $HF (L , L ; J_L )$. $\square$

\begin{cor}
\label{corb1}
Let $(X, \omega) $ be a closed or convex at infinity symplectic manifold of dimension $2n \geq 4$ 
with vanishing first Chern class. Then,

1) No smooth Lagrangian sphere of $(X, \omega) $ can be displaced from itself by a Hamiltonian
isotopy.

2) As soon as $m \geq 1$, $(X \times \C^m , \omega \oplus \omega_{std})$ does not contain any
Lagrangian sphere. $\square$
\end{cor}

The first part of Corollary \ref{corb1} is a particular case of a result obtained by
Fukaya, Oh, Ohta and Ono with different methods, see the new version of \cite{FOOO}. Recall \cite{Beau} that a projective Calabi-Yau 
manifold with non-vanishing first Betti number has a covering of the form $X \times \C^m$, $m \geq 1$,
so that from Corollary \ref{corb1} it does not contain any Lagrangian sphere. This result has been observed by
Paul Biran and kindly communicated to me. Not that the product $\C P^n \times \C^{n+1}$ does contain Lagrangian
sphere, see \cite{BirCie}.

Let us finally denote by $ \text{Ob} ( \widetilde{{\cal L}ag}_0 (X, \omega))$ the set of (graded) vanishing cycles 
given by Corollary \ref{corvanishingcycle} and by $ \text{Hom} ( \widetilde{{\cal L}ag}_0 (X, \omega))$
the set of Floer cohomology modules between two such objects. Let $\text{For}_J : (L , J_L ) \in
 \text{Ob} ( \widetilde{{\cal L}ag}_0 (X, \omega)) \mapsto L \in {\cal L}ag_0 (X, \omega)$ and
 $\text{For}_L : (L , J_L ) \in
 \text{Ob} ( \widetilde{{\cal L}ag}_0 (X, \omega)) \mapsto J_L \in \cup_{L \in {\cal L}ag_0 (X, \omega)}
 ({\cal J}_\omega^\infty (L) \setminus \text{Im} (\pi^\infty))$. Here, ${\cal L}ag_0 (X, \omega)$  denotes the
 set of smooth Lagrangian spheres of $(X, \omega)$ ; a section of $\text{For}_J $ is a map
 $\sigma : {\cal L}ag_0 (X, \omega) \to  \text{Ob} ( \widetilde{{\cal L}ag}_0 (X, \omega))$ such that
 $\text{For}_J  \circ \sigma =  id$. We denote by $\widetilde{{\cal L}ag}_0 (X, \omega)$ the pair
 $\big( \text{Ob} ( \widetilde{{\cal L}ag}_0 (X, \omega)) , $ $ \text{Hom} ( \widetilde{{\cal L}ag}_0 (X, \omega)) \big)$.

\begin{theo}
\label{theoDonaldson}
Let $(X, \omega) $ be a closed or convex at infinity symplectic manifold of dimension $2n \geq 4$ 
whose first Chern class vanishes. Then, $\widetilde{{\cal L}ag}_0 (X, \omega)$
has the structure of a small preadditive category. Moreover, if $\sigma $, $\sigma'$ are two sections
of $\text{For}_J $, the images of $\sigma $, $\sigma'$ are isomorphic subcategories of 
$\widetilde{{\cal L}ag}_0 (X, \omega) $. $\square$
\end{theo}

\begin{cor}
Under the hypothesis of Theorem \ref{theoDonaldson},  $\widetilde{{\cal L}ag}_0 (X, \omega) $ induces
a structure of small  preadditive  category on $\cup_{L \in {\cal L}ag_0 (X, \omega)}
 ({\cal J}_\omega^\infty (L) \setminus \text{Im} (\pi^\infty))$ and a structure of isomorphism class of
 small preadditive  category on ${\cal L}ag_0 (X, \omega)$. $\square$
\end{cor}

\addcontentsline{toc}{part}{\hspace*{\indentation}Bibliography}

\bibliographystyle{abbrv}

\vspace{0.7cm}
\noindent 
Unit\'e de math\'ematiques pures et appliqu\'ees de l'\'Ecole normale sup\'erieure de Lyon,\\
 CNRS - Universit\'e de Lyon.

\end{document}